\documentclass[onecolumn]{IEEEtran}
\usepackage[dvips]{graphicx}
\usepackage{times}
\usepackage{fancybox,latexsym,amsmath,psfrag}
\usepackage{epsf,epsfig}
\usepackage{amssymb}
\usepackage{graphicx}
\usepackage{labelfig}
\usepackage{graphics,psfig,times}
\usepackage{rotate,float,tabularx,supertabular,hhline,multirow}
\usepackage{times, multicol}
\usepackage{fancyhdr}
\def\BibTeX{{\rm B\kern-.05em{\sc i\kern-.025em b}\kern-.08em
    T\kern-.1667em\lower.7ex\hbox{E}\kern-.125emX}}

\DeclareMathOperator*{\argmin}{arg\,min}

\begin{document}

\renewcommand{\headrulewidth}{0pt}

\title{Non-Myopic Sensor Control for Target Search and Track Using a Sample-Based GOSPA  Implementation}

 \author{
 $\ $ \\
 \begin{center}
 Marcel Hernandez, \'Angel F. Garc\'ia-Fern\'andez, and Simon Maskell
 \end{center}
%
%\newline\newline
$\ $ \\
School of Electrical Engineering, Electronics and Computer Science,
\\ The University of Liverpool,
Brownlow Hill,
Liverpool, L69 3GJ}
%Email: Marcel.Hernandez@liverpool.ac.uk, Angel.Garcia-Fernandez@liverpool.ac.uk, smaskell@liverpool.ac.uk}

\maketitle

%\vspace{0.5cm}
% \centerline{} \centerline{Revised manuscript submitted to \textit{IEEE Transactions on Aerospace and Electronic Systems}, 24 July, 2023.}
 %  \centerline{} \centerline{} \centerline{}

\begin{abstract}
This paper  is concerned with sensor management for target search and track using the
generalised optimal subpattern assignment (GOSPA) metric. Utilising the GOSPA metric to
predict future system performance is
computationally challenging, because of the need to account for uncertainties within the
scenario, notably the number of targets, the locations of targets, and the measurements generated by the targets subsequent to performing sensing actions.
In this paper, efficient sample-based techniques are developed to calculate the predicted mean square GOSPA metric. These techniques allow for missed detections and false alarms, and thereby enable the metric to be exploited in scenarios
more complex than those previously considered. Furthermore, the GOSPA methodology is extended to perform non-myopic (i.e. multi-step) sensor management via the development of a Bellman-type recursion that
optimises a conditional GOSPA-based metric.
Simulations for scenarios with missed detections, false alarms, and planning horizons of up to three time steps demonstrate the approach, in particular showing that optimal plans align with an intuitive understanding of how taking into
account the opportunity to make future observations should influence the current action.
It is concluded that the GOSPA-based, non-myopic search and track algorithm offers a powerful mechanism for sensor
management.
\end{abstract}

\vspace{0.25cm}

%\begin{keywords}
%\noindent

\section{Introduction}

Recently, there has been great interest in sensor management \cite{Krishnamurthy_book16} for  search and track of multiple targets (e.g. see
\cite{bostrom_rost_2021,bostrom_rost_2022} and references therein). In \cite{bostrom_rost_2021,bostrom_rost_2022}, non-myopic (i.e.
multi-step) planning was performed using a Poisson multi-Bernoulli mixture (PMBM) filter \cite{williams_2015,Garcia-Fernandez_2018} as
the basis for search, and the  ``predicted ideal measurement set approach'' \cite{mahler_2004}\footnote{The ideal measurement set approach
assumes that the measurements are error-free and there are no missed detections or false alarms.\vspace{0.2cm}}
to
predict the benefit of observing/updating a target track. The optimal control problem then used a cost function
that was a weighted sum of the integrated intensity of the PMBM density (for search) and the sum total of the trace
of each updated target covariance (for tracking).
In other recent work (see \cite{Charrow_2015,van_nguyen_2020}), non-myopic planning approaches were developed that
assumed the origin of measurements was known, allowing each
detected target to be tracked  independently. An  occupancy grid
filter was then used to represent the presence of undetected targets.

In \cite{tharman_2007}, a sensor management approach was developed for scenarios with an unknown number of targets
that allowed for noisy measurements, with measurement origin uncertainty, and potentially both missed detections
and false alarms. The basis of the approach was to use the posterior Cram\'er-Rao bound (PCRB)
\cite{van_trees_2013} to predict multi-target tracking performance, accounting for measurement origin uncertainty
via a matrix of ``information reduction factors'' (e.g. see \cite{zhang_2001}). In conducting search to detect new
targets, ``particles'' (i.e. target hypotheses) were distributed uniformly along the perimeter of the surveillance
region, and the probability of detecting a new target was estimated if the sensor observed a region containing one
or more hypothesis. The bi-criterion optimisation problem then used a cost function that was a weighted total of
the multi-target PCRB and the probability of detecting a new target. It is noted that the   approach developed in
\cite{tharman_2007} built on earlier  sensor management research for target tracking that used the PCRB as the
objective function\footnote{The PCRB was shown to offer an accurate mechanism for predicting tracking performance,
allowing sensor management to be performed in time-critical applications.\vspace{0.2cm}} \cite{hernandez_2004} (see also
\cite{hernandez_2013}).

\newpage
\noindent
The generalised optimal subpattern assignment (GOSPA) metric \cite{Rahmathullah17} provides a mechanism for
combining together the costs corresponding to localisation errors for properly
detected targets, and errors for missed and false targets. These three types of errors are of major interest in multiple target estimation. Importantly, the GOSPA metric is not prone to two erroneous effects that
hinder its predecessor,  the optimal subpattern assignment metric (OSPA)  \cite{Schuhmacher08b,schuhmacher_2008}, specifically:
(i): the OSPA metric does not necessarily increase as the number of false targets increases \cite{Rahmathullah17};
and: (ii):
the OSPA metric is prone to the ``spooky effect'' in optimal metric-based estimation \cite{garcia_fernandez_2019}.

In order to use the GOSPA metric as the basis for \textit{predictive} sensor management, it is necessary to
account for uncertainties within the scenario, namely the number of targets, the locations of targets,
and the measurements generated by the targets as a result of performing sensing actions. To this end, the GOSPA
metric is averaged over these uncertainties to  provide the average minimum mean squared GOSPA (AMMS-GOSPA) error \cite{garcia_fernandez_2021}.

Calculating the AMMS-GOSPA is computationally challenging, because of the required averaging across the
uncertainties, and the minimisation (within the GOSPA metric).  Consequently, in
\cite{garcia_fernandez_2021}, sensor management was demonstrated for relatively simple
scenarios with just a single action (i.e. myopic planning), and a single target state hypothesis,
(represented by a Dirac delta function) per target. Either one target or two well separated targets were
considered\footnote{The separation of the targets allowed the GOSPA metric to be summed across the two
targets.\vspace{0.2cm}}. The AMMS-GOSPA was calculated analytically, and the resulting sensor management strategies
demonstrated the tradeoff between sensing costs and the probability of target existence, with (all else being
equal) the target(s) more likely to be observed as the probability of existence increased or the sensing costs
decreased. Also of note, in \cite[Chapter~6]{Ubeda-Medina_2018}, myopic GOSPA-based sensor management was performed using a Bernoulli-Gaussian approximation of the conditional squared GOSPA error. However, the algorithm developed in \cite[Chapter~6]{Ubeda-Medina_2018} is for track-before-detect applications, and is therefore not suitable for the detection-based measurement model that is the focus of the current paper.

Building on this previous research, this article provides details of the following theoretical advancements:
 \begin{enumerate}
 \vspace{0.15cm}
\item  Analytical calculation of the MS-GOSPA for different combinations of measurements and actions.
\vspace{0.15cm}
 \item The development of efficient sampling techniques (i.e. of the measurements), to reduce the computational
     complexity of calculating the  AMMS-GOSPA error.
     \vspace{0.15cm}
 %\newpage
 \item The development of an optimal non-myopic planning (Bellman type, e.g. \cite{Bellman_1952}) recursion that exploits the conditional AMMS-GOSPA error.
 \end{enumerate}

%\noindent
%As a baseline test, the efficient sampling techniques were applied to the single target scenario from
%\cite{garcia_fernandez_2021} and generated identical results.

\vspace{0.15cm}
\noindent
Efficient calculation of the AMMS-GOSPA  allows the GOSPA metric  to be exploited in scenarios significantly
more complex than those considered in \cite{garcia_fernandez_2019,garcia_fernandez_2021}. To this end, the approach is
demonstrated in performing sensor management with:
(i): a high degree of uncertainty in each target location\footnote{E.g. representing a search track, or a
target that has not yet been accurately geo-located.\vspace{0.2cm}}, with the prior distribution represented by a mixture of
weighted Dirac delta functions, e.g. as in a particle filter estimate \cite{arulampalam_2002}; and (ii): time
horizons with multiple time steps (i.e. non-myopic planning).

The remainder of this paper is as follows.
In Section \ref{sec:performance_metrics}, a review of the GOSPA metric is provided.
In Section \ref{sec:exploit}, details are provided of the existing state-of-the-art regarding utilising  the GOSPA metric  for predictive sensor management via optimising the AMMS-GOSPA error.
In Section \ref{sec:multi_time_error}, the AMMS-GOSPA is generalised to multiple time step scenarios, and analytical
equations for the MMS-GOSPA error are calculated.
In Section \ref{sec:AMMS_GOSPA}, efficient sample-based approximations of the
AMMS-GOSPA error are presented, along with baseline tests applying the approaches to  a scenario from
\cite{garcia_fernandez_2021} in which the optimal solution was determined analytically. It is shown that an efficient sampling approach determines the optimal solution with a circa 250 times reduction in computational expensive compared to a  sampling approach that does not take into account the possible measurement sequences.
In Section \ref{sec_non_myopic}, the optimal non-myopic planning approach is presented, along with a suboptimal approach and a baseline approach that minimises target localisation errors.
In Section \ref{sec_sim}, simulation results demonstrate the multi-step planning approaches, and in particular show
that the optimal actions align with an intuitive understanding of how taking into account the ability to make
further observations should influence the current observation.
In Section \ref{sec:conc}, a Summary and Conclusions are provided. Finally, in an Appendix, it is proven that when the ideal measurement set approximation is  used (but allowing $P_d<1$), the optimal non-myopic planning approach generates identical solutions to a commonly implemented  suboptimal approach.
%

%\newpage

\section{GOSPA Metric}
\label{sec:performance_metrics}
This section reviews the GOSPA metric. We first present the notation and then its definition.

\newpage

\subsection{Notation}

Let $c$ and $p$ be two real numbers such that $c>0$ and $1\leq p<\infty$.
Let $d\left(\cdot,\cdot\right)$ denote a metric on the single
target space\footnote{In this paper, the metric $d\left(\cdot,\cdot\right)$ will denote the geo-location distance
(error) between the target (typically denoted by $X$) and the target state estimate (typically denoted by
$\hat{X}$).}, which is typically ${\mathbb R}^{n_{x}}$.

Let $X=\left\{ x_{1},...,x_{\left|X\right|}\right\} $ and $Y=\left\{ y_{1},...,y_{\left|Y\right|}\right\} $
denote two finite sets of targets, with $\left|X\right|\leq\left|Y\right|$,
and $\left|X\right|$ being the cardinality (number of elements) of the
set $X$.
In the context of target tracking, $X$ typically represents the set of target ground-truth states, and $Y$ the set
of target state estimates.

Let $\gamma$ be an assignment set between $\left\{ 1,...,\left|X\right|\right\} $
and $\left\{ 1,...,\left|Y\right|\right\} $, which satisfies $\gamma\subseteq\left\{ 1,...,\left|X\right|\right\}
\times\left\{ 1,...,\left|Y\right|\right\} $,
$\left(i,j\right),\left(i,j'\right)\in\gamma\rightarrow j=j'$, and
$\left(i,j\right),\left(i',j\right)\in\gamma\rightarrow i=i'$. The
last two properties ensure that every $i$ and $j$ have at most one
assignment. The set of all possible $\gamma$ is denoted by $\Gamma$.

\subsection{Metric}

The GOSPA metric, with parameters $p$ and $c$, between $X$
and $Y$ (for  $\alpha=2$) is given as follows (see \cite[Proposition 1]{Rahmathullah17}):
\vspace{-0.2cm}
\begin{eqnarray}
 d_{p}^{\left(c,2\right)}\left(X,Y\right) &
  \triangleq &\min_{\gamma\in\Gamma}\Bigg(\sum_{\left(i,j\right)\in\gamma} \label{eq:GOSPA_alpha2}  \big[d\left(x_{i},y_{j}\right)\big]^{p}
  +\frac{c^{p}}{2}\big(\left|X\right|-\left|\gamma\right|+\left|Y\right|-\left|\gamma\right|\big)\Bigg)^{1/p}
\end{eqnarray}

\noindent
The first term in (\ref{eq:GOSPA_alpha2}) represents the localisation
errors (to the $p$-th power) for assigned targets (properly detected
ones), which meet $\left(i,j\right)\in\gamma$. The terms
$\frac{c^{p}}{2}\left(\left|X\right|-\left|\gamma\right|\right)$
and $\frac{c^{p}}{2}\left(\left|Y\right|-\left|\gamma\right|\right)$
represent the costs (to the $p$-th power) for missed and false targets respectively.

Compared to the OSPA metric, the GOSPA metric has an additional
parameter $\alpha$ that controls the cardinality mismatch penalty.
Importantly, as shown in equation (\ref{eq:GOSPA_alpha2}), only for $\alpha=2$ can the GOSPA metric be written
in terms of costs corresponding to localisation errors for properly
detected targets, missed and false targets, which are usually the
penalties of interest in multiple target estimation. Consequently, $\alpha=2$ is used throughout this paper.

\section{Exploitation of the  GOSPA Metric for Myopic Sensor Management}
\label{sec:exploit}

In order to use the GOSPA metric for predictive sensor management it is necessary to account for potential
target state  and measurement origin/accuracy uncertainties. To this end, the following errors are determined (see \cite{garcia_fernandez_2021}).

\subsection{Mean Squared GOSPA (MS-GOSPA) Error}

Given an action $a\in \mathbb{A}$ (e.g. sensor mode or steer direction), and a resulting measurement $z$, a
posterior multi-target state estimate $\hat{X}(z,a)$ can be determined. The resulting mean squared GOSPA (MS-GOSPA)
error given $z$ and $a$ is then defined as follows:
\begin{eqnarray}
  \mbox{MS-GOSPA}(\hat{X};z,a)
 & \triangleq &
 {\mathbb E}_X\bigg[\Big(d_{2}^{\left(c,2\right)}\big(X,\hat{X}(z,a)\big)\Big)^{2}\bigg|z;a\bigg] \\
 &  =  & \int_X\Big(d_{2}^{\left(c,2\right)}\big(X,\hat{X}(z,a)\big)\Big)^{2}p\left(X|z;a\right)d X\qquad
 \label{eq:posterior_MSGOSPA}
\end{eqnarray}

\subsection{Minimum MS-GOSPA (MMS-GOSPA) Error}

The minimum MS-GOSPA (MMS-GOSPA) error is achieved by selecting the estimate $\hat{X}(z,a)$ to minimise equation
(\ref{eq:posterior_MSGOSPA}), i.e.:\vspace{-0.2cm}
\begin{eqnarray}
  \mbox{MMS-GOSPA}(z,a) \label{eq:posterior_MMSGOSPA}
 & = & \min_{\hat{X}(z,a)}
 \int_X\Big(d_{2}^{\left(c,2\right)}\big(X,\hat{X}(z,a)\big)\Big)^{2}p\left(X|z;a\right)d X\qquad
\end{eqnarray}

\subsection{Average MMS-GOSPA (AMMS-GOSPA) Error}

The average MMS-GOSPA (AMMS-GOSPA) for action $a$, averaged over the measurement $z$ is then given as follows:
\begin{eqnarray}
  \mbox{AMMS-GOSPA}(a)
&  \triangleq  &  {\mathbb E}_z\big[\mbox{MMS-GOSPA}(z,a)  \big] \\
 &  = & \int_z
 \mbox{MMS-GOSPA}(z,a) p(z;a) dz
\qquad \label{eq:posterior_AMMSGOSPA}
\end{eqnarray}

 \noindent
For each pair $(z,a)$, the AMMS-GOSPA selects the estimate $\hat{X}(z,a)$ that minimises the MS-GOSPA metric, and
then averages over the value of $z$. The AMMS-GOSPA therefore gives the average minimum MS-GOSPA, taking into
account both uncertainty in the target state $X$, and uncertainty in the measurement $z$ for each potential
action $a$.

 \subsection{Optimal Action}

 Selecting the action $a$ that minimises the sum total of the AMMS-GOSPA error and sensing costs (if applicable) provides a
 mechanism for performing (myopic) sensor management, enabling the system to balance the objectives of (i):
 minimising localisation errors for properly
detected targets, and: (ii): minimising cardinality errors resulting from missed and false targets, (iii):
minimising sensing costs.

\subsection{Computational Complexity}

Clearly, calculating the AMMS-GOSPA is computationally challenging, because of the averaging performed (over $X$
given $(z,a)$, and then $z$ given $a$), and the minimisations (within the GOSPA metric, and then over all potential
posterior estimates).  Consequently, e.g. in \cite{garcia_fernandez_2021}, sensor management was demonstrated for a
relatively simple scenario with just a single action (i.e. myopic planning), and a single potential target location
hypothesis.

In this paper,  efficient techniques are developed to calculate the AMMS-GOSPA, enabling the metric to be used as a
basis for performing sensor management in more complex scenarios (i.e. with targets whose states are highly
uncertain, and multiple time steps).

\section{Multiple Time Step AMMS-GOSPA Error}
\label{sec:multi_time_error}

In this section the AMMS-GOSPA error is determined for a time window containing multiple time steps.
%Analytical expressions are determined for the MMS-GOSPA error on each time step.

\subsection{Decomposition Over Time Steps}

To extend the methodology to calculate the multiple time step AMMS-GOSPA error, some additional notation is required. Let
$T\geq 1$ denote the number of time steps considered (e.g. $T=1$ denotes myopic planning). Let:
\begin{eqnarray}
\hat{X}_{1:T} & \triangleq & \big(\hat{X}_1, \ldots, \hat{X}_T\big)
\end{eqnarray}
where $\hat{X}_{1:T}$ is a sequence of sets of target state estimates at times 1 through $T$. Similarly, let
$a_{1:T}\triangleq\left(a_{1},\ldots,a_{T}\right)$ denote the sensor actions at times 1 through $T$; and
let $z_{1:T}\triangleq\left(z_{1},\ldots,z_{T}\right)$ denote the sensor measurements at times 1 through $T$.

 The MS-GOSPA error at time $t$ is then given as follows:
\begin{eqnarray}
 \mbox{MS-GOSPA}(\hat{X}_t;z_{1:t}, a_{1:t})
\label{new_MS_GOSPA}
&= & \int_{X_{t}} \left(d_2^{(c,2)}(X_t,\hat{X}_t(z_{1:t},
a_{1:t}))\right)^2
%\qquad \nonumber \\
%&& \qquad \times
p(X_{t}| z_{1:t}, a_{1:t})dX_{t}\quad
\end{eqnarray}

\noindent
The multiple step  MS-GOSPA error is the discounted sum of MS-GOSPA errors at times 1 through $T$, i.e.:
\begin{eqnarray}
 \mbox{MS-GOSPA}_{1:T}(\hat{X}_{1:T};z_{1:T}, a_{1:T})
& =  & \sum_{t=1}^T\lambda^{t-1} \mbox{MS-GOSPA}(\hat{X}_t;z_{1:t}, a_{1:t})
\end{eqnarray}
where the  parameter $\lambda\in[0,1]$ is the discount factor.

The multiple step  MMS-GOSPA is then defined as follows:
\begin{eqnarray}
\mbox{MMS-GOSPA}_{1:T}(z_{1:T},a_{1:T})
&\triangleq & \min_{\hat{X}_{1:T}}\Big[\sum_{t=1}^T
\lambda^{t-1}\mbox{MS-GOSPA}(\hat{X}_{t}; z_{1:t}, a_{1:t})\Big]  \\
& = & \sum_{t=1}^T\lambda^{t-1}\min_{\hat{X}_{t}}\Big[ \mbox{MS-GOSPA}(\hat{X}_{t}; z_{1:t}, a_{1:t})\Big]
\label{new_GOSPA}\\
&  = & \sum_{t=1}^T\lambda^{t-1}\mbox{MMS-GOSPA}(z_{1:t}, a_{1:t})\label{new_GOSPA_multi}
\end{eqnarray}
where:
\begin{eqnarray}
 \mbox{MMS-GOSPA}(z_{1:t}, a_{1:t})
& \triangleq  &\min_{\hat{X}_{t}}\Big[ \mbox{MS-GOSPA}(\hat{X}_{t}; z_{1:t},
a_{1:t})\Big]
\end{eqnarray}
The simplification in (\ref{new_GOSPA}) is because given $z_{1:t}$ and  $a_{1:t}$, the MS-GOSPA at each
     time step $t$  is dependent only on the target state estimate $\hat{X}_{t}$ at that time, and critically
     is independent of $\hat{X}_{1:t-1}$, and also independent of $\hat{X}_{t+1:T}$ (i.e. estimates that perform smoothing are not considered).

The AMMS-GOSPA error at times 1 through  $T$ is then given as
follows:
\begin{eqnarray}
 \mbox{AMMS-GOSPA}_{1:T}(a_{1:T})
&\triangleq & \int_{z_{1:T}} \mbox{MMS-GOSPA}_{1:T}(z_{1:T},a_{1:T})
p(z_{1:T};a_{1:T}) dz_{1:T}  \\
& = & \sum_{t=1}^T\lambda^{t-1}\mbox{AMMS-GOSPA}(a_{1:t}) \label{simp_multi}
\end{eqnarray}
where:
\begin{eqnarray}
 \mbox{AMMS-GOSPA}(a_{1:t})
 & \triangleq  & \int_{z_{1:t}}\mbox{MMS-GOSPA}(z_{1:t}, a_{1:t})  p(z_{1:t};a_{1:t}) dz_{1:t}
\end{eqnarray}

\noindent
It is noted that equation (\ref{simp_multi}) does not rely on any assumptions, e.g. regarding the number of targets or the number of sensors.

\subsection{Target Hypotheses}

For the remainder of this article, it is assumed that
there is either zero or one target in the focal scenario\footnote{The generalisation to scenarios to multiple
targets is straightforward, provided that the targets are well separated. Indeed, in \cite{garcia_fernandez_2021},
it was shown that the AMMS-GOSPA was additive across the Bernoulli components for well separated targets.
%Indeed,
%in a future paper on multi-target search and track, it will be shown that provided the separation between targets
%is greater than the cardinality error, this result is true.
In scenarios with targets that are in close proximity
to each other, the analytical results presented later in Section \ref{sec:analytical_multi} no longer hold, and the AMMS-GOSPA
error must also consider the optimal assignment between the hypothesised targets and the target state estimates.\vspace{0.25cm}}.
Let the set of potential target hypotheses be denoted by:
\begin{eqnarray}
\mbox{target hypotheses, } X = \left\{x_0, x_1, \ldots, x_n\right\}
\end{eqnarray}

\noindent
where ``$x_0=\phi$'' denotes the hypothesis that a target is not present, and $x_i$, $i=1, \ldots,n$ denotes the
hypothesis that a target exists and its state is given by $x_i$. The prior probability of $X$ is given as follows:
\begin{eqnarray}
f(X) &=& \left\{
\begin{array}{ll}
rw_i\delta(X-x_i) & \quad \mbox{for }   i=1,\ldots,n\qquad\qquad \\
1-r  & \quad \mbox{otherwise}
\end{array}
\right.
\end{eqnarray}
where $r$ is the probability of target existence, and $w_i$ denotes the weight of each location hypothesis, with
$\displaystyle \sum_{i=1}^n w_i = 1$.

Note that $f(X)$ is a Bernoulli density with probability of existence $r$ and whose single target density is a
mixture of weighted Dirac delta functions (e.g. as in a particle filter implementation \cite{arulampalam_2002}). The probability of target existence remains fixed throughout (i.e. it is assumed that the  probability of survival is unity).
It is also noted that $x_i$ ($i\geq 1$) will generally depend  on $t$ (unless the target is stationary). For brevity, this dependency is omitted from the notation.

In applications in which a particle filter is not used to estimate the target state(s), the hypotheses represent samples drawn from the target probability density functions (in tracking applications), or samples drawn from a region under surveillance (in search applications).

\subsection{Calculation of the MS-GOSPA Error at Each Time Step}

 The MS-GOSPA error at time $t$ (given by equation (\ref{new_MS_GOSPA})) can be written in terms of the target hypotheses as follows:
\begin{eqnarray}
 \mbox{MS-GOSPA}(\hat{X}_t;z_{1:t}, a_{1:t}) \label{new_GOSPA_2}
& = & \sum_{i=0}^n \left(d_2^{(c,2)}(x_i,\hat{X}_t(z_{1:t}, a_{1:t}))\right)^2 p(x_i| z_{1:t},
a_{1:t})
\end{eqnarray}

\noindent
The probability $p(x_i| z_{1:t}, a_{1:t})$  can be written as follows (for $i\geq 0$):
\begin{eqnarray}
p(x_i| z_{1:t}, a_{1:t}) & = & \frac{p(z_{1:t}|x_i,a_{1:t}) p(x_i)}{\gamma_t} \\
&=& \frac{p(x_i)}{\gamma_t}\prod_{j=1}^t p(z_j|x_i,a_j) \qquad \label{prob}
\end{eqnarray}

%\newpage

\noindent
where $p(x_i)$ denotes the prior probability of hypothesis $i$ being true, i.e.:
\begin{eqnarray}
p(x_i) &=& \left\{
\begin{array}{ll}
rw_i & \quad \mbox{for }i>0 \qquad \\
(1-r) & \quad \mbox{otherwise}
\end{array}
\right.\label{px_i}
\end{eqnarray}

%\newpage
\noindent
and $\gamma_t$ is a normalising constant that ensures that $\displaystyle \sum_{i=0}^n p(x_i| z_{1:t}, a_{1:t})=1$.

\newpage

\subsection{Analytical Calculation of the MMS-GOSPA Error}

\label{sec:analytical_multi}

Based on the measurements $z_{1:t}$ and action sequence $a_{1:t}$, the following two target state estimates
$\hat{X}_t$ are considered:
\begin{eqnarray}
\hat{X}_t &=& \phi \qquad \mbox{(i.e. no target is present)}\qquad \\
\hat{X}_t &=&  \sum_{i=1}^n w_{1:t}^i x_i  \ \triangleq \ \hat{X}_e \label{estimate}
\end{eqnarray}
where $w_{1:t}^i$ are the posterior hypothesis weights based on the measurements $z_{1:t}$ and action sequence
$a_{1:t}$, given as follows (for $i \geq 1$):
\begin{eqnarray}
 w_{1:t}^i & \propto & p(x_i| z_{1:t}, a_{1:t}) \\
 &=&\frac{p(z_{1:t}|x_i,a_{1:t}) p(x_i)}{\lambda_t} \\
&=& \frac{p(x_i)}{\lambda_t}\prod_{j=1}^t p(z_j|x_i,a_j) \qquad \label{new_weights}
\end{eqnarray}

\noindent
and $\lambda_t$ is a normalising constant that ensures that $\displaystyle \sum_{i=1}^n w_{1:t}^i =1$.
%
%
%Note that $\lambda_t$ is smaller than $\gamma_t$ due to the %fact that the summation is for $i\geq 1$ and not $i\geq 0 $.
%
It then follows from equations (\ref{prob}) and (\ref{new_weights}) that:
\begin{eqnarray}
 p(x_i| z_{1:t}, a_{1:t})  &=& \lambda_t w_{1:t}^i / \gamma_t \qquad \label{new_lambda}
\end{eqnarray}

\noindent
Summing equation  (\ref{new_lambda}) over $i=1,\ldots,n$, it follows that $\lambda_t/\gamma_t$ gives the posterior probability that a target is not present.

To proceed, it is straightforward to show the following:
\begin{eqnarray}
&& d_2^{(c,2)}(x_0,\phi) \ = \ 0 \label{analytical1} \qquad
  \mbox{(i.e. no target is present or estimated to be present)}  \\
&& d_2^{(c,2)}(x_i,\phi) \ = \ c/\sqrt{2}\qquad \mbox{for }i>0 \label{analytical2} \qquad
  \mbox{(i.e. the target is missed, cardinality error
$=$ 1)} \qquad\\
&& d_2^{(c,2)}(x_0,\hat{X}_e) \ = \ c/\sqrt{2} \label{analytical3} \qquad
 \mbox{(i.e. a target is incorrectly estimated to be present, cardinality error
$=$ 1)}\qquad
 \\
&& d_2^{(c,2)}(x_i,\hat{X}_e) \ = \ \min\big(d(x_i,\hat{X}_e), c\big) \label{cold_1}
\end{eqnarray}
where $d(x_i,\hat{X}_e)$ is the distance between the hypothesis $x_i$ and the estimate $\hat{X}_e$, and $c$ is the cut-off. It is noted that
$d_2^{(c,2)}()$ is actually a distance between sets, so one could write $\{x_i\}$ instead of $x_i$, and  $\{\hat{X}_e\}$ instead of $\hat{X}_e$ etc. The simplified notation is used for brevity.

The square of the truncated distance error in equation (\ref{cold_1}) can be written as follows:
\begin{eqnarray}
\left[d_2^{(c,2)}(x_i,\hat{X}_e)\right]^2 &=& i_d(x_i,\hat{X}_e) d(x_i,\hat{X}_e)^2 \label{analytical3_new}
+ \big(1 -
i_d(x_i,\hat{X}_e)\big)c^2
\end{eqnarray}
where $i_d(x_i,\hat{X}_e) = 1$ if $d(x_i,\hat{X}_e)\leq c$, and $i_d(x_i,\hat{X}_e) = 0$ otherwise.

It then follows from equations (\ref{new_GOSPA_2}) and (\ref{analytical1}) -- (\ref{analytical2})  that:
\begin{eqnarray}
 \mbox{MS-GOSPA}(\hat{X}_t=\phi;z_{1:t}, a_{1:t})
& = & \frac{c^2}{2} \sum_{i=1}^n  p(x_i| z_{1:t}, a_{1:t})
\label{multi_MSGOSPA_1} \\
&  =  & \frac{c^2}{2} \big(1-  p(x_0| z_{1:t}, a_{1:t})\big)
\quad\label{multi_MSGOSPA_1b}
\end{eqnarray}
Equation (\ref{multi_MSGOSPA_1b}) gives the squared cost of a cardinality error, weighted  by the posterior
probability that a target is present, based on the measurements $z_{1:t}$ and action sequence $a_{1:t}$.

To continue, from equations (\ref{new_GOSPA_2}) and (\ref{analytical3}) -- (\ref{analytical3_new})  it follows
that:
\begin{eqnarray}
 \mbox{MS-GOSPA}(\hat{X}_t=\hat{X}_e;z_{1:t}, a_{1:t})
&= & \frac{c^2}{2}p(x_0| z_{1:t}, a_{1:t}) + \sum_{i=1}^n \Big[i_d(x_i, \hat{X}_e) d(x_i, \hat{X}_e)^2 \qquad \\
&& \qquad + \big(1-
i_d(x_i, \hat{X}_e)c^2\big)\Big]p(x_i| z_{1:t}, a_{1:t})\qquad\qquad \nonumber
\end{eqnarray}

\noindent
This can be written as follows:
%\begin{eqnarray}
%&& \mbox{MS-GOSPA}(\hat{X}_t=\hat{X}_e;z_{1:t}, %a_{1:t})\label{multi_MSGOSPA_2} \\
%
%&& =  \frac{c^2}{2}p(x_0| z_{1:t}, a_{1:t})  \nonumber \\
%&&
%\qquad + \frac{\lambda_t}{\gamma_t} %\Bigg(\sum_{l=1}^2\hat{{\mathcal V}}_l(x) %+c^2T_{{w}_{1:t}}(\hat{X}_e) \Bigg)\qquad\qquad
%nonumber
%\end{eqnarray}
%
\begin{eqnarray}
 \mbox{MS-GOSPA}(\hat{X}_t=\hat{X}_e;z_{1:t}, a_{1:t}) \label{multi_MSGOSPA_1c}
\ = \ \frac{c^2}{2}p(x_0| z_{1:t}, a_{1:t})  + \big(1-p(x_0|z_{1:t},a_{1:t})\big)
\Bigg(\sum_{l=1}^2\hat{{\mathcal V}}_l(x)
+c^2T_{{w}_{1:t}}(\hat{X}_e) \Bigg)
\end{eqnarray}

\noindent
where:
\begin{eqnarray}
 \hat{{\mathcal V}}_l(x) &= &{\mathcal E}_{\underline{{w}}_{1:t}}\left[x(l)^2\right]  \label{variance_new} +
{\mathcal E}_{{w}_{1:t}}\left[x(l)\right]^2 -
2{\mathcal E}_{\underline{{{w}}}_{1:t}}\left[x(l)\right]{\mathcal E}_{{w}_{1:t}}\left[x(l)\right]   \\
\underline{{w}}^i_{1:t} &=&  {w}_{1:t}^i i_d(x_i,\hat{X}_e) \quad \mbox{for }i=1, \ldots, n\qquad \\
{\mathcal E}_{\underline{{{w}}}_{1:t}}\left[x(l)^m\right]
& = & \sum_{i=1}^n \underline{{w}}^i_{1:t} x_i(l)^m \quad \mbox{for }m=1,2\\
%}
%
T_{{w}_{1:t}}(\hat{X}_e) &=& \left(1-\sum_{i=1}^n i_d(x_i, \hat{X}_e){w}_{1:t}^i\right) \label{T_value}
\end{eqnarray}
where $\hat{{\mathcal V}}_l(x)$ denotes the cross-covariance of the (posterior) weighted samples of the $l$-th dimension (denoted $x(l)$)
of the target state (herein $x(1)$ and $x(2)$ denote the $x$- and $y$- coordinates of the target state respectively), based on the measurements $z_{1:t}$ and action sequence $a_{1:t}$. $T_{{w}_{1:t}}(\hat{X}_e)$
denotes the sum of weighted samples for which the estimation error exceeds $c$, given that a target exists.

The first term on the right-hand side of equation (\ref{multi_MSGOSPA_1c}) is the squared cost of a cardinality
error, multiplied  by the posterior probability that a target is not present. The second term is the average
estimation error squared (adjusted for a maximum estimation error squared of $c^2$) of the posterior estimate
$\hat{X}_e$, multiplied  by the posterior probability that a target is present.

The MMS-GOSPA$(z_{1:t}, a_{1:t})$ is then given by the minimum of   (\ref{multi_MSGOSPA_1b}) and
(\ref{multi_MSGOSPA_1c}).

\subsection{Target State Estimation in the Presence of Clutter}

In the presence of clutter, the expected likelihood particle filter (ELPF) approach is used \cite{marrs_2022}  to
determine the posterior hypothesis weights and target state estimates. These are used in the calculation of the performance metrics for the three (i.e. optimal, suboptimal, and baseline) control approaches.
The ELPF approach is analogous to the probabilistic data association Kalman filter (e.g. see \cite{barshalom_2001}).

If the target is within the field-of-view (FOV) of the sensor, it is detected with probability $P_d$. Conditional on target hypothesis $i>0$, each target generated measurement is sampled from a Gaussian distribution with mean $x_i$ and covariance $\Sigma$.
It is assumed that the number of false alarms at each sampling time has a Poisson distribution with mean $\lambda_{FA} V$, where $\lambda_{FA}$ is the false alarm rate per unit volume of the observation region, and $V$ is the volume of the FOV of the sensor. False alarms are uniformly distributed within the FOV of the sensor.

 Excusing an abuse of notation, let $Z=\left\{z_1,\ldots, z_N\right\}$ denote a vector of $N$ measurements (i.e. a maximum of one target generated measurement plus false alarms) at any given sampling time.
Using the ELPF approach, the measurement likelihood (used to calculate the posterior hypothesis weights, and the subsequent posterior target state estimate), for $N>0$, is given as follows \cite{marrs_2022}:
\begin{eqnarray}
p(Z|x_i, a) & \propto & \left\{
\begin{array}{ll}
\displaystyle \lambda_{FA} (1-P_d) + P_d\sum_{k=1}^N {\mathcal N}(z_k; x_i, \Sigma)  & \qquad \mbox{if }i>0 \mbox{ and }x_i \in \mbox{FOV}(a)
\\
\lambda_{FA} & \qquad \mbox{otherwise}
\end{array}
\right.
\end{eqnarray}
where ${\mathcal N}(z_k; x_i, \Sigma)$ denotes the Gaussian probability density function evaluated at $z_k$ for a distribution with mean $x_i$ and error covariance $\Sigma$. For $N=0$, $Pr(N=0|x_i,a) \propto (1-P_d)$ if $i>0$ and $ x_i \in \mbox{FOV}(a)$; and $Pr(N=0|x_i,a) \propto 1$ otherwise. The posterior hypothesis weights at each sampling time $t$ are then given by equation (\ref{new_weights}), from which the posterior target state estimate can be determined via equation (\ref{estimate}).

\section{Sample-Based Estimation of the AMMS-GOSPA Error}
\label{sec:AMMS_GOSPA}

\subsection{General Sampling Approach}
\label{general_sampling_approach}

The
 AMMS-GOSPA at time $t$ can be written as follows:
 \begin{eqnarray}
\mbox{AMMS-GOSPA}(a_{1:t})  & = & \int_{z_{1:t}} \mbox{MMS-GOSPA}(z_{1:t}, a_{1:t}) p(z_{1:t}; a_{1:t})
dz_{1:t}  \label{mulitcols_1}  \\
 & = & \int_X\int_{z_{1:t}} \mbox{MMS-GOSPA}(z_{1:t}, a_{1:t}) p(z_{1:t}| X, a_{1:t})p(X) dz_{1:t} dX \qquad \\
 & = & \sum_{i=0}^n p(x_i) \int_{z_{1:t}} \mbox{MMS-GOSPA}(z_{1:t}, a_{1:t}) p(z_{1:t}| x_i, a_{1:t})dz_{1:t}
 \label{125}
\end{eqnarray}
where $p(x_i)$ is given by equation (\ref{px_i}).

\newpage

\noindent
 A sample-based approximation is then given as follows:
 \begin{eqnarray}
 \mbox{AMMS-GOSPA}(a_{1:t}) \label{gen_AMMS_GOSPA}
 & \approx & \frac{1}{m}\sum_{i=0}^np(x_i)\sum_{l=1}^{m}\mbox{MMS-GOSPA}(z_{1:t}^l(x_i), a_{1:t}) \qquad
 \end{eqnarray}

\noindent
where $z_{1:t}^l(x_i)$ are $m$ measurement samples, conditional on the target hypothesis $x_i$.
In the most general case there can be multiple measurements per time step, which can include both  target generated measurements and false alarms.

 \subsection{Efficient Sampling Approach -- Conditioning on  the Measurement Sequence}

 \label{efficient_sampling_approach}

In this section, the sample-based AMMS-GOSPA error is simplified to take into account the sequence of target detections and missed detections.
The measurement vector at each time step can include a maximum of one target generated measurement plus false alarms. For notational brevity and ease of understanding, the derivations are presented for a
single sensor scenario. However, the approach can be readily applied  to scenarios in which there are multiple sensors.

 Let $s_{1:t}(m)=\left(s_1(m),\ldots, s_t(m)\right)$ denote a sequence of target detections and missed detections at times 1 through $t$, with $s_{k}(m)= 0 \mbox{ or } 1$.
The number of
potential detection sequences is $2^t$. $Pr\big(s_{1:t}(m)|x_i, a_{1:t}\big)$ is the probability of the sequence
$s_{1:t}(m)$ occurring, conditional on the target hypothesis $i$ and actions $a_{1:t}$. This is given as follows:
\begin{eqnarray}
Pr\big(s_{1:t}(m)|x_i, a_{1:t}\big) &=& \prod_{k=1}^t Pr(s_k(m)|x_i, a_k)\qquad
\end{eqnarray}
where:
\begin{eqnarray}
\begin{array}{rlll}
Pr(s_k(m)=0|x_i,a_k) & = & 1 & \quad\mbox{if }i=0, \mbox{ or } i>0 \mbox{ and } x_i\notin \mbox{FOV}(a_k)\\
Pr(s_k(m)=1|x_i,a_k) & = & 0 &  \quad\mbox{if }i=0, \mbox{ or } i>0 \mbox{ and } x_i\notin \mbox{FOV}(a_k) \\
Pr(s_k(m)|x_i,a_k)   & = & P_d^{s_k(m)}(1-P_d)^{1-s_k(m)}  & \quad \mbox{if }i > 0\mbox{ and } x_i\in \mbox{FOV}(a_k), \ \mbox{for }s_k(m)=\left\{0,1\right\}\label{prob_target_meas}
\end{array}
\end{eqnarray}
and $P_d$ is again the probability of detection on each sensor scan for a the target that is within the FOV of the sensor.
Equation (\ref{125}) can then be approximated as follows:
\begin{eqnarray}
 \mbox{AMMS-GOSPA}(a_{1:t})
 & \approx &   \frac{1}{n_h}\sum_{i=0}^n p(x_i) \sum_{m=1}^{2^t} Pr\big(s_{1:t}(m)|x_i, a_{1:t}\big)
 \sum_{j=1}^{n_h} \mbox{MMS-GOSPA}\big(z_{1:t}^{ij}(m), a_{1:t}\big) \label{new_det_approach_0} \qquad
 \end{eqnarray}
where $n_h$ is the number of measurement samples per target hypothesis and measurement sequence pair.

  The measurement sequences $z_{1:t}^{ij}(m)\triangleq\left(z_1^{ij}(m),\ldots,z_t^{ij}(m)\right)$, $i=0,\ldots,n$, $j=1,\ldots, n_h$ and  $m=1,\ldots,2^t$ can include both target generated measurements and false measurements, with a target generated measurement occurring at time  $k$ only if $s_k(m)=1$.

In the approximation (\ref{new_det_approach_0}), the computation is simplified by the fact that the summations need only be performed for detection sequences for which $Pr\big(s_{1:t}(m)|x_i, a_{1:t}\big)>0$. For example, for $i=0$ the only possible sequence has no target detection on any time step (because there is no target under this hypothesis).
Equally, for $i>0$, sequences for which $s_k=1$ when  $x_i \notin \mbox{FOV}(a_k)$ can also be discounted (i.e. a target cannot be detected if it is outside of the sensor FOV).

Therefore, the total number of samples of $z_{1:t}$ used in calculating the AMMS-GOSPA via equation (\ref{new_det_approach_0})
is $(2^tn +1)n_h$, and increases exponentially with the number of time steps $t$.

\subsection{Computational Complexity}

On first viewing, the approximation (\ref{new_det_approach_0}) would appear to be significantly more complex, and
require a greater number of samples, than the approximation (\ref{gen_AMMS_GOSPA}). However, the approximation
(\ref{new_det_approach_0}) accounts for all potential sequences of target detections and missed detections, whereas the
approximation (\ref{gen_AMMS_GOSPA}) does not. Consequently, the number of samples $m$ required in the
approximation (\ref{gen_AMMS_GOSPA}) is significantly greater than the number of samples $n_h$ in the approximation
(\ref{new_det_approach_0}). This is because the approximation (\ref{gen_AMMS_GOSPA}) must contain enough
samples to account for the target probability of detection.

For  example, if the measurements are extremely accurate, it is feasible to allow $n_h=1$ in
(\ref{new_det_approach_0}), whereas  $m>100$ is required for an accurate
estimate using (\ref{gen_AMMS_GOSPA}), even in the simplest case with $t=1$.
If measurements are less accurate, a large number of samples $n_h$ may be required irrespective of the
approximation used. This significantly  increases the  computational expense of the optimisation algorithm, which
is already computationally expensive for $t>1$.

\subsection{Demonstration}

This scenario is identical to analysis I in \cite{garcia_fernandez_2021}. There is one potential target location hypothesis $x_1$ (therefore, $n=1$), and
just a single time step (i.e. $t=1$).
If the sensor observes the potential location of the target, it generates a perfect measurement of the target state\footnote{In \cite{garcia_fernandez_2021}, the measurement probability density function  was specified via the Dirac delta function. Herein, the measurement has a Gaussian distribution with mean $x_1$ and error covariance $\Sigma = \mbox{diag}(10^{-10}, 10^{-10})$.\vspace{0.2cm}} with probability $P_d=0.6$. There are no false alarms. The cost of a cardinality error is $c=10$km. It is noted that the units of cardinality errors and sensing costs are both kilometres, in order to balance the units with the distance metric $d(x_i,\hat{X}_e)$, which is measured in kilometres.

There are two possibilities for the target state estimate, $\hat{X}(z,a)$:
\begin{eqnarray}
\begin{array}{llll}
\hat{X}(z,a) &=& \phi & \qquad \mbox{(no target is present)}  \label{est1} \\
\hat{X}(z,a) &=& x_1 & \qquad \mbox{(the potential target location)} \qquad \label{est2}
\end{array}
\end{eqnarray}
There are also two possible actions $a$, these being ``do not attempt an observation/measurement''  and ``observe the potential target location''.

The analytical solution determined in \cite{garcia_fernandez_2021} is compared to the solutions generated via the following three sample-based approaches:
\begin{enumerate}
\vspace{0.15cm}
\item \textbf{Approach 1:} The general sampling approach described in Section
\ref{general_sampling_approach}, with the MMS-GOSPA for each sample calculated via equation (\ref{eq:posterior_MMSGOSPA}).
\vspace{0.15cm}
\item \textbf{Approach 2:} The general sampling approach, exploiting the analytical calculations of the MMS-GOSPA provided in Section \ref{sec:analytical_multi}.
    \vspace{0.15cm}
\item \textbf{Approach 3:} The efficient sampling approach described in Section \ref{efficient_sampling_approach}, exploiting the analytical calculations of the MMS-GOSPA and the known probability of detection.
\end{enumerate}

\begin{figure}[H]
\vspace{-1.0cm}
\[
\begin{array}{ccc}%
%\SetLabels
%\L (0.12*0.825) \small{\textbf{(a):} $m=100$ samples} \\
%\endSetLabels
%\psfrag{y}[]{ sensing cost ($s$)}
%\psfrag{x}[]{ probability of existence ($r$)}
%\strut\AffixLabels{\includegraphics[width=5cm]{optimal_action_Pd_0.6_100.ps}} &
%
\SetLabels
\L (0.12*0.815) \textbf{(a):} $m=10$ samples \\
\endSetLabels
\psfrag{y}[]{ sensing cost ($s$)}
\psfrag{x}[]{ probability of existence ($r$)}\hspace{0.5cm}
\strut\AffixLabels{\includegraphics[width=6cm]{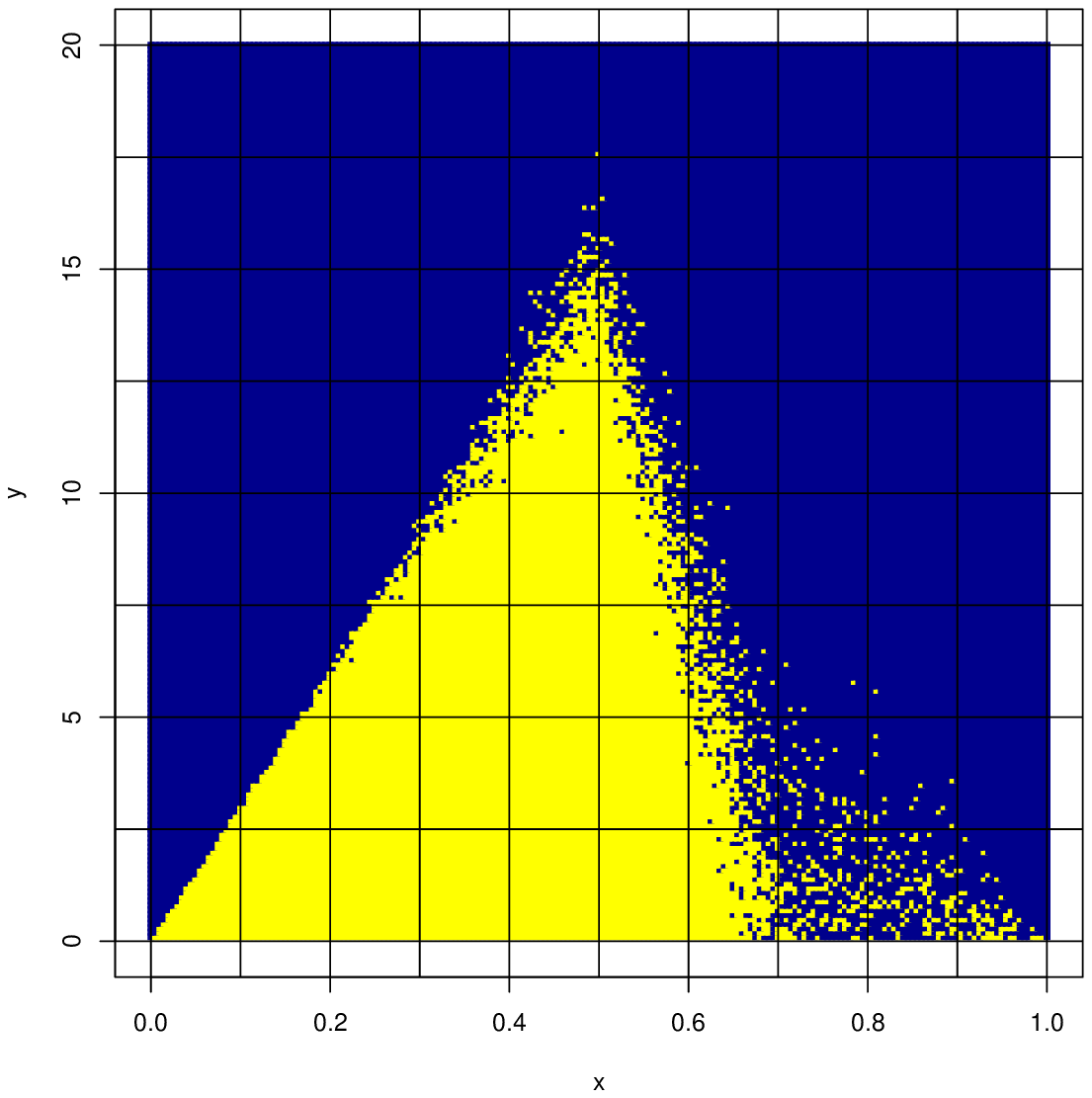}} &
\SetLabels
\L (0.12*0.815) \textbf{(b):} $m = 100$  samples \\
\endSetLabels
\psfrag{y}[]{ sensing cost ($s$)}
\psfrag{x}[]{ probability of existence ($r$)}
\vspace{-1.75cm}\hspace{0.5cm}
\strut\AffixLabels{\includegraphics[width=6cm]{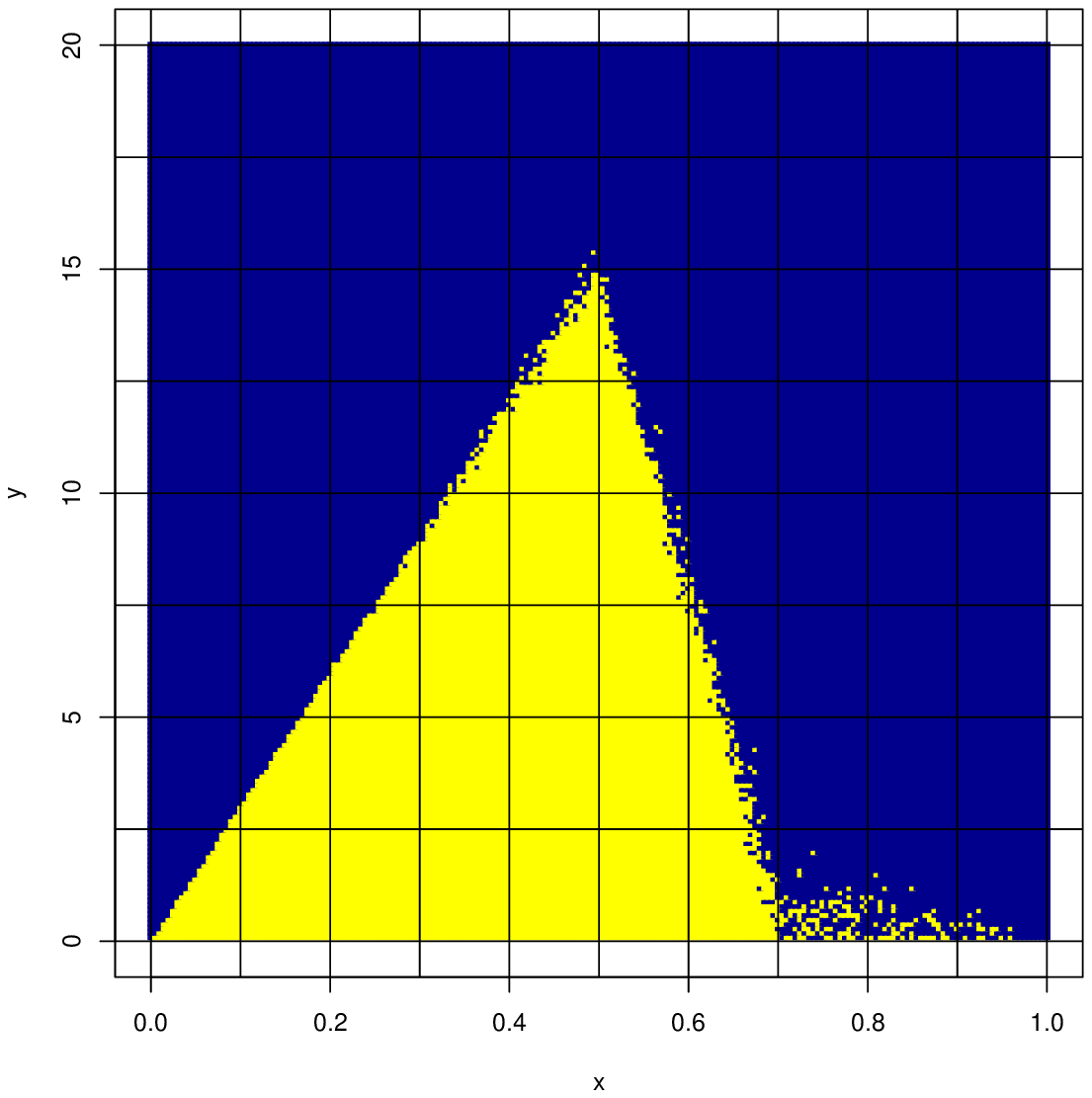}} \\
%
%
%
%\SetLabels
%\L (0.12*0.825) \small{\textbf{(d):} $m=1000$ samples} \\
%\endSetLabels
%\psfrag{y}[]{ sensing cost ($s$)}
%\psfrag{x}[]{ probability of existence ($r$)}
%\strut\AffixLabels{\includegraphics[width=5cm]{optimal_action_Pd_0.6_1000.ps}} &
%
\SetLabels
\L (0.12*0.815) \textbf{(c):} $m=1\mbox{,}000$ samples \\
\endSetLabels
\psfrag{y}[]{ sensing cost ($s$)}
\psfrag{x}[]{ probability of existence ($r$)}\hspace{0.5cm}
\strut\AffixLabels{\includegraphics[width=6cm]{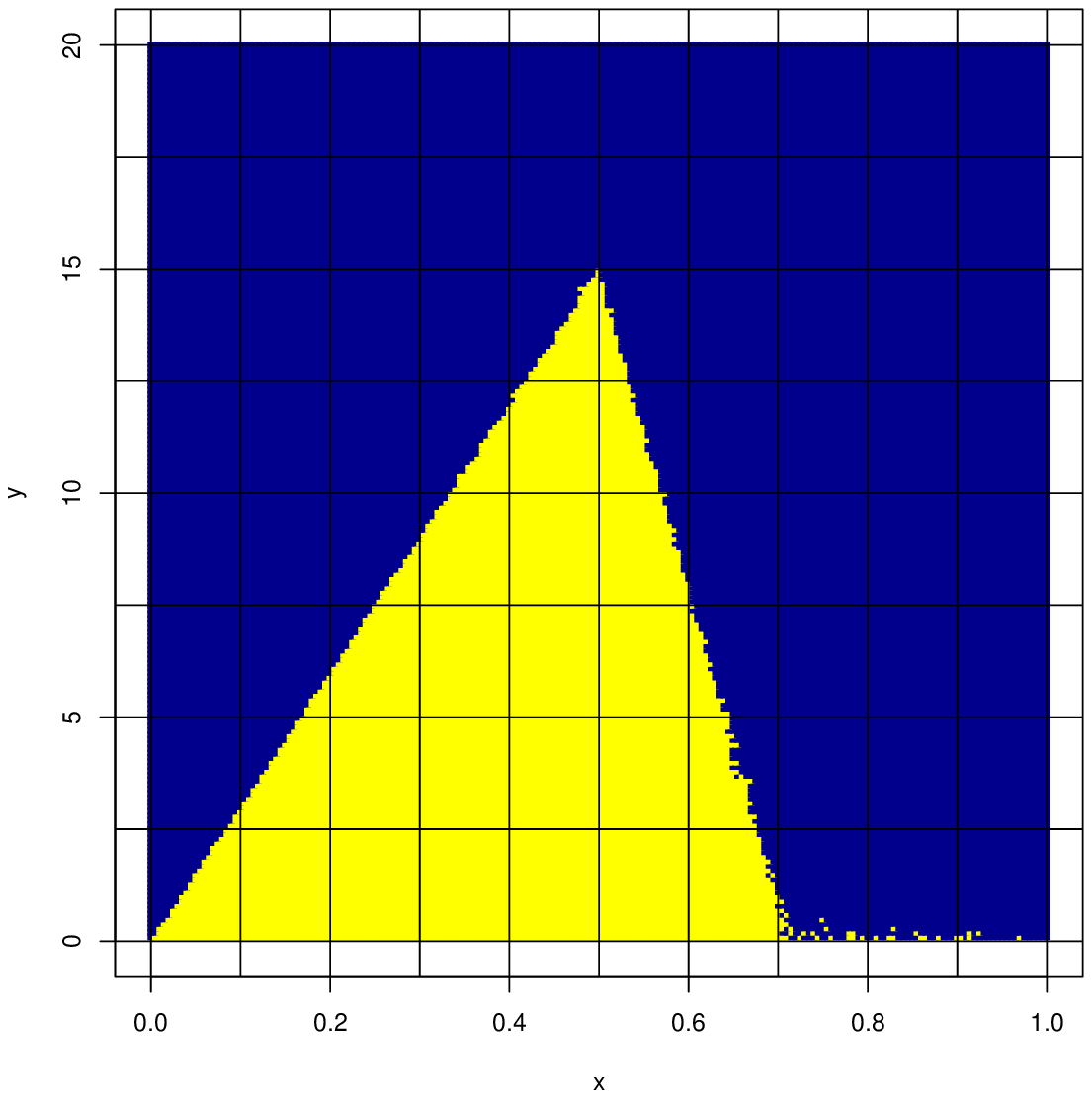}} &
\SetLabels
\L (0.12*0.815) \textbf{(d):} $n_h = 1$ sample \\
\endSetLabels
\psfrag{y}[]{ sensing cost ($s$)}
\psfrag{x}[]{ probability of existence ($r$)}
\vspace{-1.5cm}\hspace{0.5cm}
\strut\AffixLabels{\includegraphics[width=6cm]{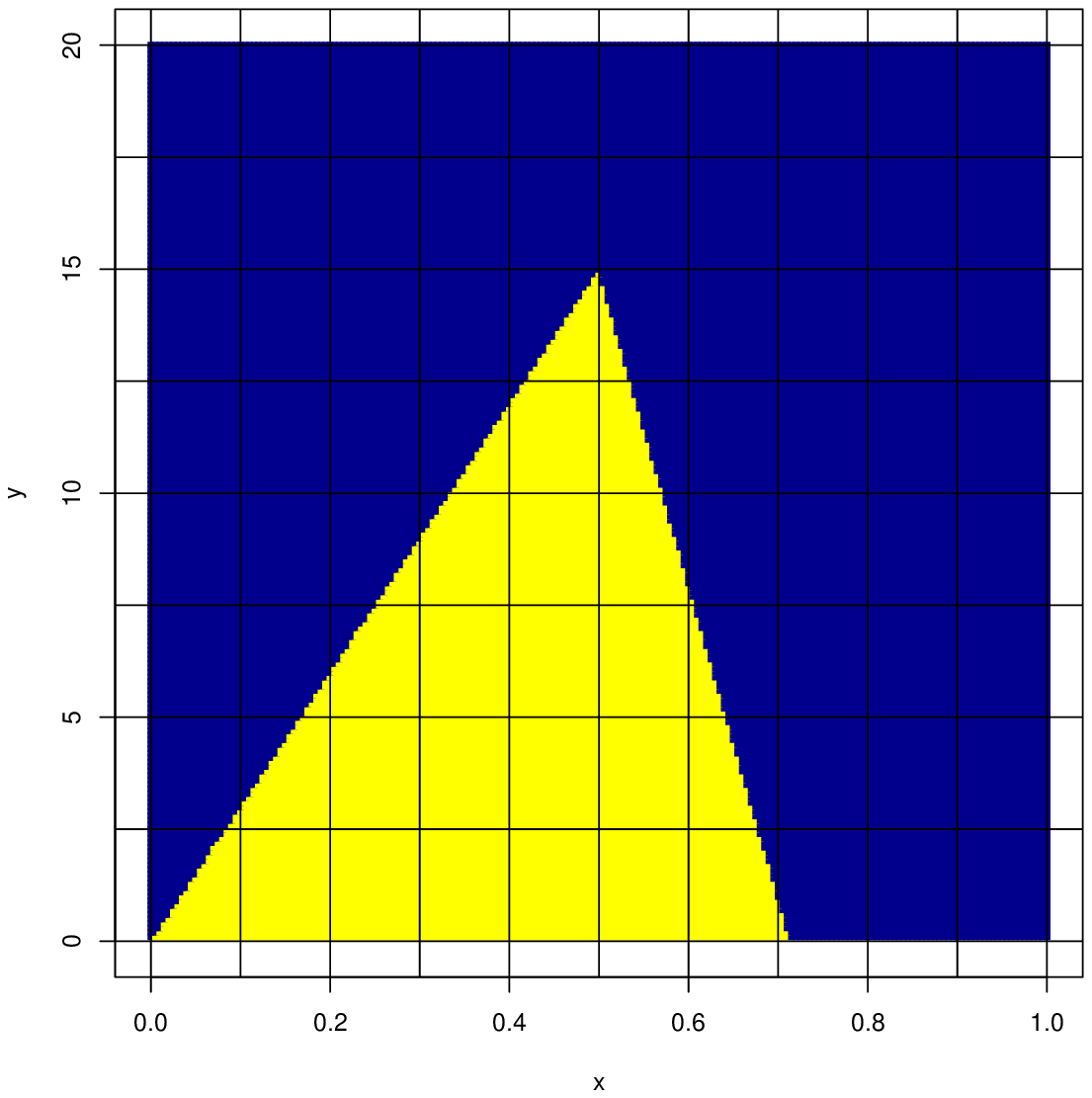}} \\
\vspace{0.2cm}
\end{array}
\]
\caption{Optimal action for each sensing cost and probability of target existence. Key: blue regions: do not make an observation/measurement, yellow regions: observe the potential target location. In \textbf{(a) -- (c):} the general sampling approach of Section \ref{general_sampling_approach}  are used. In \textbf{(d):} the efficient sampling approach of Section \ref{efficient_sampling_approach} is used.}
\label{figure_0}
\end{figure}

\vspace{0.1cm}
\noindent
The results are shown in Table \ref{table_0} and Figure \ref{figure_0}, with a range of values of the sensing cost $s\in (0, 20\mbox]$ km (which are added to the AMM-GOSPA error) and $r\in[0,1]$ evaluated. The solutions  generated using the efficient sampling approach shown in Figure \ref{figure_0}(d) are identical to the analytical results shown in Figure 2 in \cite{garcia_fernandez_2021}\footnote{The analytical results demonstrate that it is optimal to attempt a measurement if the probability of existence is not too low or too high, with the decision boundary dependent on the magnitude of the sensing cost $s$. For very high sensing costs, it is never optimal to attempt a measurement. Conversely. for $s=0$ (not shown, as $s\geq 0.1$ in Figure \ref{figure_0}), is is optimal to always attempt a measurement irrespective of the existence probability.}. Moreover, these solutions require only a single sample (i.e. $n_h=1$).

However, the general sampling  approaches 1--2 generate solutions that can be different at the decision boundary (see Figures \ref{figure_0}(a)--(c)), particularly if a relatively small number of samples are used.
Furthermore, the general  sampling approach that calculates the MMS-GOSPA from first principles, via equation (\ref{eq:posterior_MMSGOSPA}) (i.e. Approach 1) has a computational time  circa 250 times greater than the efficient sampling approach, in order to determine comparable (but still occasionally sub-optimal) solutions (e.g. compare Table \ref{table_0}(A): $m=1000$ with Table \ref{table_0}(C): $n_h=1$).

\begin{table}[H]
\vspace{1.0cm}
 \caption{Run-time per optimisation (averaged across each combination of $r$ and $s$ that is considered) as a function of the number of samples.
 \textbf{(a):} General sampling approach. \textbf{(b):} As in (a), but exploiting the analytical results of Section \ref{sec:analytical_multi}.  \textbf{(c):} The efficient sampling   approach again exploiting the analytical results.  All  computations were performed using an Intel\textregistered $\mbox{ }$Core\textsuperscript{TM}  i7-8750H (2.6GHz) processor. }
\begin{center}
%\vspace{-0.5cm}
\begin{tabular}{|c|c||c|c||c|c|} \hline
 \multicolumn{2}{|c||}{\textbf{(A)}} & \multicolumn{2}{c||}{\textbf{(B)}} & \multicolumn{2}{c|}{\textbf{(C)}}  \\ \hline
\multicolumn{1}{|c|}{$\ \ $\textbf{Samples $ (m)\ \ $}} & \multicolumn{1}{c||}{$\ $\textbf{Run-time (ms)$\ $}} &
\multicolumn{1}{c|}{$\ \ $\textbf{Samples $ (m) \ \ $}} & \multicolumn{1}{c||}{$\ $\textbf{Run-time (ms)$\ $}} &
\multicolumn{1}{c|}{$\ \ $\textbf{Samples $ (n_h)\ \ $}} & \multicolumn{1}{c|}{$\ $\textbf{Run-time (ms)$\ $}}
\\ \hline\hline
10 & 0.3 & 10 & 0.1 & 1 & 0.08 \\ \hline
100 & 2    & 100 &  0.3  %\multicolumn{2}{c|}{} % & 10  &  0.1
\\ \cline{1-4}
1,000 & 20 &   1,000   & 3   % & \multicolumn{2}{c|}{} % &     &
\\
\cline{1-4}
\end{tabular}
%\vspace{0.5cm}
 \label{table_0}
\end{center}
\end{table}

\newpage

%\newpage

\section{Non-Myopic (Multi-Step) Sensor Planning}
\label{sec_non_myopic}

\subsection{Suboptimal Control Approach}

\subsubsection{Motivation}
\label{sec:motivation}

In this section, a suboptimal planning approach is introduced. The primary reason it is suboptimal is that  due to
the way in which averaging across measurements is performed, the action calculated at each sampling time
depends only on the actions at the previous sampling times, and is independent of the measurements generated at
those times. Expressing this another way, the first action does not anticipate how the measurements that are subsequently generated will impact on future actions that are dependent on the realisation of the measurements.

This formulation is shown later (in Section \ref{sec:non_myopic_2}) to generate a solution with an
overall cost that is an upper bound on the overall cost of a solution calculated via a second recursive formulation
that conditions on previous measurements.
The second formulation is referred to as the ``optimal control
approach''.

\subsubsection{Overall Cost Function}

The cost function is the AMMS-GOSPA error, and it is assumed that there are no sensing costs. Let:
\begin{eqnarray}
 r_t(z_{1:t},a_{1:t}) & \triangleq &  \mbox{MMS-GOSPA}(z_{1:t}, a_{1:t}) \\
 V_{t:T}(a_{1:T})  & \triangleq  & {\mathbb E}_{z_{t:T}|z_{1:t-1}}\left[\sum_{l=t}^T\lambda^{l-1}r_l(z_{1:l},
a_{1:l})\right] \qquad
\end{eqnarray}
where $\mbox{MMS-GOSPA}(z_{1:t},
a_{1:t})$ is given by the minimum of (\ref{multi_MSGOSPA_1b}) and (\ref{multi_MSGOSPA_1c}).

\newpage
\noindent
The overall cost function  of the actions $a_{1:T}$ is then given as follows:
\begin{eqnarray}
C(a_{1:T})  & = & \mbox{AMMS-GOSPA}_{1:T}(a_{1:T})
\end{eqnarray}
It then follows from equation (\ref{simp_multi}) that $C(a_{1:T}) = V_{1:T}(a_{1:T})$.

\subsubsection{Optimisation}

An action ${a}^\star_{1}$ at each time step can be determined to minimise the  overall cost function $C(a_{1:T})$,
i.e.:
\begin{eqnarray}
{a}^\star_{1} & = & \argmin_{a_{1}}\Big[\min_{a_{2:T}}\big[ C(a_{1:T}) \big]\Big] \label{optimal_action_multi}
%a^\star_{1:t} &=& \argmin_{a_{1:T}}\Big[ C(a_{1:T}) \Big] \label{optimal_action_multi}
\end{eqnarray}
The minimum value of the cost function $C(a_{1:T})$  can then  be written as follows:
\begin{eqnarray}
\min_{a_{1:T}}\Big[C(a_{1:T})\Big] & = &
\min_{a_{1:T}}\Big[V_{1:T}(a_{1:T})\Big]  \ \triangleq \  V_{1:T}^\star \label{simon_10}
\end{eqnarray}
This is a commonly used approach (e.g. see \cite[equation~(24)]{bostrom_rost_2021} and \cite[equation~(28)]{bostrom_rost_2022}).

It is noted that each action $a_l$, $l=2,\ldots,T$ is dependent on the previous actions $a_1,\ldots,a_{l-1}$ but is
independent of the potential measurements, due to the expectation over $z_{1:T}$ in the cost function $C(a_{1:T})$.
It is also noted that actions $a_l$, $l=2,\ldots,T$ may never be realised, because at the next time step, the
situation awareness picture will be updated (it is dependent on the measurement generated following the action
$a_1^\star$),  and the optimisation (\ref{optimal_action_multi}) will be repeated for the new sliding window of $T$
time steps. This is the basis of standard receding horizon planning (e.g. see \cite{mac_2002}).

%\newpage

\subsection{Optimal Control Approach}
\label{sec:non_myopic_2}

\subsubsection{Formulation}

$V_{t:T}^\star$ can be written as follows:
\begin{eqnarray}
V_{t:T}^\star & = & \min_{a_{t:T}}{\mathbb E}_{z_{t}|z_{1:t-1}}\left[
{\mathbb E}_{z_{t+1:T}|z_{1:t}}\bigg[\sum_{l=t}^T\lambda^{l-1}r_l(z_{1:l}, a_{1:l})\bigg]\right] \label{new_multicol_1} \\
&& \qquad
 \mbox{(using Bayes' rule)}  \qquad\qquad \nonumber
 \\
&=& \min_{a_{t:T}}{\mathbb E}_{z_{t}|z_{1:t-1}}\left[\lambda^{t-1}r_t(z_{1:t}, a_{1:t}) +
{\mathbb E}_{z_{t+1:T}|z_{1:t}}\bigg[\sum_{l=t+1}^T\lambda^{l-1}r_l(z_{1:l}, a_{1:l})\bigg]\right] \\
&& \qquad \mbox{(using the fact that $r_t(z_{1:t}, a_{1:t})$ is conditional only on $z_{1:t}$ and $a_{1:t}$)}
\nonumber \\
&=& \min_{a_{t}}\Bigg\{{\mathbb E}_{z_{t}|z_{1:t-1}}\left[\lambda^{t-1}r_t(z_{1:t}, a_{1:t}) \right]
\label{eq:SM_10}
% \\
%&& \qquad \qquad
+
\min_{a_{t+1:T}}\bigg\{{\mathbb E}_{z_{t}|z_{1:t-1}}\Big[{\mathbb E}_{z_{t+1:T}|z_{1:t}}\big[\sum_{l=t+1}^T\lambda^{l-1}r_l(z_{1:l},
a_{1:l})\big]\Big]\bigg\}\Bigg\} \qquad  \\
&& \qquad \mbox{(using the fact that ${\mathbb E}_{z_{t}|z_{1:t-1}}\left[r_t(z_{1:t}, a_{1:t})\right]$ is independent of $a_{t+1:T}$)} \nonumber
\end{eqnarray}

\noindent
The inequalities (\ref{eq:SM_9}) -- (\ref{new_bellman_full}) are then true, via interchanging $\displaystyle\min_{a_{t+1:T}}$ and
${\mathbb E}_{z_t|z_{1:t-1}}$ in the second term of equation (\ref{eq:SM_10})\footnote{These inequalities follow from
the fact that $\displaystyle
 \min_{a_{t+1:T}}{\mathbb E}_{z_t|z_{1:t-1}}\left[ \Phi\right]  \geq
 {\mathbb E}_{z_t|z_{1:t-1}}\left[\min_{a_{t+1:T}}\left[ \Phi\right]\right] $
with equality only if the function $\Phi$ is independent of the variable $z_t$. Putting it another way, the minimum
of a sum is greater than the sum of the individual minimums.}:
\begin{eqnarray}
%
%\begin{eqnarray}
V_{t:T}^\star
&\geq & \min_{a_{t}}{\mathbb E}_{z_{t}|z_{1:t-1}}\left[\lambda^{t-1}r_t(z_{1:t}, a_{1:t}) +
\min_{a_{t+1:T}}{\mathbb E}_{z_{t+1:T}|z_{1:t}}\bigg[\sum_{l=t+1}^T\lambda^{l-1}r_l(z_{1:l}, a_{1:l})\bigg]\right]
\qquad \label{eq:SM_9} \\
&=& \min_{a_{t}}{\mathbb E}_{z_{t}|z_{1:t-1}}\Big[\lambda^{t-1}r_t(z_{1:t}, a_{1:t}) + V_{t+1:T}^\star\Big] \\
& \geq & \min_{a_{t}}{\mathbb E}_{z_{t}|z_{1:t-1}}\bigg[\lambda^{t-1}r_t(z_{1:t}, a_{1:t}) +
\min_{a_{t+1}}{\mathbb E}_{z_{t+1}|z_{1:t}}\Big[\lambda^{t}r_{t+1}(z_{1:t+1}, a_{1:t+1}) +
V_{t+2:T}^\star\Big]\bigg]\label{new_bellman} \qquad \\
& \vdots & %\qquad \qquad \mbox{(by recursion of inequality (\ref{eq:SM_9}), for $T\geq (t+2)$)}
\nonumber \\
& \geq & \lambda^{t-1}\min_{a_{t}}{\mathbb E}_{z_{t}|z_{1:t-1}}\Bigg[r_t(z_{1:t}, a_{1:t}) +
\lambda\min_{a_{t+1}}{\mathbb E}_{z_{t+1}|z_{1:t}}\bigg[r_{t+1}(z_{1:t+1}, a_{1:t+1})  \label{new_bellman_full} \\
&& \ + \lambda\min_{a_{t+2}}{\mathbb E}_{z_{t+2}|z_{1:t+1}}\Big[ r_{t+2}(z_{1:t+2},a_{1:t+2}) + \ldots +
\lambda\min_{a_T}{\mathbb E}_{z_{T}|z_{1:T-1}}\big[r_{T}(z_{1:T}, a_{1:T})\big]\Big]\bigg]\Bigg] \nonumber \\
&\triangleq & \hat{V}_{t:T}^\star \qquad \mbox{($=$ optimal overall cost)} \label{opt_cost_function}
\end{eqnarray}

\newpage

\noindent
Inequalities (\ref{new_bellman}) and (\ref{new_bellman_full}) (of Bellman type, e.g. see \cite{Bellman_1952})
follow by recursion of inequality (\ref{eq:SM_9}).
For the suboptimal control approach, $C_{1:T}^\star = V_{1:T}^\star$. Clearly this cost in an upper bound
on the Bellman type equation  (\ref{new_bellman_full}) for all values of $T$.

Motivated by this finding,  the ``optimal control approach'' determines each action $\hat{a}_1^\star$ as follows:
 \begin{eqnarray}
\hat{a}^\star_1 & = & \argmin_{a_{1}}{\mathbb E}_{z_{1}}\Bigg[r_1(z_{1}, a_{1}) +
\lambda\min_{a_{2}}{\mathbb E}_{z_{2}|z_{1}}\bigg[r_{2}(z_{1:2}, a_{1:2}) \label{optimal_Bellman_act}  \\
&& \ + \lambda\min_{a_{3}}{\mathbb E}_{z_{3}|z_{1:2}}\Big[ r_{3}(z_{1:3},a_{1:3}) + \ldots +
\lambda\min_{a_T}{\mathbb E}_{z_{T}|z_{1:T-1}}\big[r_{T}(z_{1:T}, a_{1:T})\big]\Big]\bigg]\Bigg]\qquad \nonumber
\end{eqnarray}

 \noindent
 It is important to note that each minimisation (i.e. minimum with regard to $a_{t}$) in equation
(\ref{optimal_Bellman_act}) is conditional on both the previous actions $a_{1:t-1}$ and the previous measurements
$z_{1:t-1}$.

As an illustrative example of the approach, for $T=2$:
\begin{eqnarray}
\hat{a}^\star_{1} &=&
\argmin_{a_{1}}{\mathbb E}_{z_{1}}\Big[r_1(z_{1}, a_{1}) \label{astar_2step}
+ \lambda\min_{a_2} {\mathbb E}_{z_{2}|z_1}
\big[r_2(z_{1:2}, a_{1:2})   \big]\Big]
\end{eqnarray}

\noindent
where:
\begin{eqnarray}
\hat{a}_2^\star(z_1,a_1) & = & \argmin_{a_2}{\mathbb E}_{z_{2}|z_1} \big[r_2(z_{1:2}, a_{1:2})
 \big]\qquad\qquad\label{a_2_hat}
 \end{eqnarray}

 %\newpage
 \noindent
 and:
 \begin{eqnarray}
 \hat{V}_{2:2}^\star &=& \lambda\min_{a_2}{\mathbb E}_{z_{2}|z_1} \big[ r_2(z_{1:2}, a_{1:2})\big] \\
 \hat{V}_{1:2}^\star &=& \min_{a_{1}}{\mathbb E}_{z_{1}}\Big[r_1(z_{1}, a_{1}) + \hat{V}_{2:2}^\star \Big] \qquad
\end{eqnarray}

\noindent
When using the ideal measurement  approach (but allowing $P_d<1$), the optimal and suboptimal approaches generate identical actions, i.e. $a^\star_1 = \hat{a}^\star_1$. This follows from the fact that in this case, the MMS-GOSPA error at each time $t$ is only non-zero if measurements are not generated at times 1 through $t$ (see Appendix for full details).

%
%
%\item Three-step planning  ($T=3$):
%\begin{eqnarray}
%\hat{a}^\star_{1} &=&
%\argmin_{a_{1}}{\mathbb E}_{z_{1}}\Big[r_1(z_{1}, a_{1}) + \lambda\min_{a_2} {\mathbb E}_{z_{2}|z_1} \big[
%r_2(z_{1:2}, a_{1:2})  + {\hat V}_{3:3}^\star \big]\Big]\label{astar_3step} \qquad \end{eqnarray}
%where:
%\begin{eqnarray}
%{\hat V}_{3:3}^\star &=& \min_{a_3}{\mathbb E}_{z_{3}|z_{1:2}}\Big[\lambda^2 r_3(z_{1:3},a_{1:3})\Big]
%\label{astar_3step2}  \\
%%
%\hat{a}_3^\star(z_{1:2}, a_{1:2}) &=&
%\argmin_{a_3}{\mathbb E}_{z_{3}|z_{1:2}}\Big[r_3(z_{1:3},a_{1:3})\Big]\label{astar_3step4}
%\\
%%
%{\hat V}_{2:2}^\star &=& \min_{a_2}{\mathbb E}_{z_{2}|z_1} \big[\lambda r_2(z_{1:2}, a_{1:2})  + {\hat V}_{3:3}^\star
%\big]\label{astar_3step3}
%\\
%%
%\hat{a}_2^\star(z_{1}, a_{1}) &=& \argmin_{a_2} {\mathbb E}_{z_{2}|z_1} \big[ r_2(z_{1:2}, a_{1:2})  +
%{\hat V}_{3:3}^\star \big]\qquad\qquad\\
%{\hat V}_{1:1}^\star &=& \min_{a_{1}}{\mathbb E}_{z_{1}}\Big[r_1(z_{1}, a_{1}) + {\hat V}_{2:2}^\star\Big]\qquad
%\mbox{($=$ optimal overall cost)}\qquad \label{3step_final}
%\end{eqnarray}
%\end{itemize}
%

\subsubsection{Calculating the Conditional AMMS-GOSPA Error}
\label{sec:cond_AMMS_GOSPA}

The optimal control approach requires calculation of the AMMS-GOSPA at each time $t+1$, conditional on the previous
measurements $z_{1:t}$  and all actions $a_{1:t+1}$.
Similar to the calculation of the (unconditional) AMMS-GOSPA in Section \ref{sec:AMMS_GOSPA}, the conditional
AMMS-GOSPA is given as follows:
\begin{eqnarray}
 && {\mathbb E}_{z_{t+1}|z_{1:t}}\big[\mbox{MMS-GOSPA}(z_{1:t+1}, a_{1:t+1})\big] \nonumber \\
 &&  =  \ \int_{z_{t+1}} \mbox{MMS-GOSPA}(z_{1:t+1}, a_{1:t+1})  p(z_{t+1}|z_{1:t}, a_{1:t+1}) dz_{t+1} \\
 && =  \ \sum_{i=0}^n \int_{z_{t+1}} \mbox{MMS-GOSPA}(z_{1:t+1}, a_{1:t+1})  \label{170}  p(z_{t+1}|x_i,z_{1:t}, a_{1:t+1}) p(x_i|z_{1:t},a_{1:t}) dz_{t+1}  \\
&&  =  \    \sum_{i=0}^n \hat{w}_i
\int_{z_{t+1}} \mbox{MMS-GOSPA}(z_{1:t+1}, a_{1:t+1})\label{125b}  p(z_{t+1}|x_i,a_{t+1}) dz_{t+1}
\end{eqnarray}

\noindent
where $\hat{w}_i = p(x_i|z_{1:t},a_{1:t})$ is the posterior probability of each target hypothesis (including ``no
target present''), conditional on the previous measurements $z_{1:t}$  and previous  actions $a_{1:t}$.
Equation (\ref{125b}) is approximated as in equation (\ref{new_det_approach_0}) but only requires samples to be
generated from $p(z_{t+1}| x_i, a_{t+1})$,
This approximation of the conditional AMMS-GOSPA error is given as follows:
\begin{eqnarray}
 && {\mathbb E}_{z_{t+1}|z_{1:t}}\left[\mbox{MMS-GOSPA}(z_{1:t+1}, a_{1:t+1})\right] \\
  &&  \approx \  \frac{1}{n_h}\sum_{i=0}^n\hat{w}_i \sum_{m=0}^{1} Pr\big(s_{t+1}=m|x_i,a_{t+1}\big)
 \sum_{j=1}^{n_h} \mbox{MMS-GOSPA}\big(z_{1:t}, z_{t+1}^{ij}(m), a_{1:t+1}\big) \qquad \nonumber
\end{eqnarray}

\noindent
where to remind the reader, $s_{t+1}=1$ denotes a target  detection  at sampling time $t+1$.
$Pr(s_{t+1}=m|x_i,a_{t+1})$ is the probability that $m=0,1$ target measurements are generated on the $(t+1)$-th time step, conditional on the target state being given by $x_i$ and action $a_{t+1}$. $Pr\big(s_{t+1}(m)|x_i,a_{t+1}\big)$ is again calculated via equation (\ref{prob_target_meas}).
Each measurement sample
$z_{t+1}^{ij}(m)$ can be a vector of multiple measurements, which includes a maximum of one target generated measurement (sampled from ${\mathcal N}(x_i, \Sigma)$ if $s_{t+1}=1$) plus false alarms.

\subsubsection{Implementation of the Optimal Control Approach}

Using the efficient sampling approach,
for two-step planning, the optimal action $a_1^\star$ (given by equation (\ref{optimal_Bellman_act}))  is as follows:
\begin{eqnarray}
\hat{a}_1^\star &=& \argmin_{a_1} \frac{1}{n_h}\sum_{j=1}^{n_h}\sum_{m=0}^1\sum_{i=0}^n\Bigg\{
p(x_i) Pr(s_1=m|x_i,a_1) \bigg\{r_1(z_1^{ij}(m), a_1)
\qquad\label{optimal_action_clutter} \\
&& \qquad \qquad \qquad \qquad \qquad \qquad \qquad + \lambda  \min_{a_2} {\mathbb E}_{z_2|z_1^{ij}(m)}
\Big[r_2(z_1^{ij}(m),z_2, a_{1:2}))\Big]\bigg\}\Bigg\} \nonumber \\
  & =  &
 \argmin_{a_1} \frac{1}{n_h}\sum_{j=1}^{n_h}\sum_{m=0}^1\sum_{i=0}^n
p(x_i) Pr(s_1=m|x_i,a_1) \bigg\{r_1(z_1^{i j}(m), a_1)  \label{complete_Bellman}  \\
&& \qquad  + \min_{a_2} \frac{\lambda}{n_h}\sum_{\underline{j}=1}^{n_h}\sum_{\underline{m}=0}^1\sum_{\underline{i}=0}^n
p(x_i|z_1^{ij}(m)) Pr(s_2=\underline m|x_{\underline i},a_2)
r_2\big(z_1^{ij}(m),z_2^{\underline{i}\underline{j}}(\underline{m}), a_{1:2}\big)
\bigg\}\qquad\nonumber
\end{eqnarray}
where $r_t$ denotes the cost at time $t$ conditional on the  actions $a_{1:t}$ and measurements $z_{1:t}$ generated at times 1 through $t$. $Pr(s_1=m|x_i,a_1)$ is the probability that $m=0 \mbox{ or }1$ target measurements are generated on the first time step, conditional on the target state being given by $x_i$ and action $a_1$. This is calculated via equation (\ref{prob_target_meas}).
The measurement samples $z_1^{ij}(m)$ include a maximum of one target generated measurement (sampled from ${\mathcal N}(x_i, \Sigma)$) plus a Poisson distributed number of false alarms.

The corresponding optimal total cost incurred (the AMMS-GOSPA error, given by equations (\ref{new_bellman_full}) -- (\ref{opt_cost_function})) is then given as follows:
\begin{eqnarray}
\hat{V}^\star_{1:2} &=& \frac{1}{n_h}\sum_{j=1}^{n_h}\sum_{m=0}^1\sum_{i=0}^n\Bigg\{
p(x_i) Pr(s_1=m |x_i,a_1) \bigg\{r_1(z_1^{ij}(m), a_1^\star) \label{V_hat_T_2} \\
&& \qquad \qquad \qquad \qquad \qquad
+ \lambda  {\mathbb E}_{z_2|z_1^{ij}(m)}\Big[
r_2(z_1^{ij}(m),z_2, a_1^\star, a_2^\star\big(a_1^\star, z_1^{ij}(m))\big)\Big]
\bigg\}\Bigg\}\nonumber
\end{eqnarray}

\noindent
Equations (\ref{optimal_action_clutter}) -- (\ref{V_hat_T_2}) generalise in an obvious manner for $T>2$.

\subsubsection{Computational Complexity}

The optimisation  (\ref{optimal_action_clutter}) is computationally expensive, because for each potential action $a_1$, the minimisation of each second time step action $a_2$ must be performed for each sampled first time step measurement $z_1^{ij}(m)$. Hence, the number of minimisations that must be performed at the second time step is $(2n+1)n_h$ for each potential action $a_1$ (i.e. $2n_h$ minimisations corresponding to hypotheses $i>0$; and $n_h$ minimisations for $i=0$, because only $s_1=0$ is possible for $i=0$). Hence the total number of minimisations is $n_a(2n+1)n_h$, where $n_a$ is the number of possible actions. This generally makes the approach computationally prohibitive for time horizons greater than two time steps\footnote{E.g. for three-step planning, the number of minimisations performed on the third step alone is $\left((2n+1)n_hn_a\right)^2$.}, unless the algorithm is parallelised.

By way of comparison, the suboptimal approach performs just a single minimisation (given by equation (\ref{optimal_action_multi})), with $n_a^T$ potential combinations of actions for $T\geq 1$. Hence, for $T=2$ the computational expense of the suboptimal algorithm is typically much lower (assuming that $n_a \ll (2n+1)n_h$), though it should be noted that in the suboptimal approach, the cost function in each minimisation is computationally more expensive as it is the sum total of the costs incurred across all  time steps.

\vspace{0.25cm}

\noindent
\textit{Special case} -- If $\lambda_{FA}=0$ and $\Sigma\approx 0$ it can easily be shown that $Pr(s_1=1|x_i,a_1)r_1(z_1^{ij}(1),a) = 0$ for all values of $i$. This is due to the fact that whenever a measurement is generated by a target hypothesis, $r_1(z_1^{ij}(1),a)=0$, as a target can be inferred to be present  without geo-location or cardinality errors. As a result, in the optimisation (\ref{complete_Bellman}) one can set $n_h=1$, $m=\underline{m}=0$ and $z_1=z_2=\phi$. Consequently, the minimisation on the second time step need only to be performed if no measurement is generated on the first time step, and the optimisation  (\ref{complete_Bellman}) reduces to equation (\ref{optimal_Bellman_act_app1}) in the Appendix.
If a measurement is generated, it is not necessary to attempt a second observation.
Therefore, only one minimisation need  be performed on the second time step for each potential action $a_1$. To remind the reader, it is shown in the Appendix that in this case, the optimal and suboptimal actions are the same on the first time step (i.e. $\hat{a}_1^\star = a_1^\star$), irrespective of the length of the planning horizon.

\subsection{Baseline Approach -- Minimisation of  Target Localisation Error}

As a baseline for comparison, the optimal action sequence is determined in order to minimise the target localisation error within the time window under consideration. Minimisation of the localisation error is a widely used approach in target tracking, with the PCRB often used to predict future performance (e.g. again see \cite{hernandez_2013}).
In the current paper, the average estimated target location root mean squared error (RMSE),  conditional on a target existing, is used to quantify performance. This metric is similar to the PCRB, and takes into account the potential for missed detections, as well as the impact of false alarms and measurement errors.

Similar to the calculation of the AMMS-GOSPA  error, but not considering cardinality errors,
the average target estimated location mean squared error (MSE) for $T=2$ is calculated as follows:
\begin{eqnarray}
 \bar{\mathcal{MSE}}_{1:2}(a_{1:2}) & = &
\frac{1}{n_h}\sum_{j=1}^{n_h}\sum_{m=0}^1\sum_{i=1}^n
w_i Pr(s_1=m|x_i,a_1) \bigg\{  {\mathcal SE}_1(z_1^{i j}(m), a_1, i) \label{eq_MSE}   \\
&& \qquad +  \frac{\lambda}{n_h}\sum_{\underline{j}=1}^{n_h}\sum_{\underline{m}=0}^1\sum_{\underline{i}=1}^n
p(x_i|z_1^{ij}(m), x_i\neq x_0) Pr(s_2=\underline m|x_{\underline i},a_2)
{\mathcal SE}_2(z_1^{ij}(m),z_2^{\underline{i}\underline{j}}( \underline{m}), a_{1:2}, \underline{i})\bigg\}
\nonumber
\end{eqnarray}
where ${\mathcal SE}_t(z_{1:t}, a_{1:t}, i)$ is the squared distance between the particle filter based posterior state estimate $p(X|z_{1:t},a_{1:t})$ and the hypothesis $x_i$.
It is noted that unlike, e.g., (\ref{complete_Bellman}), the average mean squared error (\ref{eq_MSE}) only considers hypotheses under which a target exists (i.e. giving $i>0$ and $\bar{i}>0$ in the summations in (\ref{eq_MSE})). The first term on the right-hand side of equation (\ref{eq_MSE}) gives the average MSE on the first time step (denoted  $\bar{\mathcal{MSE}}_{1}(a^\star_{1})$), with the second term giving the average MSE on the second time step (denoted  $\bar{\mathcal{MSE}}_{2}(a^\star_{1:2})$). Equation (\ref{eq_MSE}) generalises for  $T>2$ in a straightforward manner via adding the average MSE  at subsequent time steps.

The baseline control approach then determines the action ${a}^b_{1}$ at each time step to minimise $\bar{\mathcal{MSE}}_{1:T}(a_{1:T})$ ,
i.e.:
\begin{eqnarray}
{a}^b_{1} & = & \argmin_{a_{1}}\Big[\min_{a_{2:T}}\big[ \bar{\mathcal{MSE}}_{1:T}(a_{1:T}) \big]\Big] \label{optimal_action_baseline}
\end{eqnarray}
This is analogous to the suboptimal control approach (i.e. it does not use a full Bellman recursion).

%$\bar{\mathcal{MSE}}_{1:2}(a^\star_{1:2})$, $\bar{\mathcal{MSE}}_{1}(a^\star_{1})$ and $\bar{\mathcal{MSE}}_{2}(a^\star_{1:2})$  are each averaged over a number of simulations (nominally 20) with the average RMSEs (denoted by $\bar{\mathcal{RMSE}}_{1:2}(a^\star_{1:2})$, $\bar{\mathcal{RMSE}}_{1}(a^\star_{1})$ and $\bar{\mathcal{RMSE}}_{2}(a^\star_{1:2})$ respectively) given by the square roots of the corresponding  MSE values.

\section{Simulations}
\label{sec_sim}

\subsection{Scenario Specification}

Three target distributions are considered, unimodal, bimodal and trimodal, representing scenarios concerned with sensor control for wide-area search. Target hypotheses are sampled from these distributions as  outlined in Table \ref{table_dist}. In each case, the hypotheses are given equal weights (i.e. $w_i=1/n$) and are time-invariant (i.e. stationary).

\begin{figure}[H]
\vspace{-1.5cm}
\[
\begin{array}{ccc}%
\SetLabels
\L (0.12*0.825) \textbf{(a):} Unimodal \\
\endSetLabels
\psfrag{y}[]{\small{Northings (km)}}
\psfrag{x}[]{\small{Eastings (km)}}
\strut\AffixLabels{\includegraphics[width=4cm]{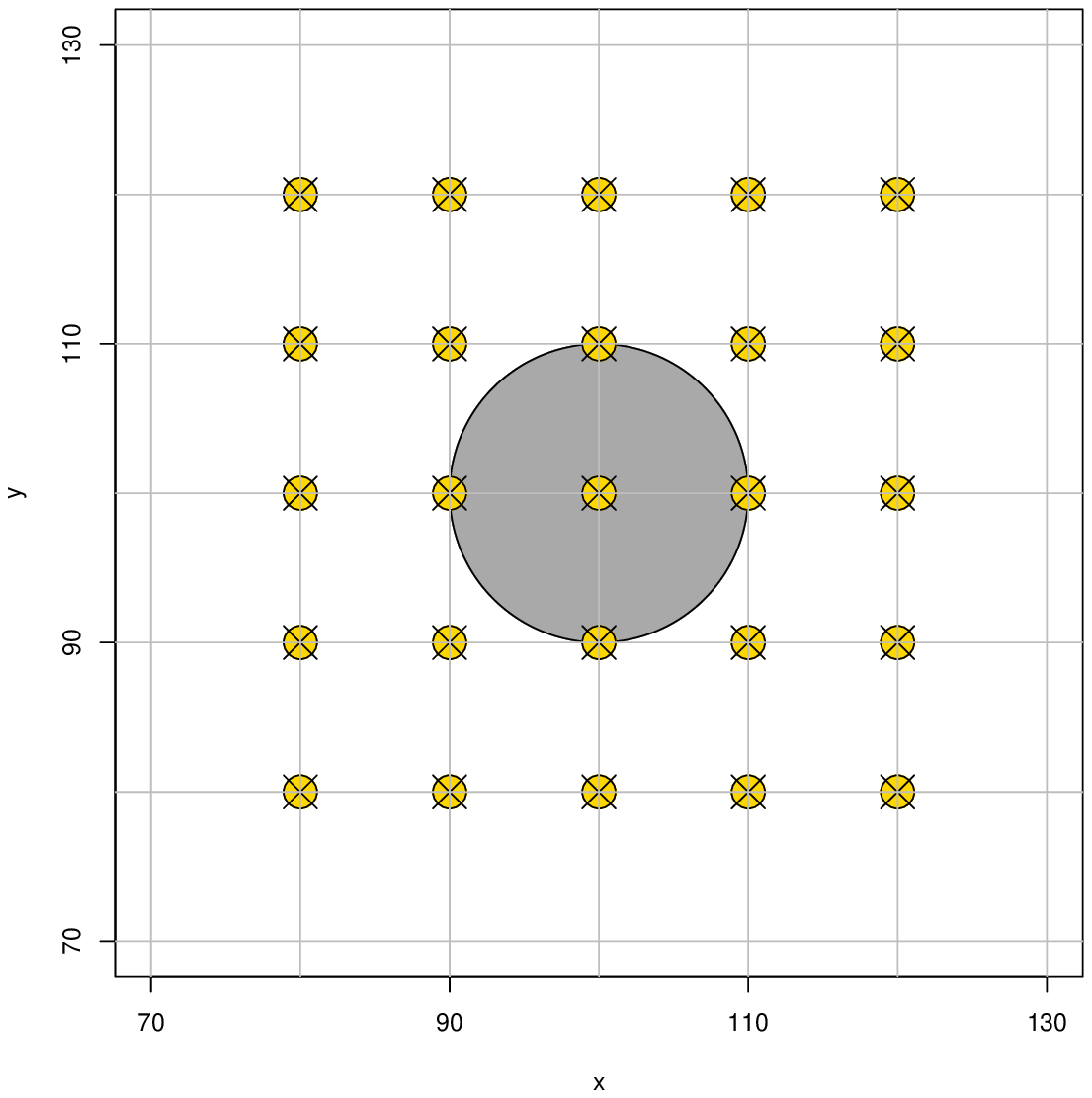}} &
\SetLabels
\L (0.12*0.825) \textbf{(b):} Bimodal \\
\endSetLabels
\psfrag{y}[]{\small{Northings (km)}}
\psfrag{x}[]{\small{Eastings (km)}}\hspace{0.5cm}
\strut\AffixLabels{\includegraphics[width=4cm]{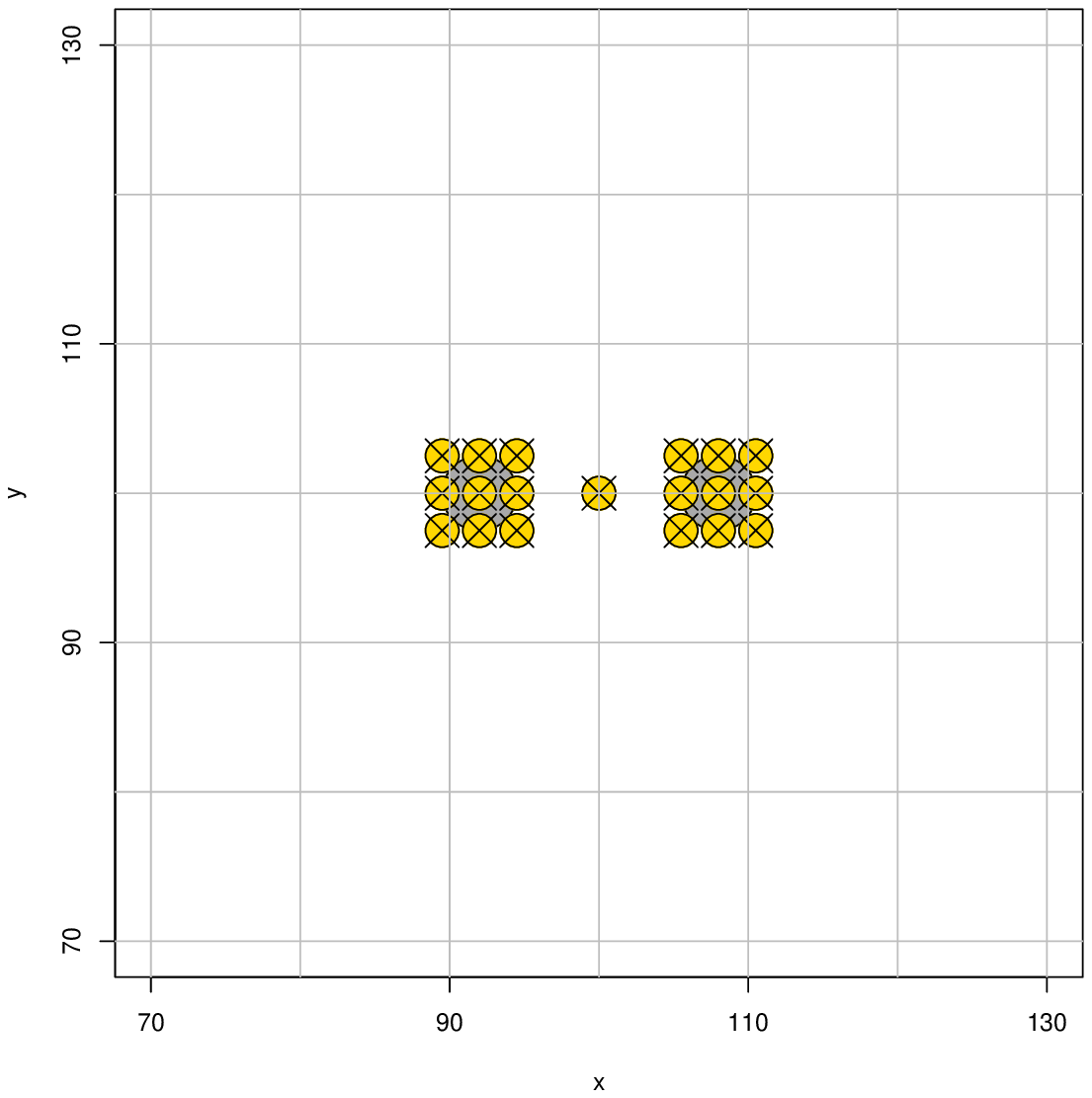}} &
\SetLabels
\L (0.12*0.825) \textbf{(c):} Trimodal \\
\endSetLabels
\psfrag{y}[]{\small{Northings (km)}}
\psfrag{x}[]{\small{Eastings (km)}}
\vspace{-0.5cm}\hspace{0.5cm}
\strut\AffixLabels{\includegraphics[width=4cm]{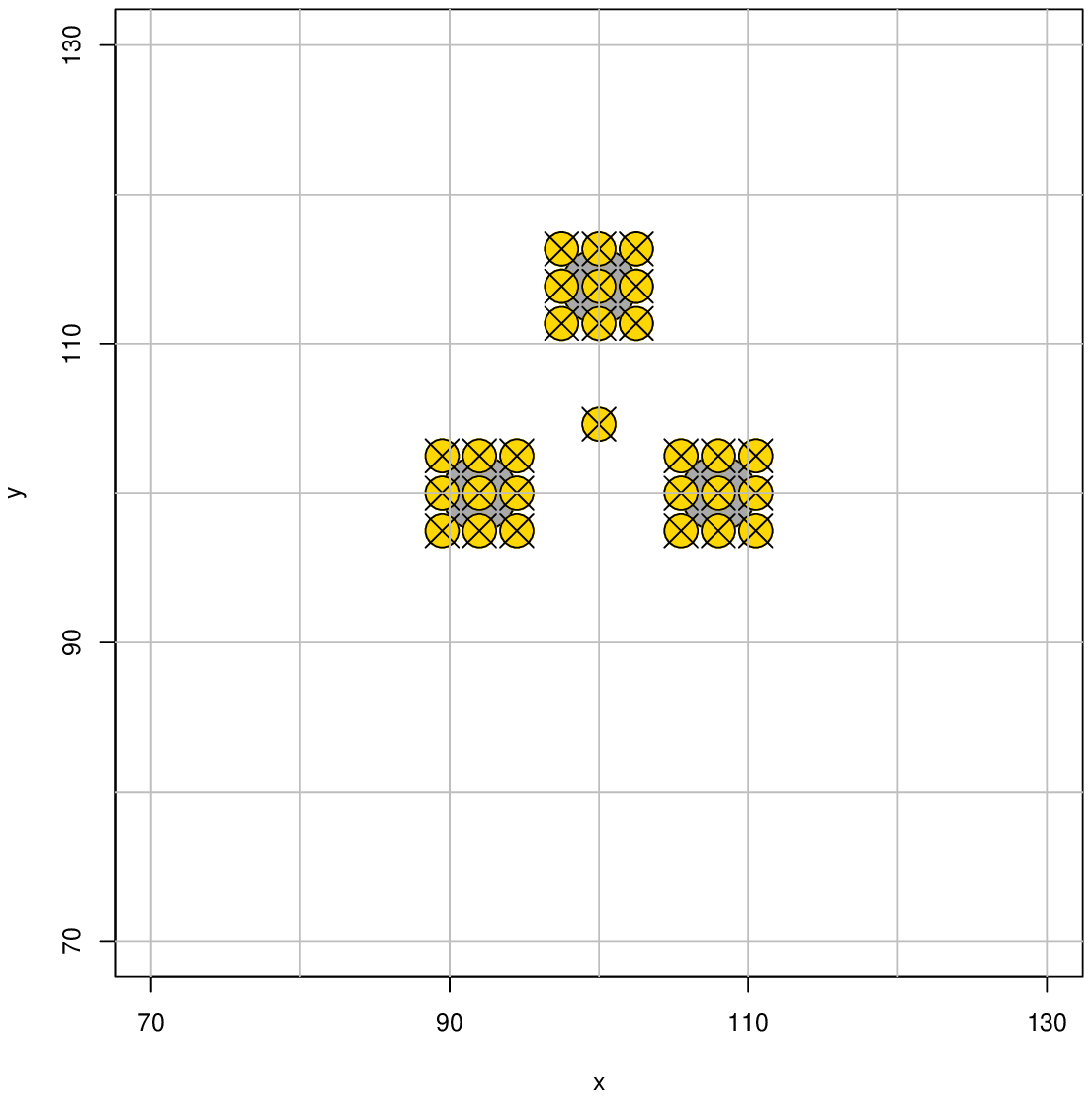}} \\
\end{array}
\]
\caption{Action hypotheses for the three scenarios. In each case, the yellow circles show the potential actions, with the cross showing the ``spotlight'' centre. The grey regions show the one standard deviation uncertainty regions for each target distribution. Note that the sensor FOV is circular, centred on the spotlight centre with a FOV of 10km.}
\label{figure_2}
\end{figure}

A single  sensor has a circular FOV (i.e. operates in a ``spotlight'' mode) of radius 10km centred on the action $a\in A$ (which specifies the coordinates of the centre of the spotlight).
The action hypotheses consist of  uniformly spaced spotlight centres that overlay each target mode, as shown
 in Figure \ref{figure_2}. In each scenario
there is also the action hypothesis of not attempting  to make a target observation\footnote{The default  is to not  attempt an observation unless doing so offers a performance improvement. As noted earlier, when $\lambda_{FA}=0$ and $\Sigma\approx 0$, observations are only required up to the time instance at which the target is detected.}.
Conditional on target hypothesis $x_i$ being true and action $a$ being performed, a measurement of the Cartesian coordinates of the target is generated
with probability $P_d$ if the hypothesis is within the FOV of the sensor. Each measurement error has a zero-mean Gaussian distribution with covariance $\Sigma$. There are either no  false alarms (i.e. $\lambda_{FA}=0$) or false alarms are generated with rate $\lambda_{FA}=0.01$ per unit volume of the sensor FOV.
Parameter settings are summarised in Table \ref{table_par}.

\begin{table}[H]
 \vspace{0.75cm}
  \caption{Sampling distributions of the target hypotheses in the three scenarios considered.}
\begin{center}
%\vspace{-0.5cm}
\begin{tabular}{|c|c|c|c|} \hline
 &  \textbf{ Number of }
& \multicolumn{2}{c|}{\textbf{ Sampling distribution }}  \\
\raisebox{1.25ex}[0cm][0cm]{\textbf{ Target distribution }}  & \multicolumn{1}{c|}{\textbf{hypotheses ($n$)}} &    \multicolumn{2}{c|}{\textbf{$\left(x_i  \sim   {\mathcal N}(\bar{x}, \Sigma_x)\right)$}} \\
\hline \hline
Unimodal &  100 & $\bar{x} = (100\mbox{km} \ 100\mbox{km})'$ & $\Sigma_x = \mbox{diag}(100\mbox{km}^2, 100\mbox{km}^2)$\\
\hline\hline
&  50 & $\bar{x} = (92\mbox{km} \ 100\mbox{km})'$  &  \\
\cline{2-3}
\raisebox{1.25ex}[0cm][0cm]{Bimodal} &  50 & $\bar{x} = (108\mbox{km} \ 100\mbox{km})'$ & \raisebox{1.25ex}[0cm][0cm]{$\Sigma_x = \mbox{diag}(2.5^2\mbox{km}^2, 2.5^2\mbox{km}^2)$}
\\
  \hline\hline
&  33 & $\bar{x} = (92\mbox{km} \ 100\mbox{km})'$ & \\
\cline{2-3}
Trimodal  &  33 & $\bar{x} = (108\mbox{km} \ 100\mbox{km})'$  &  $\quad\Sigma_x = \mbox{diag}(2.5^2\mbox{km}^2, 2.5^2\mbox{km}^2)\quad$ \\
\cline{2-3}
  &  33 & $\quad\bar{x} = (100\mbox{km} \ (100 + \sqrt{192})\mbox{km})'\quad$ & \\
\hline
\end{tabular}
\vspace{0.0cm}
 \label{table_dist}
\end{center}
\end{table}

\newpage

\begin{table}[H]
 \vspace{-0.0cm}
 \caption{Summary of parameter settings used in the simulations.}
\begin{center}
%\vspace{-0.5cm}
\begin{tabular}{|l|l|} \hline
\multicolumn{1}{|c|}{\textbf{Parameter}} & \multicolumn{1}{c|}{\textbf{Value}} \\ \hline \hline
Number of time steps, $T$ & 1, 2 or 3   \\ \hline
Number of target hypotheses, $n$ & 100 (unimodal, bimodal) or 99 (trimodal) $\qquad\qquad$ \\ \hline
Number of possible actions, $n_a$ & 26 (unimodal), 20 (bimodal), 29 (trimodal) \\\hline
Probability that a target exists, $r$ & 0.8  \\\hline
Probability of detection, $P_d$ & 0.6, 0.9 or  1.0  \\ \hline
False alarm rate, $\lambda_{FA}$ (m$^{-2}$)& 0.0 or 0.01 $\quad$  \\ \hline
Measurement error standard deviation, $\sigma$  &  10$^{-5}$km $\qquad$ (when $\lambda=0$)  $\qquad\qquad$ \\
$\left(\Sigma=\mbox{diag}(\sigma^2,\sigma^2)\right)$   &  10$^{-2}$km $\qquad $ (when $\lambda=0.01$)\\ \hline
Sensor FOV & circular, with radius 10km  \\ \hline
Measurement samples $n_h$ per target hypothesis$\qquad\qquad$  & 1 \ \ \  (when $\lambda_{FA}=0$) \\
& 10 \   (when $\lambda_{FA}=0.01$, baseline/suboptimal approaches)$\qquad$ \\
& 1 \ \ \  (when $\lambda_{FA}=0.01$, optimal approaches)
 \\ \hline
Cardinality error cost, $c$ & 10km  \\ \hline
Discount factor  for non-myopic planning, $\lambda$ & 1.0 $\qquad$ \\
%GOSPA metric parameter, $p$  & 2.0 \\
%GOSPA metric parameter, $\alpha$  & 2.0 \\ \hline
%
\hline
\end{tabular}
\vspace{0.5cm}
 \label{table_par}
\end{center}
\end{table}

%\hrule

%%%%%%%%%%%%%%%%%%%%%%%%%%%%%%%%%%%%%%%%%%%%%%%%%%%%%%%%%%%%%%%%%%%%%%%%%%%%%%%%%%%%%%%%%%%

\subsection{Simulation Results}
\label{sec:multistep_results}

\subsubsection{Baseline vs. GOSPA-Based Control Approaches}

Exemplar optimal actions for $P_d \in\{0.6,  0.9, 1.0\}$ and $T=\{2,3\}$ are shown in Figure \ref{figure_4}, for each of the target distributions and for scenarios with no clutter (i.e. $\lambda_{FA}=0$) and insignificant measurement errors (i.e. $\Sigma=\mbox{diag}(10^{-10}\mbox{km}^2, 10^{-10}\mbox{km}^2)$). The second and third time step actions (shown in green and blue respectively) are only necessary if a measurement is not generated on any previous time step. If a measurement is generated at any time, the target state can be estimated extremely accurately and without cardinality error (because there is
no clutter, and so the presence of a measurement indicates that a target is present), resulting in an AMMS-GOSPA error of zero.
The suboptimal approach generates identical solutions on the first time step (due to the ideal measurement set assumptions). However,  the actions on subsequent time steps  are independent of the previous measurements, as discussed in Section \ref{sec:motivation}. As a result, under the suboptimal approach the second and third time step actions are  made irrespective of whether previous measurements have been generated.

Overall location RMSE and AMMS-GOSPA values for the baseline and optimal  approaches are shown in Table \ref{table_costs}. Because there is no cost of sensing, the AMMS-GOSPA cost is identical for the suboptimal and optimal approaches. It can be seen that although, by design, the baseline approach generates slightly lower RMSE values, the corresponding AMMS-GOSPA values are significantly higher than for the GOSPA-based approaches. Hence, the baseline approach does a poor job in balancing the tradeoff between estimation accuracy and cardinality errors, with cardinality errors far greater than for the GOSPA-based approaches\footnote{It was observed that the baseline approach often favours making multiple observations of the same region (particularly if $P_d<1$ and there is clutter), thereby allowing the target to be accurately geo-located if  (by chance) it is within that region and is subsequently detected at least once. However, the lack of exploration makes it difficult to infer whether a target is present if this tactic fails (i.e. the target is not detected).}.

Exemplar suboptimal and optimal actions for $T=2$ and $P_d \in\{0.6,  0.9, 1.0\}$ are shown in Figure \ref{figure_5}, for each of the target distributions and for scenarios with clutter (i.e. $\lambda_{FA}=0.01$) and measurement errors (i.e. $\Sigma=\mbox{diag}(10^{-4}\mbox{km}^2, 10^{-4}\mbox{km}^2)$). In these cases, because of the measurement origin uncertainty, the generation of measurements does not guarantee that a target is present, and the second step actions are always required. For the optimal approach, the second step action is dependent on the measurement(s) generated on the first time step, and the four most commonly selected second step actions are shown in green in Figures \ref{figure_5}(g) -- \ref{figure_5}(l). Generally, in the scenarios considered, the first step action is the same for the optimal and suboptimal approaches. The one exception is shown in  \ref{figure_5}(j), in which the the optimal action offsets the first sensor look and subsequently does not always view the centre of the target distribution.

Overall location RMSE and AMMS-GOSPA values for the three control approaches, for $\lambda_{FA}=0.01$, are shown in Table \ref{table_costs_2}.
It is observed that the baseline approach again results in slightly lower estimation errors, but with significantly greater AMMS-GOSPA errors (and therefore significantly greater cardinality errors).
The suboptimal approach generates AMMS-GOSPA values that are greater than the optimal approach (reaffirming the  derivations (\ref{eq:SM_9}) -- (\ref{opt_cost_function})), with a maximum difference of around 7\%. Future work will identify  scenarios in which performance differences between the optimal and suboptimal approaches are more significant.

%Consequently, the overall costs incurred by the suboptimal approach ($V^\star_{1:T}$) are greater (as proven in (\ref{eq:SM_9}) -- (\ref{new_bellman_full}), and demonstrated in Table \ref{table_costs}), due to the costs incurred by  unnecessary future sensing.

%\begin{eqnarray}
%\mbox{action hypotheses, }A &\triangleq& \left\{a_0,a_1, \ldots, a_n\right\}
%\end{eqnarray}
%where $a_0=\phi$ (i.e. do not attempt a sensor observation), and $a_i$ ($i > 0$) denotes that the sensor is %steered
%to have target hypothesis $x_i$ at the centre of its spotlight.

%\newpage

\subsubsection{Myopic vs. Non-Myopic Planning}
\begin{figure}[H]
\vspace{-0.0cm}
\[
\begin{array}{ccc}%
\SetLabels
\L (0.11*0.825) \small{\textbf{(a):} $T=2$, $P_d=0.6$} \\
\endSetLabels
\psfrag{y}[]{\small{Northings (km)}}
\psfrag{x}[]{\small{Eastings (km)}}
\strut\AffixLabels{\includegraphics[width=4cm]{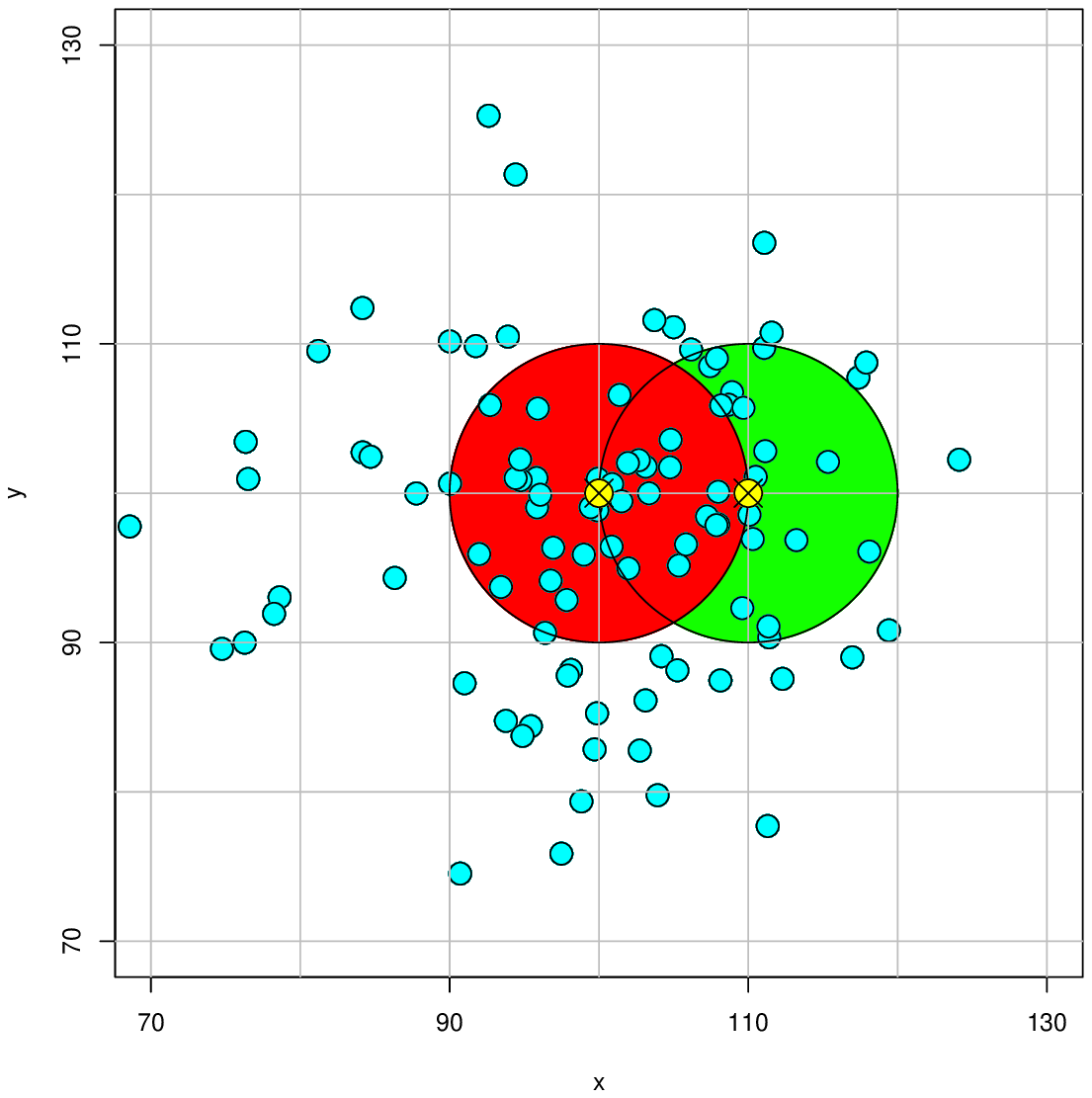}} &
\SetLabels
\L (0.11*0.825) \small{\textbf{(b):} $T=2$, $P_d=\{0.6,0.9\}$} \\
\endSetLabels
\psfrag{y}[]{\small{Northings (km)}}
\psfrag{x}[]{\small{Eastings (km)}}\hspace{0.5cm}
\strut\AffixLabels{\includegraphics[width=4cm]{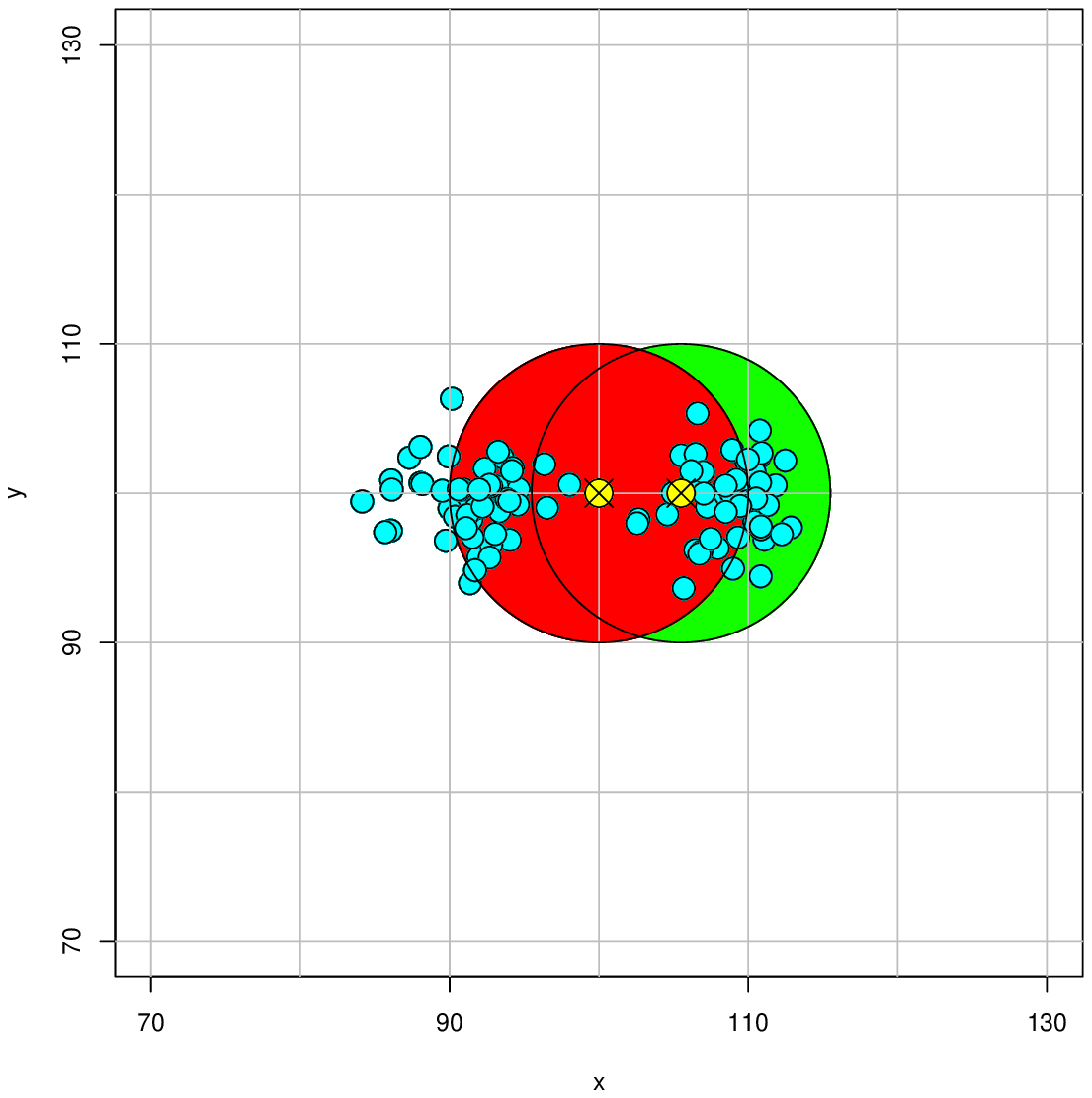}} &
\SetLabels
\L (0.11*0.825) \small{\textbf{(c):} $T=2$, $P_d=0.6$} \\
\endSetLabels
\psfrag{y}[]{\small{Northings (km)}}
\psfrag{x}[]{\small{Eastings (km)}}
\vspace{-1cm}\hspace{0.5cm}
\strut\AffixLabels{\includegraphics[width=4cm]{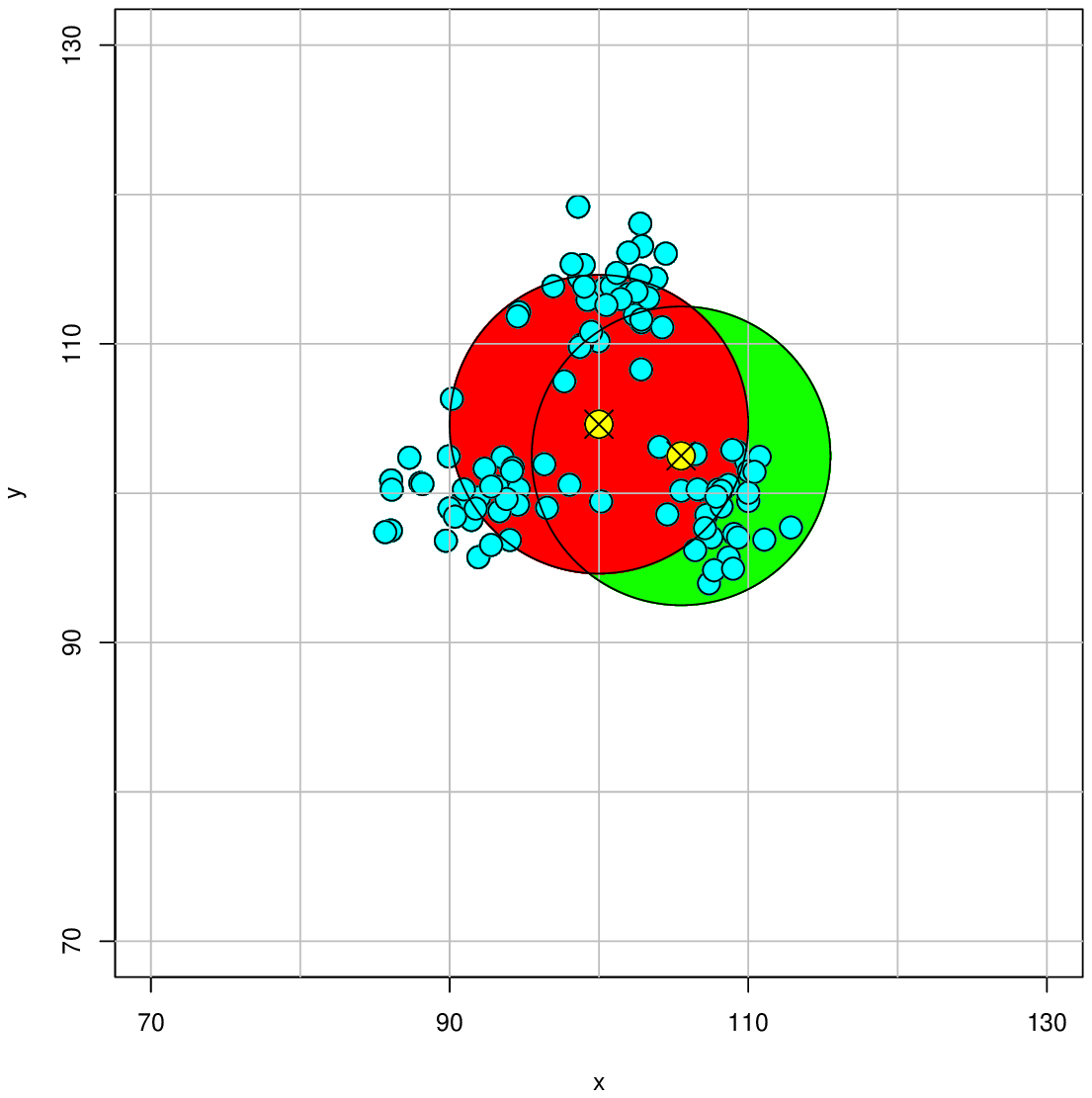}} \\
\SetLabels
\L (0.11*0.825) \small{\textbf{(d):} $T=2$, $P_d=\{0.9, 1.0\}$} \\
\endSetLabels
\psfrag{y}[]{\small{Northings (km)}}
\psfrag{x}[]{\small{Eastings (km)}}
\strut\AffixLabels{\includegraphics[width=4cm]{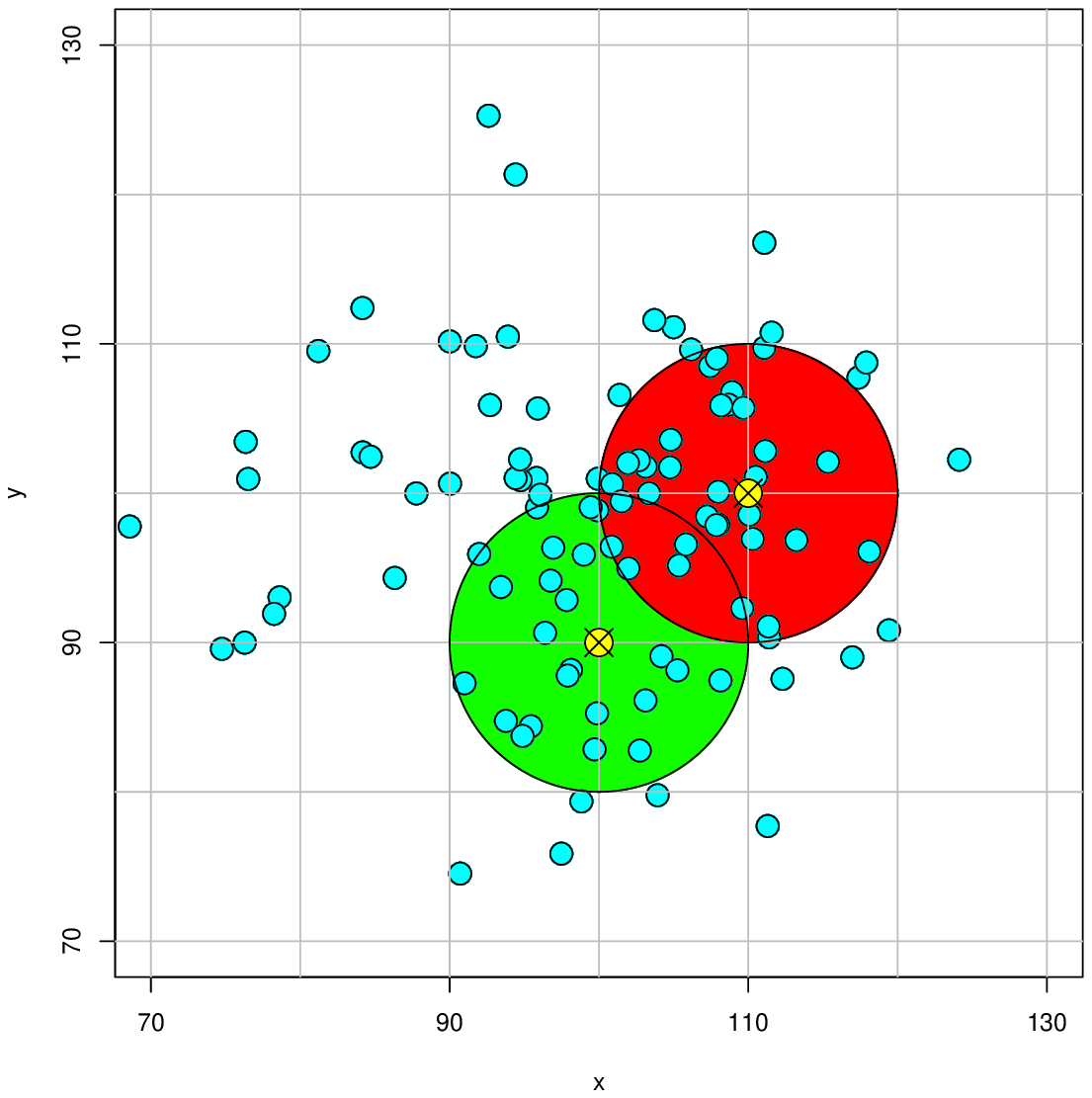}} &
\SetLabels
\L (0.11*0.825) \small{\textbf{(e):} $T=2$, $P_d= 1.0$} \\
\endSetLabels
\psfrag{y}[]{\small{Northings (km)}}
\psfrag{x}[]{\small{Eastings (km)}}\hspace{0.5cm}
\strut\AffixLabels{\includegraphics[width=4cm]{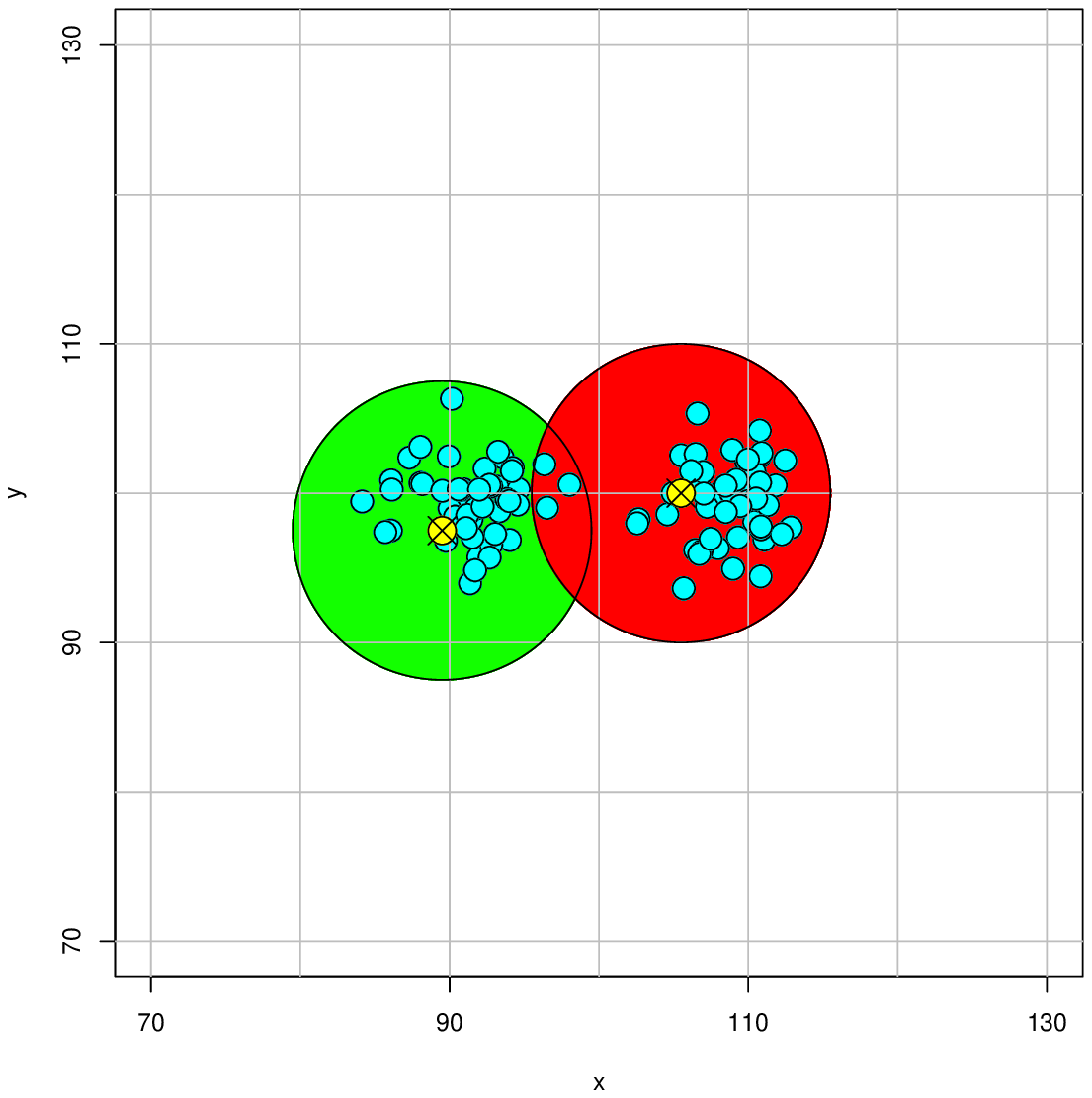}} &
\SetLabels
\L (0.11*0.825) \small{\textbf{(f):} $T=2$, $P_d=\{0.9, 1.0\}$} \\
\endSetLabels
\psfrag{y}[]{\small{Northings (km)}}
\psfrag{x}[]{\small{Eastings (km)}}
\vspace{-1cm}\hspace{0.5cm}
\strut\AffixLabels{\includegraphics[width=4cm]{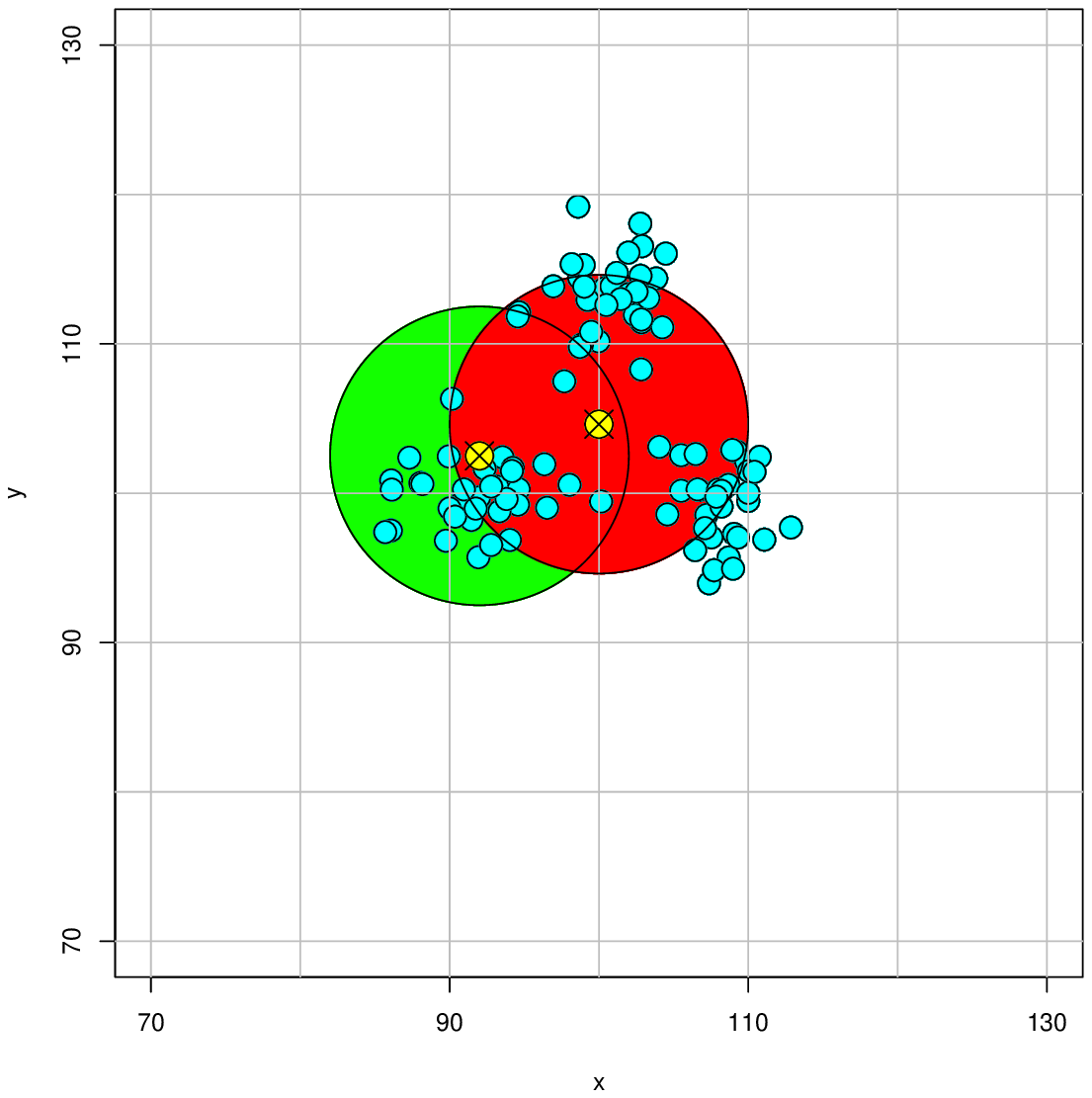}} \\
\SetLabels
\L (0.11*0.825) \small{\textbf{(g):} $T=3$, $P_d=0.6$} \\
\endSetLabels
\psfrag{y}[]{\small{Northings (km)}}
\psfrag{x}[]{\small{Eastings (km)}}
\strut\AffixLabels{\includegraphics[width=4cm]{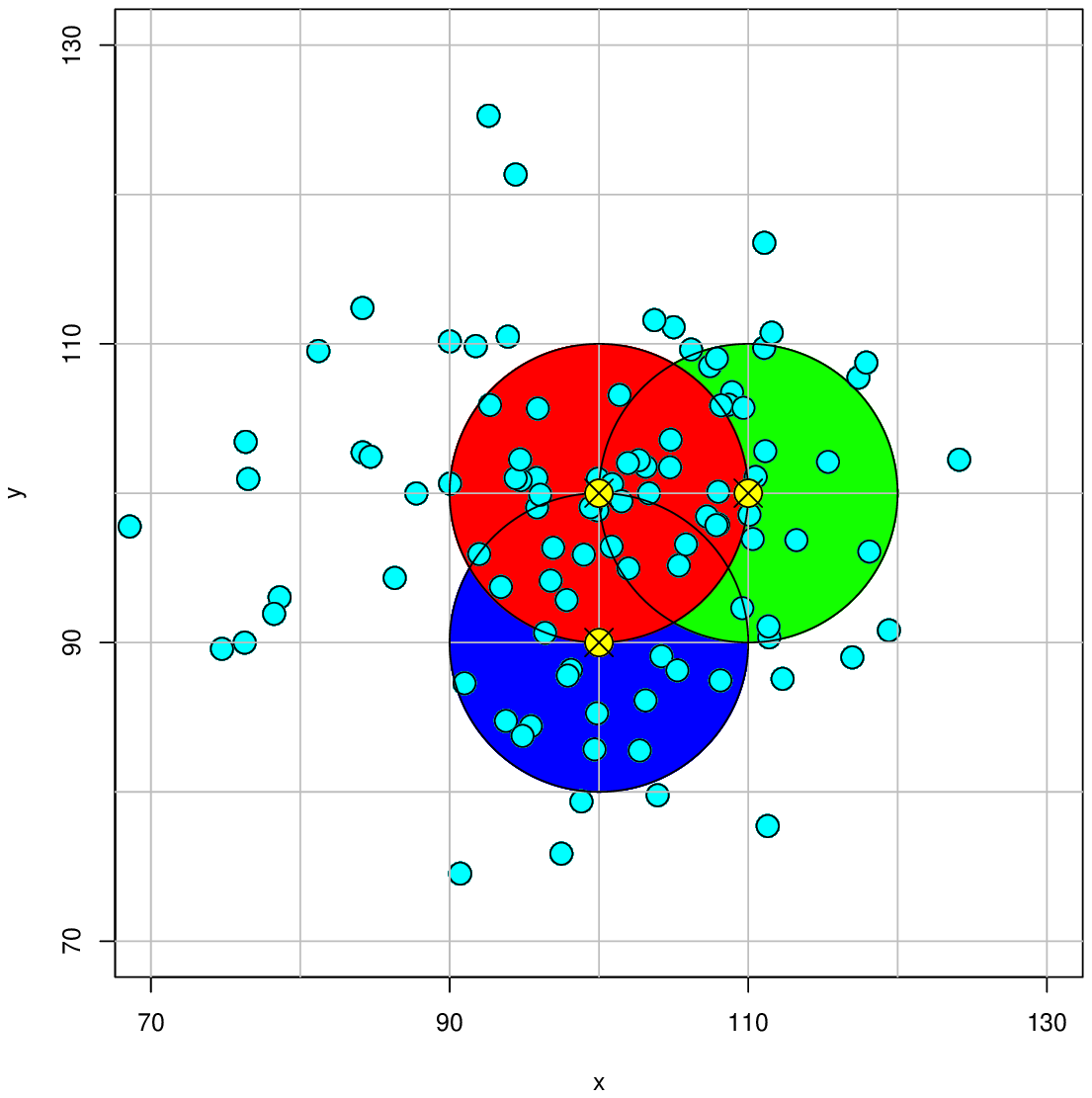}} &
\SetLabels
\L (0.11*0.825) \small{\textbf{(h):} $T=3$, $P_d=\{0.6, 0.9\}$} \\
\endSetLabels
\psfrag{y}[]{\small{Northings (km)}}
\psfrag{x}[]{\small{Eastings (km)}}\hspace{0.5cm}
\strut\AffixLabels{\includegraphics[width=4cm]{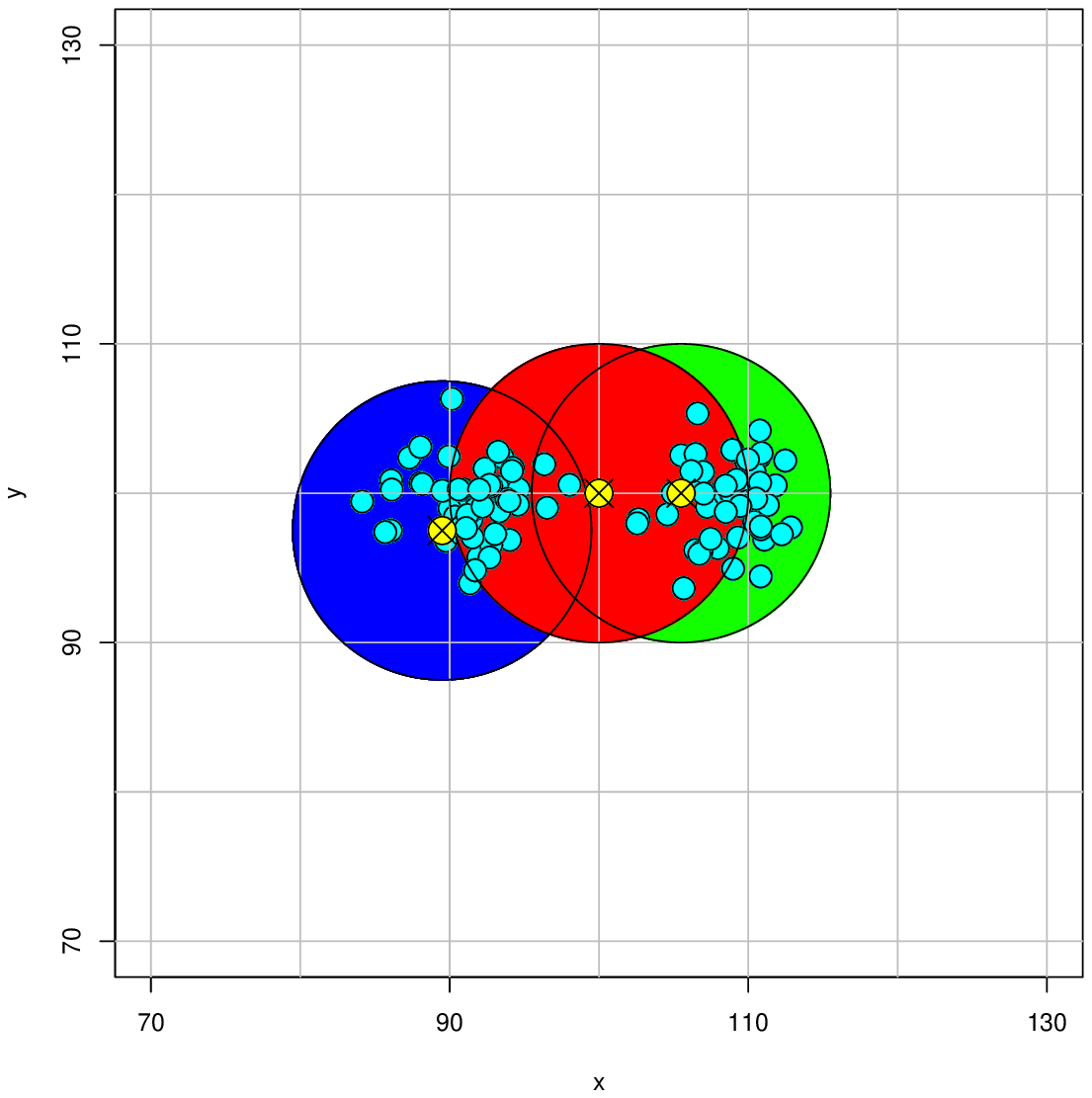}} &
\SetLabels
\L (0.11*0.825) \small{\textbf{(i):} $T=3$, $P_d=\{0.6,0.9\}$} \\
\endSetLabels
\psfrag{y}[]{\small{Northings (km)}}
\psfrag{x}[]{\small{Eastings (km)}}
\vspace{-1cm}\hspace{0.5cm}
\strut\AffixLabels{\includegraphics[width=4cm]{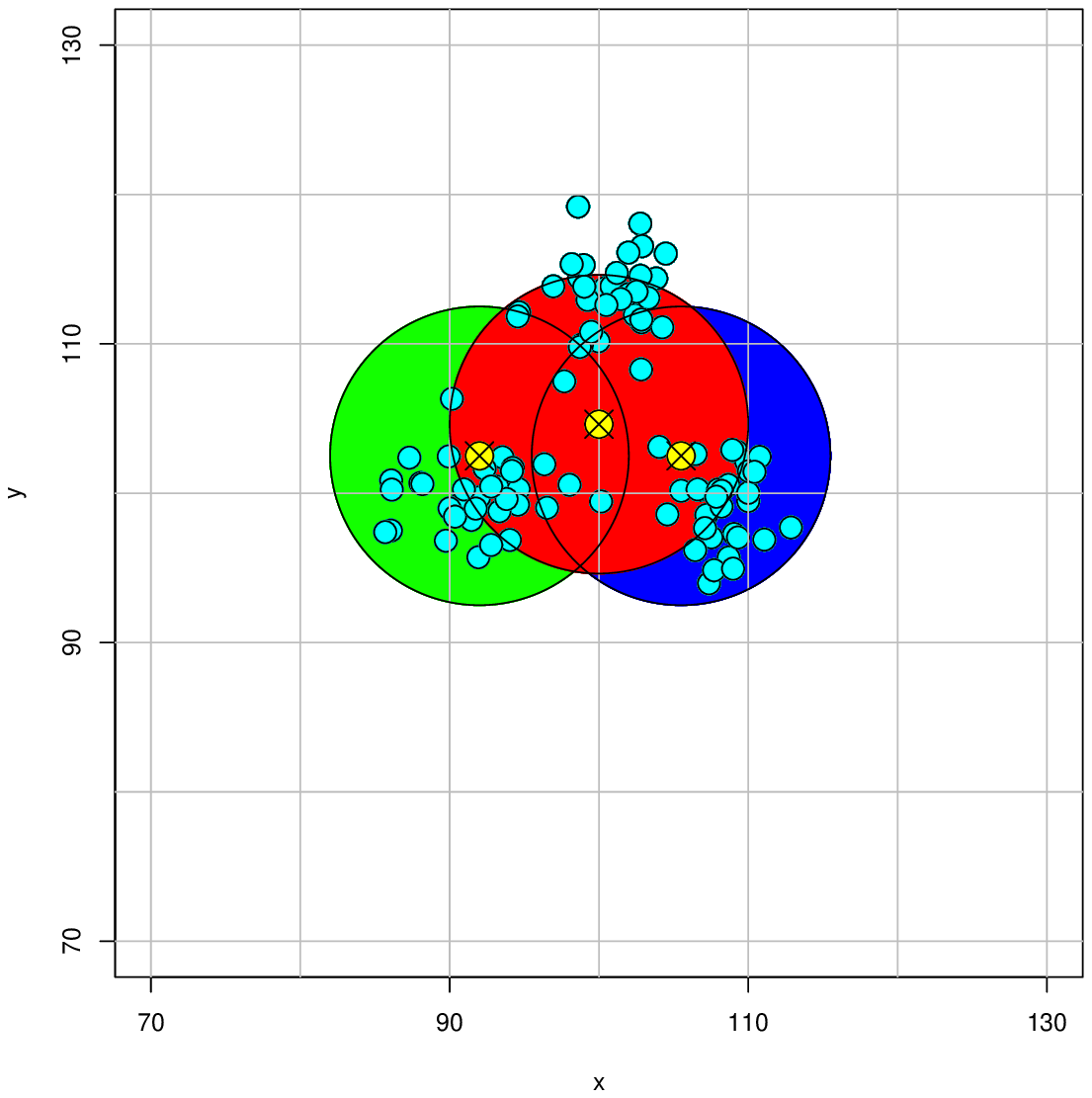}} \\
\SetLabels
\L (0.11*0.825) \small{\textbf{(j):} $T=3$, $P_d=\{0.9, 1.0\}$} \\
\endSetLabels
\psfrag{y}[]{\small{Northings (km)}}
\psfrag{x}[]{\small{Eastings (km)}}
\strut\AffixLabels{\includegraphics[width=4cm]{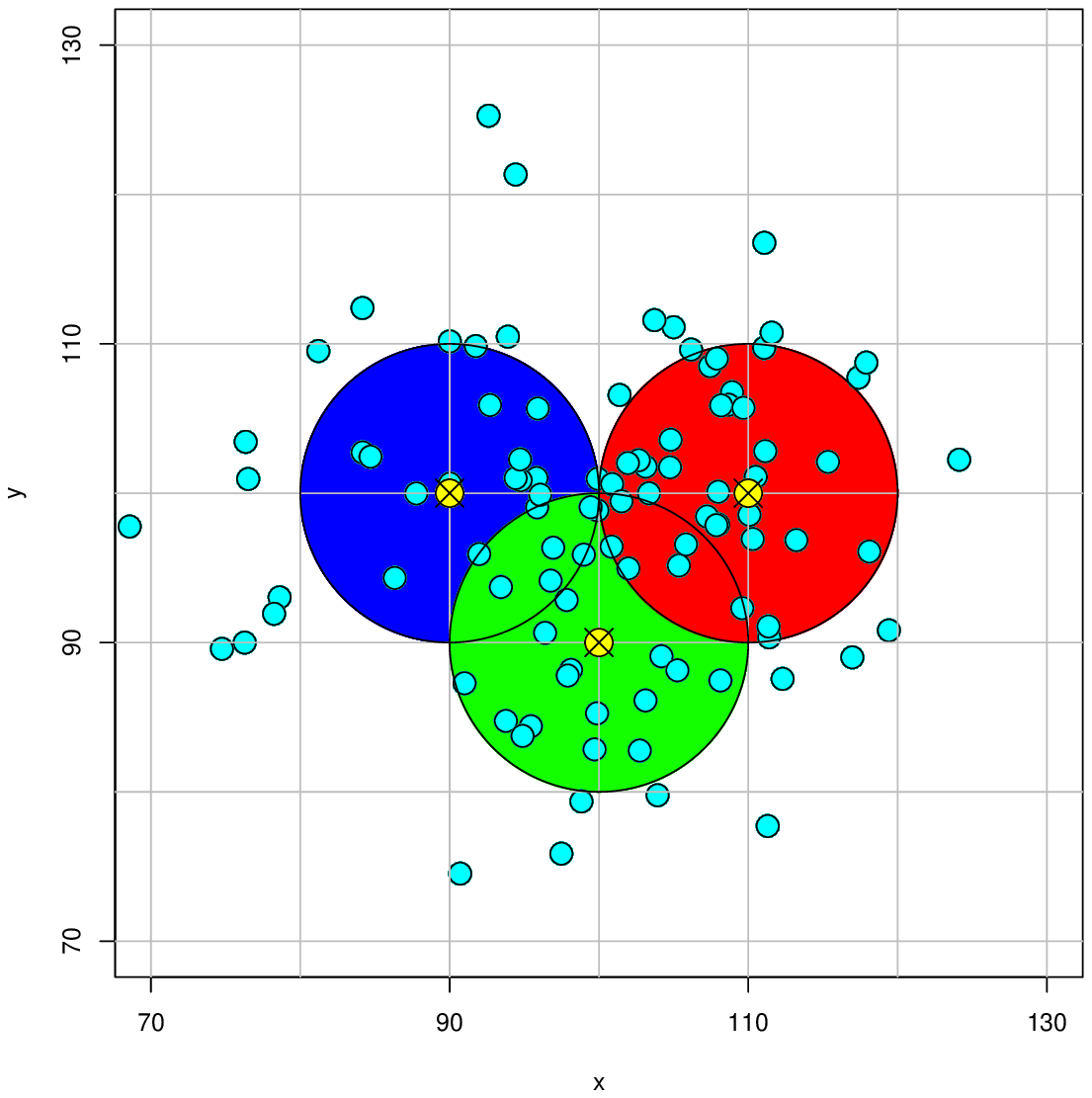}} &
\SetLabels
\L (0.11*0.825) \small{\textbf{(k):} $T=3$, $P_d=1.0$} \\
\endSetLabels
\psfrag{y}[]{\small{Northings (km)}}
\psfrag{x}[]{\small{Eastings (km)}}\hspace{0.5cm}
\strut\AffixLabels{\includegraphics[width=4cm]{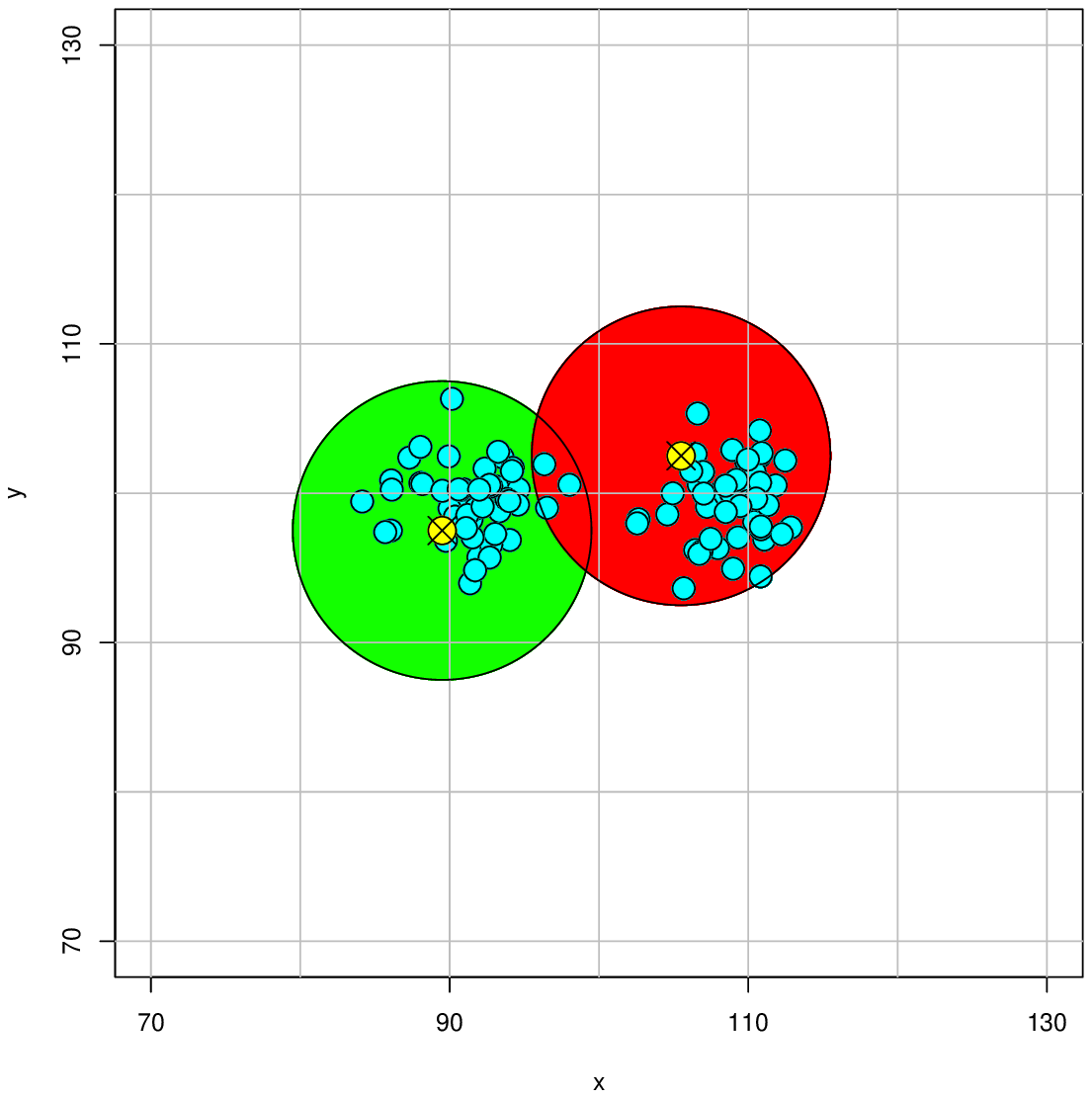}} &
\SetLabels
\L (0.11*0.825) \small{\textbf{(l):} $T=3$, $P_d=1.0$} \\
\endSetLabels
\psfrag{y}[]{\small{Northings (km)}}
\psfrag{x}[]{\small{Eastings (km)}}
\vspace{-1cm}\hspace{0.5cm}
\strut\AffixLabels{\includegraphics[width=4cm]{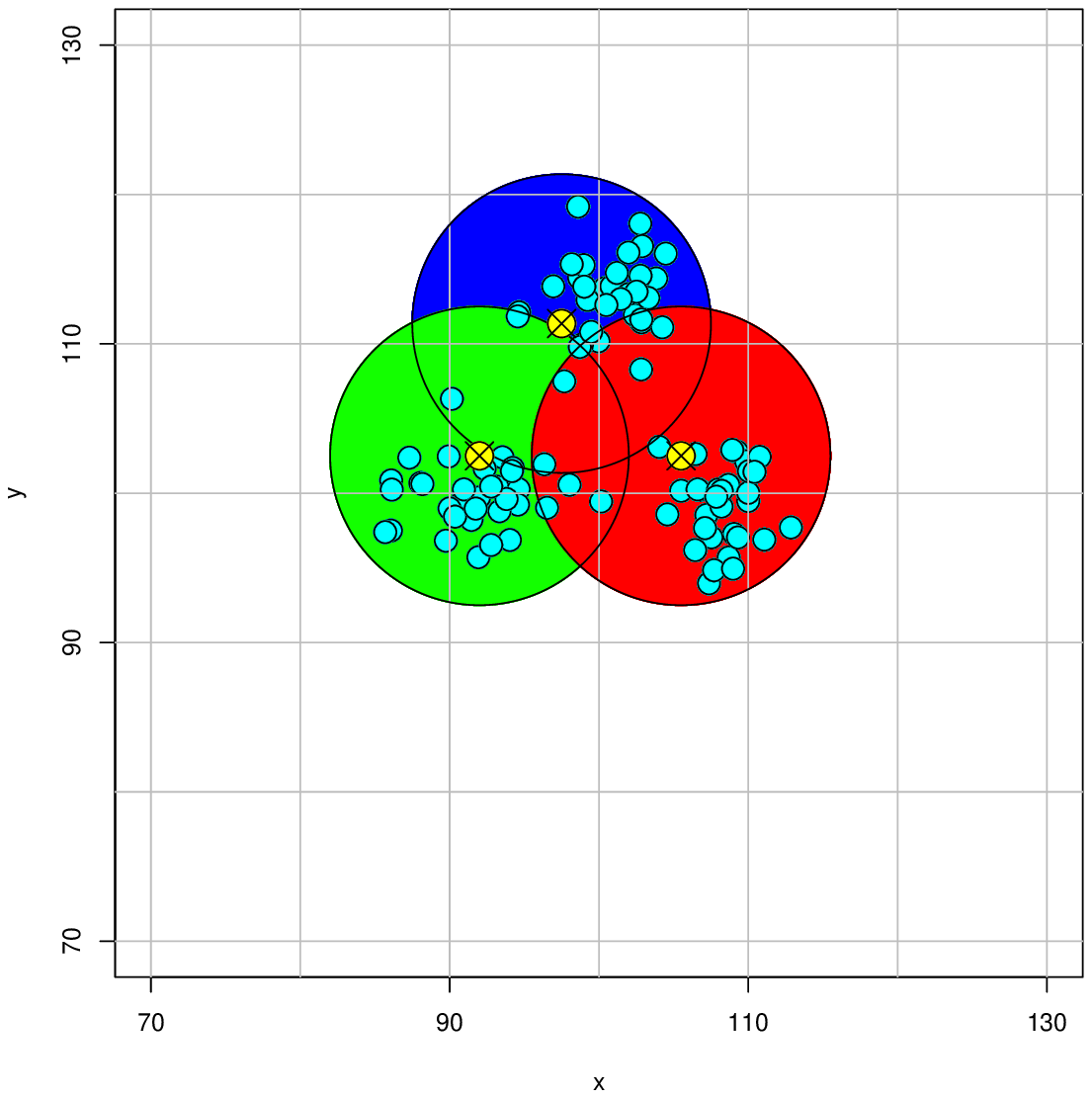}} \\
\vspace{0.3cm}
\end{array}
\]
\caption{Optimal actions  for exemplar scenarios with no clutter (i.e. $\lambda_{FA}=0.0$), $P_d \in\{0.6,  0.9, 1.0\}$, and non-myopic time horizons $T\in\{2,3\}$. Left column: unimodal target prior distribution, middle column: bimodal distribution, right column: trimodal distribution. Red circles: FOV of 1st action,
green circles: FOV of 2nd action, blue circles: FOV of 3rd
action. For the optimal approach, the 2nd and 3rd actions are only necessary if a measurement is not generated on any previous timestep. For the suboptimal control approach, the 2nd and 3rd actions occur regardless of whether previous measurements have been generated. The optimal myopic action is always to observe the centre of the distribution (i.e. with a sensor spotlight centre as shown by the red circles in (a) -- (c)).}
\label{figure_4}
\end{figure}

\begin{figure}[H]
\vspace{-0.0cm}
\[
\begin{array}{ccc}%
\SetLabels
\L (0.11*0.825) \small{\textbf{(a):} subopt, $P_d=\{0.6,0.9\}$} \\
\endSetLabels
\psfrag{y}[]{\small{Northings (km)}}
\psfrag{x}[]{\small{Eastings (km)}}
\strut\AffixLabels{\includegraphics[width=4cm]{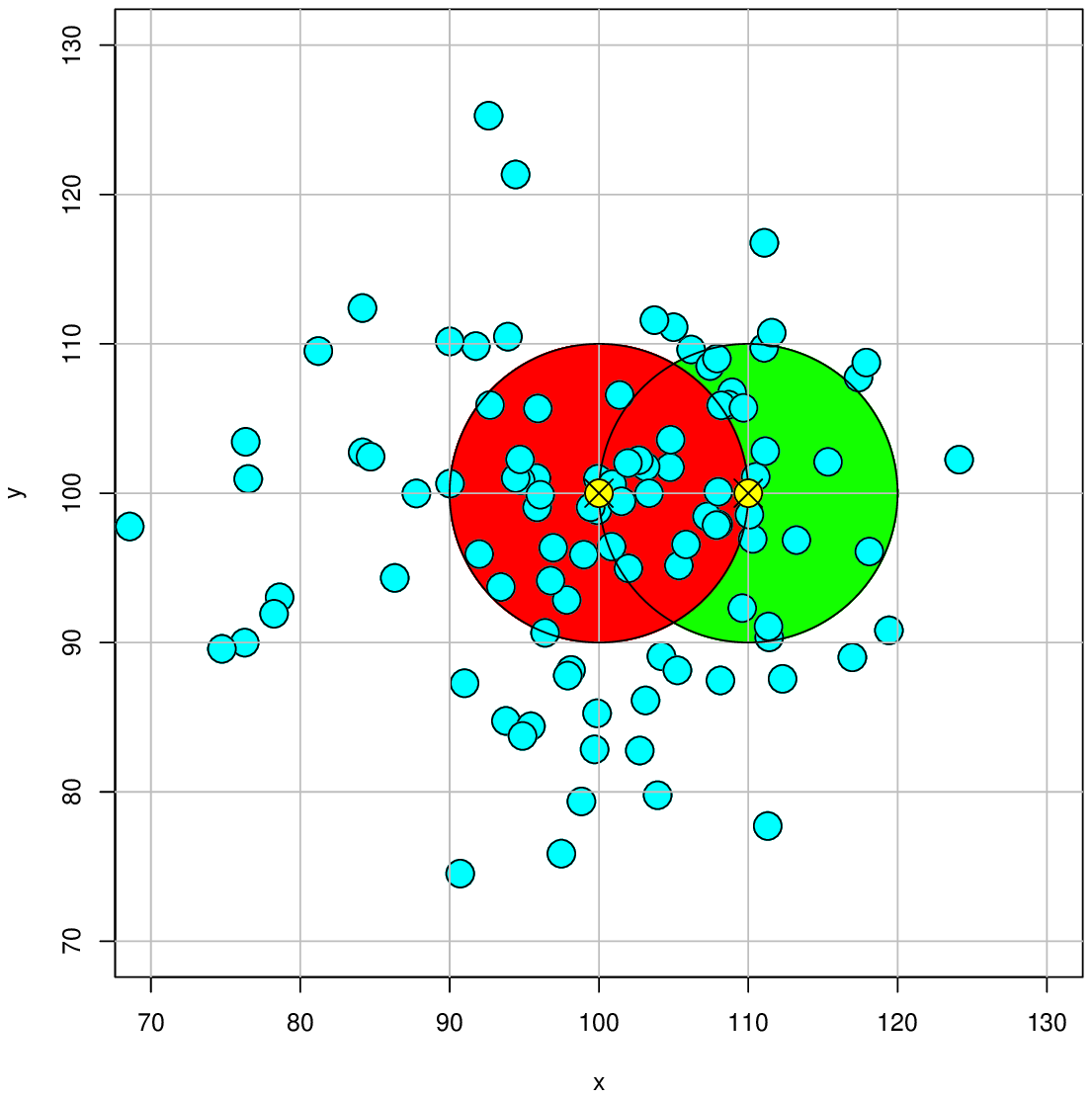}} &
\SetLabels
\L (0.11*0.825) \small{\textbf{(b):} subopt, $P_d=\{0.6,0.9\}$} \\
\endSetLabels
\psfrag{y}[]{\small{Northings (km)}}
\psfrag{x}[]{\small{Eastings (km)}}\hspace{0.5cm}
\strut\AffixLabels{\includegraphics[width=4cm]{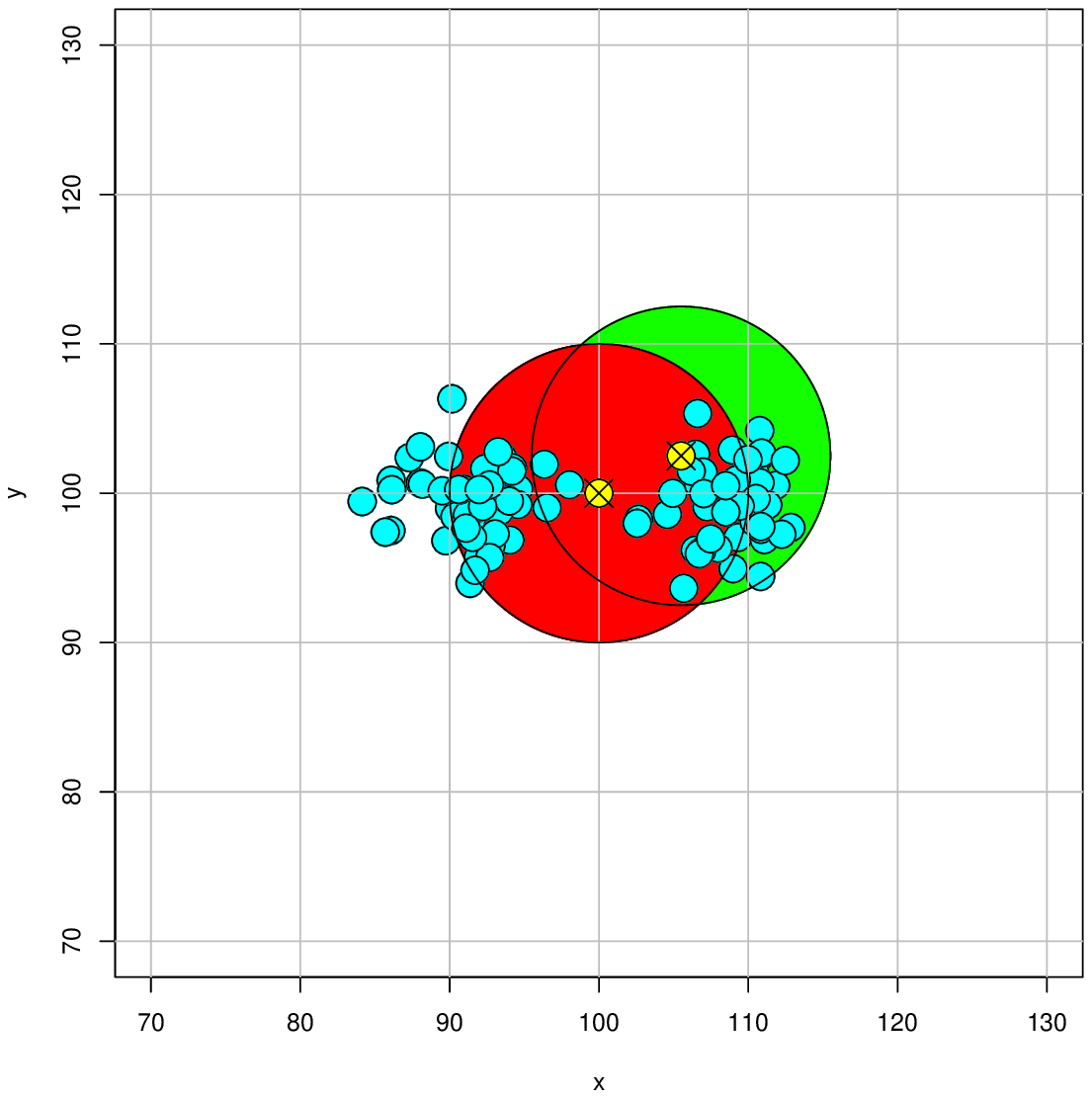}} &
\SetLabels
\L (0.11*0.825) \small{\textbf{(c):} subopt, $P_d=0.6^a$} \\
\endSetLabels
\psfrag{y}[]{\small{Northings (km)}}
\psfrag{x}[]{\small{Eastings (km)}}
\vspace{-1cm}\hspace{0.5cm}
\strut\AffixLabels{\includegraphics[width=4cm]{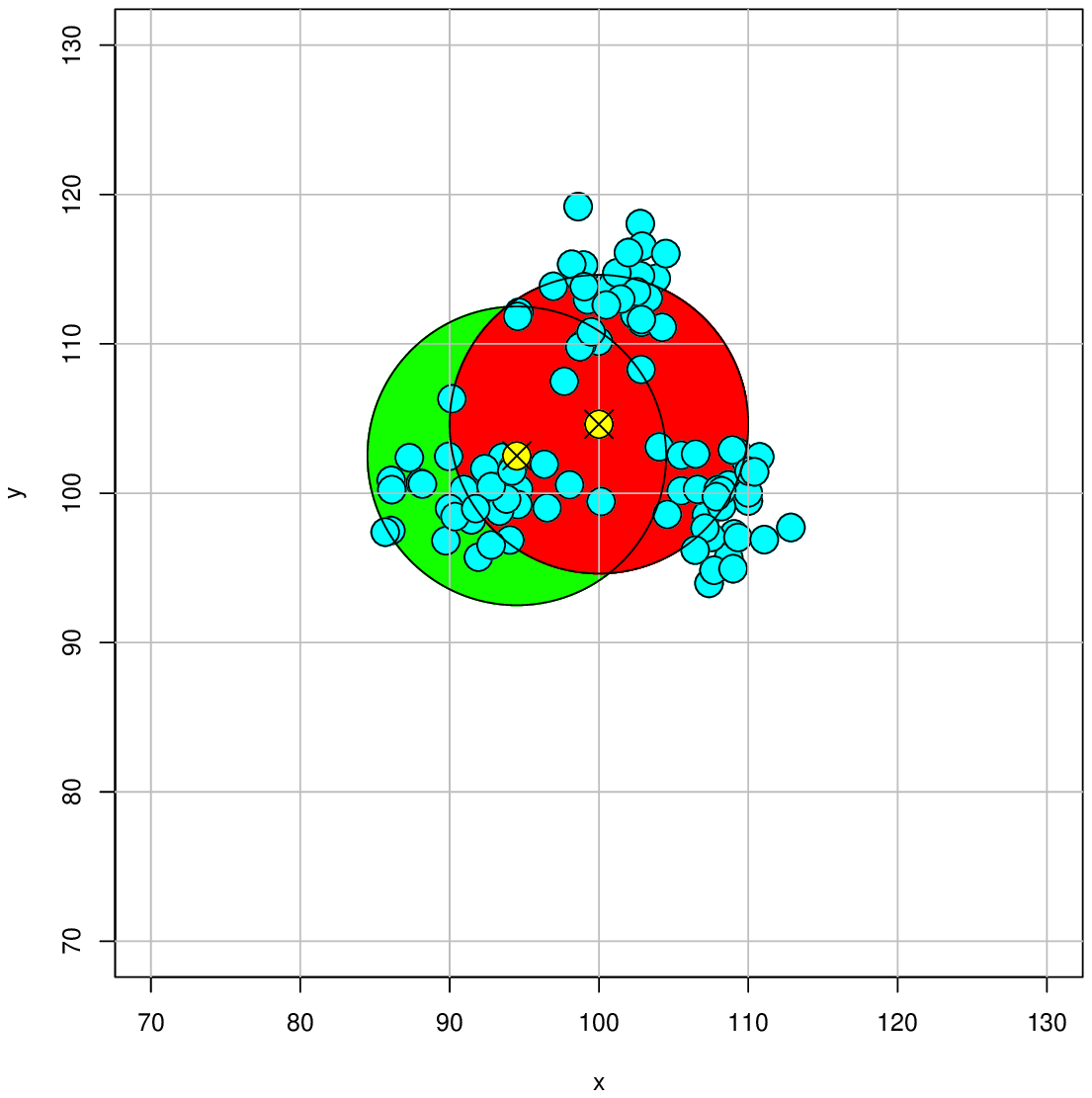}} \\
\SetLabels
\L (0.11*0.825) \small{\textbf{(d):} subopt, $P_d= 1.0$} \\
\endSetLabels
\psfrag{y}[]{\small{Northings (km)}}
\psfrag{x}[]{\small{Eastings (km)}}
\strut\AffixLabels{\includegraphics[width=4cm]{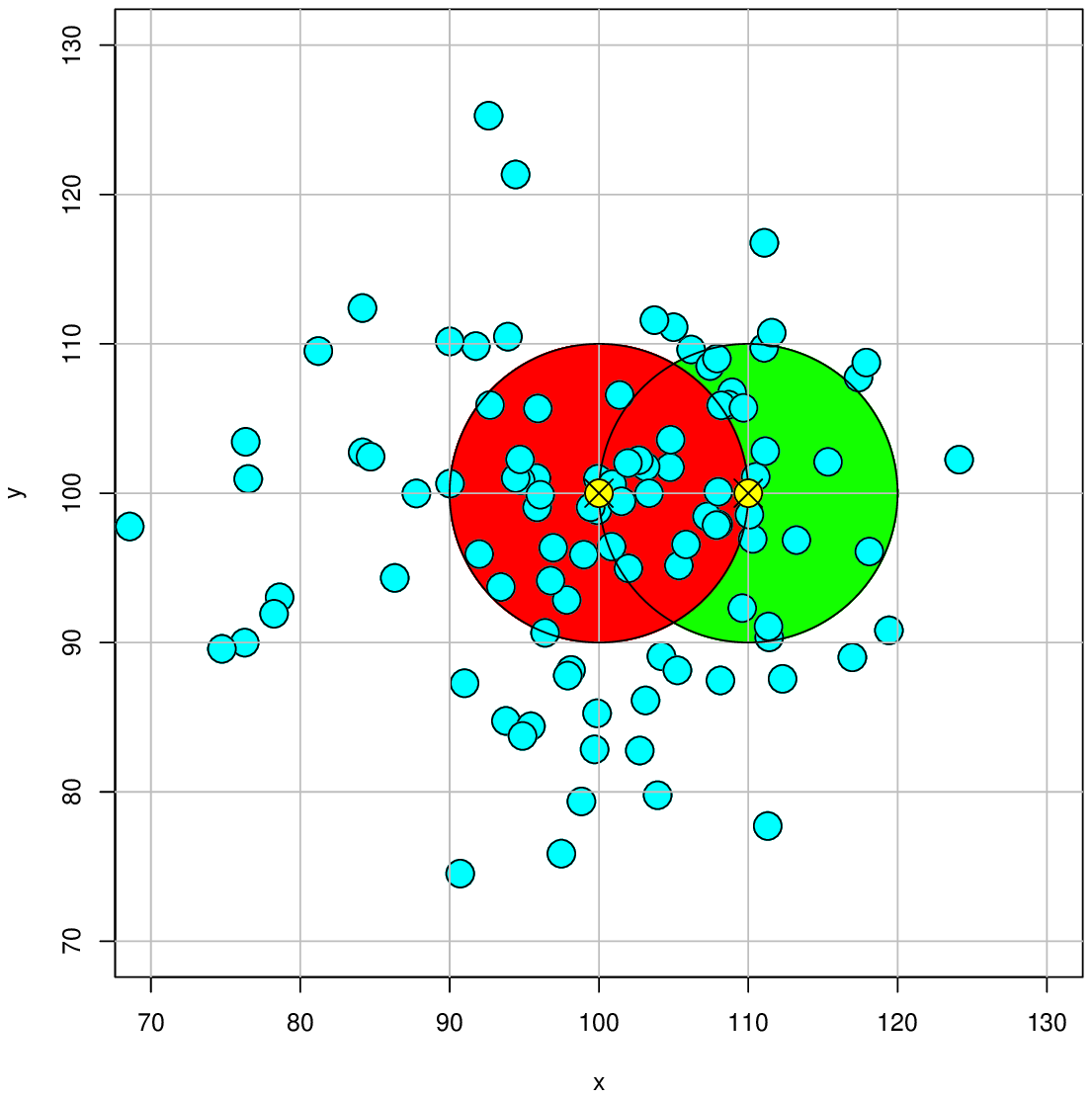}} &
\SetLabels
\L (0.11*0.825) \small{\textbf{(e):} subopt, $P_d= 1.0$} \\
\endSetLabels
\psfrag{y}[]{\small{Northings (km)}}
\psfrag{x}[]{\small{Eastings (km)}}\hspace{0.5cm}
\strut\AffixLabels{\includegraphics[width=4cm]{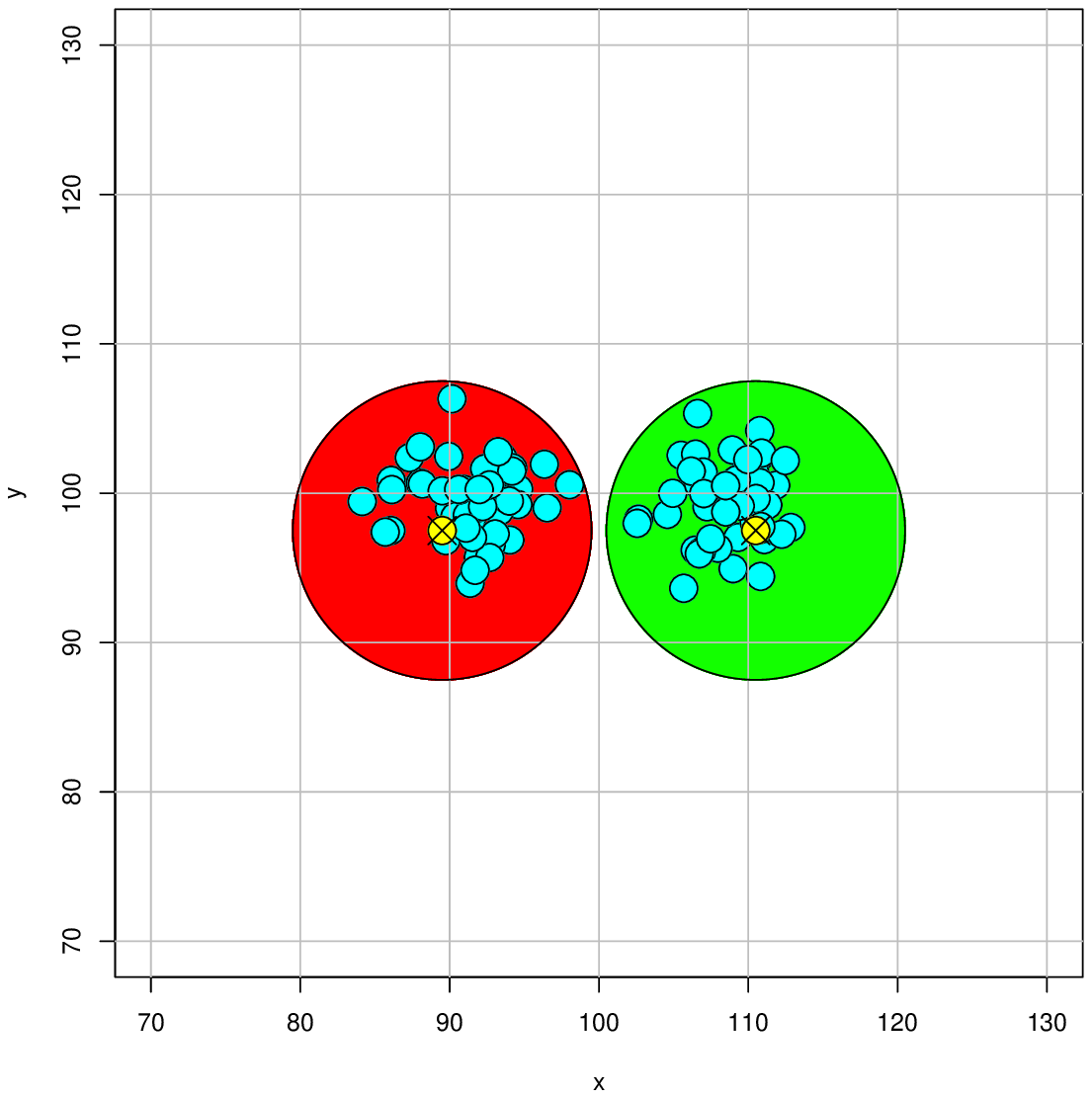}} &
\SetLabels
\L (0.11*0.825) \small{\textbf{(f):} subopt, $P_d=1.0$} \\
\endSetLabels
\psfrag{y}[]{\small{Northings (km)}}
\psfrag{x}[]{\small{Eastings (km)}}
\vspace{-1cm}\hspace{0.5cm}
\strut\AffixLabels{\includegraphics[width=4cm]{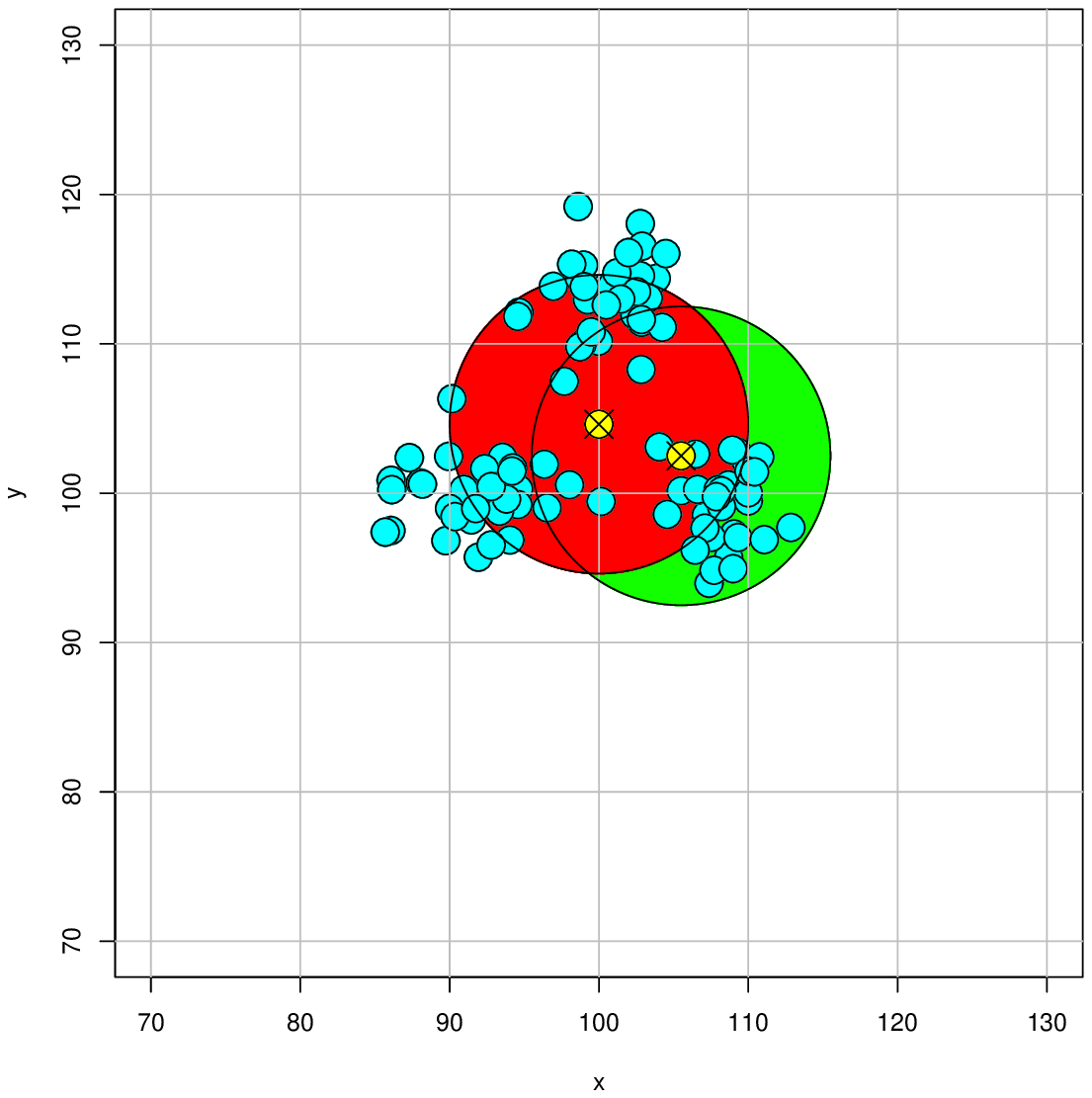}} \\
\SetLabels
\L (0.11*0.825) \small{\textbf{(g):} optimal, $P_d=\{0.6,0.9\}$} \\
\endSetLabels
\psfrag{y}[]{\small{Northings (km)}}
\psfrag{x}[]{\small{Eastings (km)}}
\strut\AffixLabels{\includegraphics[width=4cm]{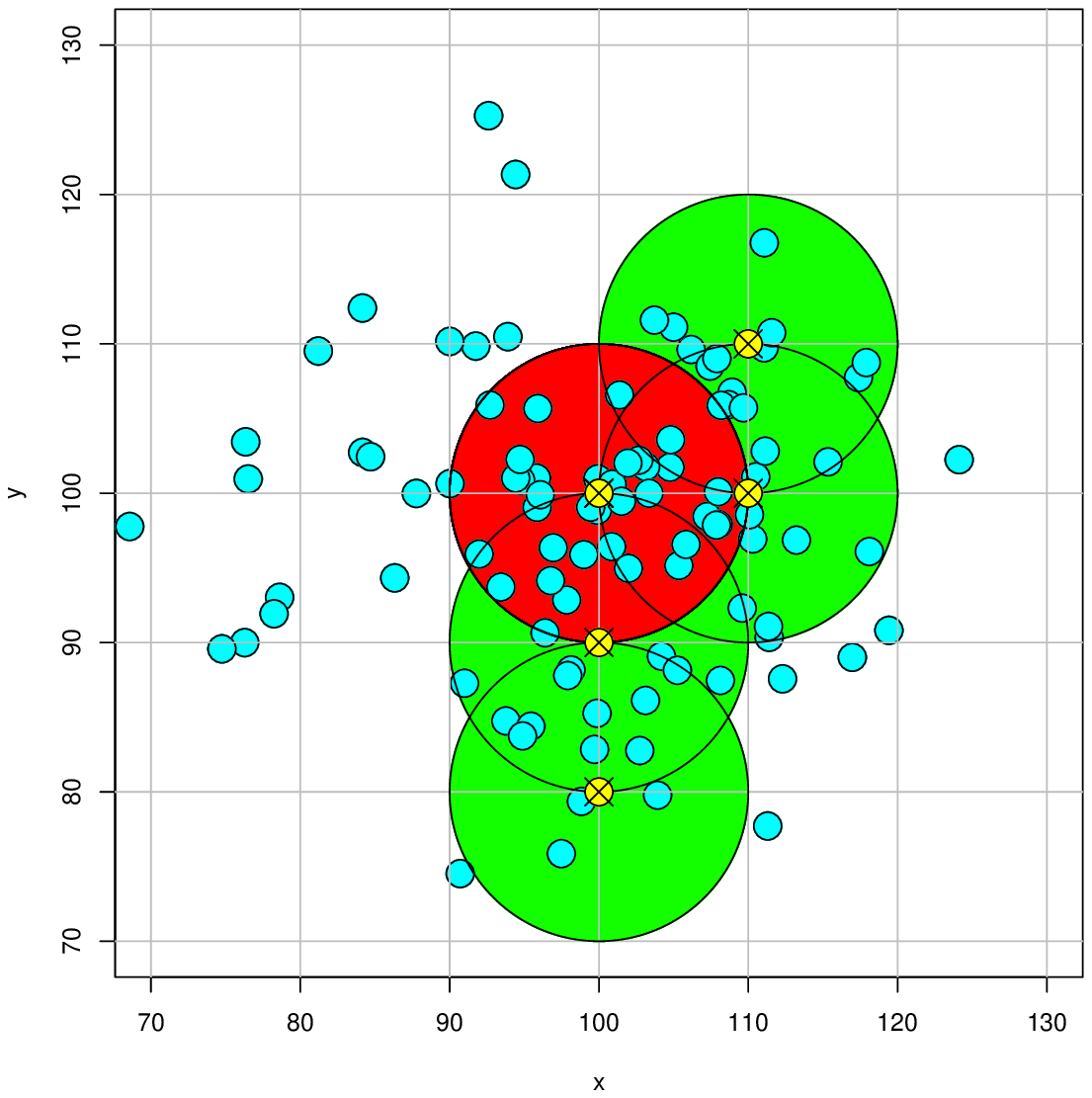}} &
\SetLabels
\L (0.11*0.825) \small{\textbf{(h):} optimal, $P_d=0.6^b$} \\
\endSetLabels
\psfrag{y}[]{\small{Northings (km)}}
\psfrag{x}[]{\small{Eastings (km)}}\hspace{0.5cm}
\strut\AffixLabels{\includegraphics[width=4cm]{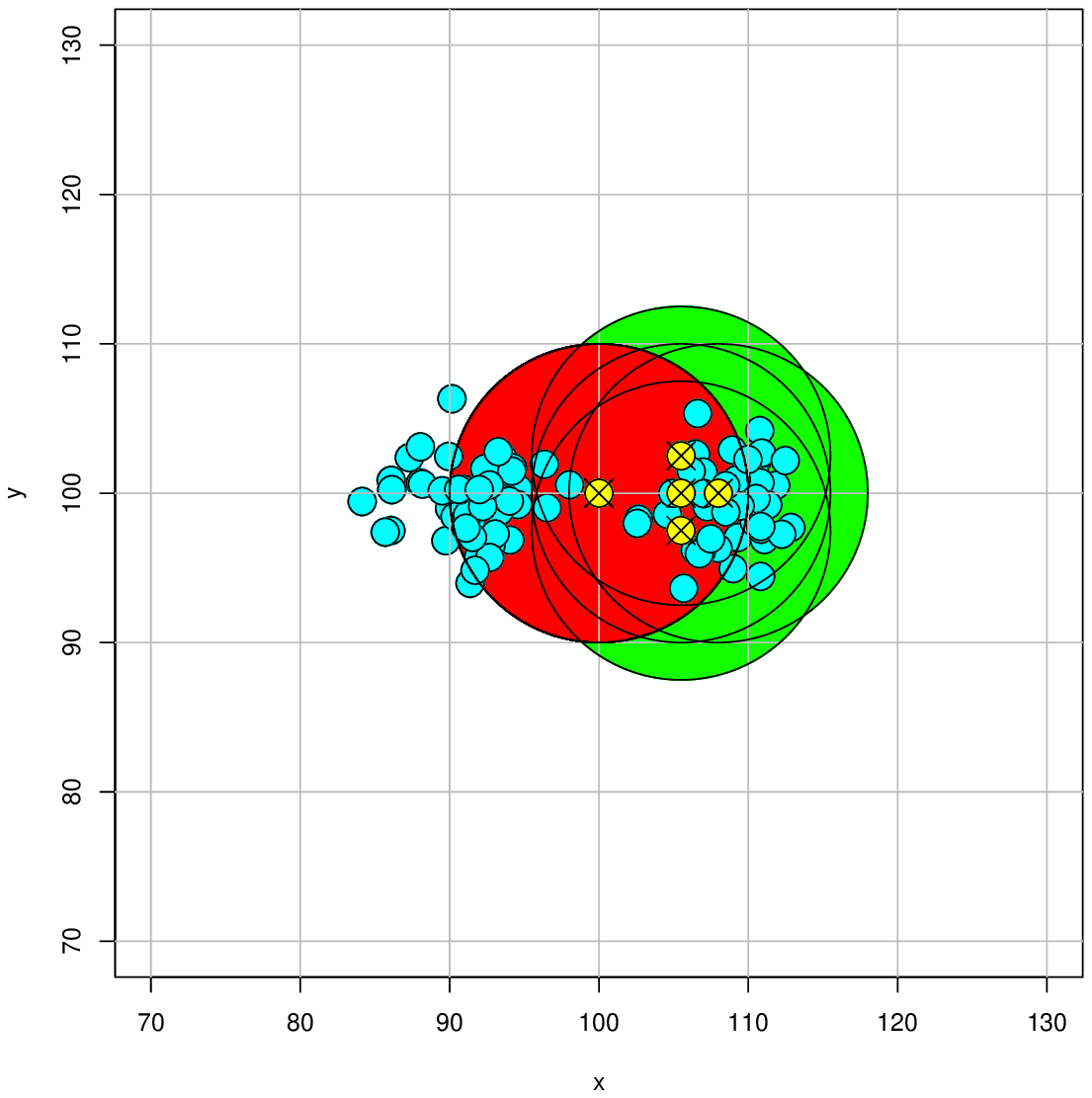}} &
\SetLabels
\L (0.11*0.825) \small{\textbf{(i):} optimal, $P_d=0.6^b$} \\
\endSetLabels
\psfrag{y}[]{\small{Northings (km)}}
\psfrag{x}[]{\small{Eastings (km)}}
\vspace{-1cm}\hspace{0.5cm}
\strut\AffixLabels{\includegraphics[width=4cm]{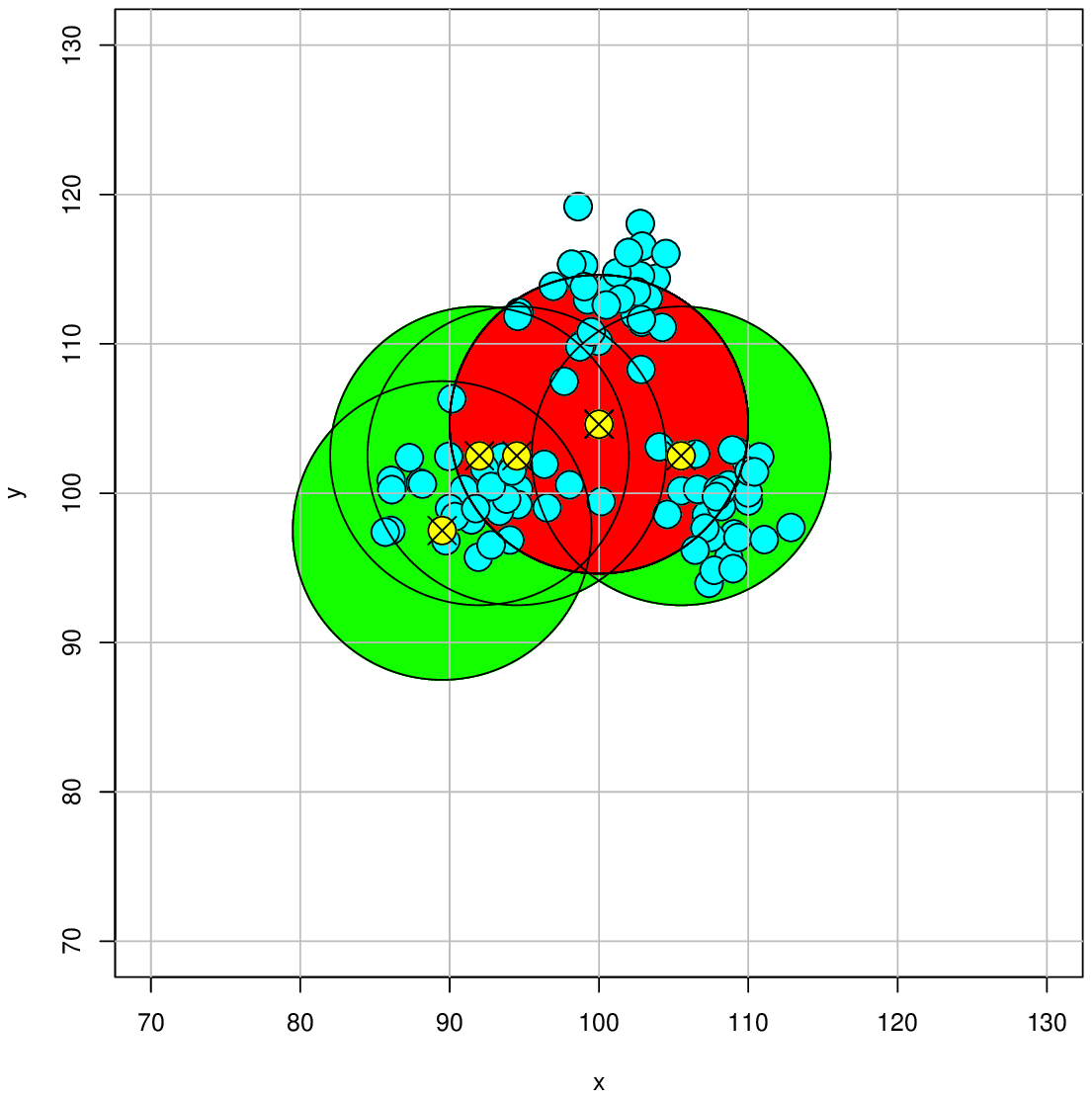}} \\
\SetLabels
\L (0.11*0.825) \small{\textbf{(j):} optimal, $P_d=1.0$} \\
\endSetLabels
\psfrag{y}[]{\small{Northings (km)}}
\psfrag{x}[]{\small{Eastings (km)}}
\strut\AffixLabels{\includegraphics[width=4cm]{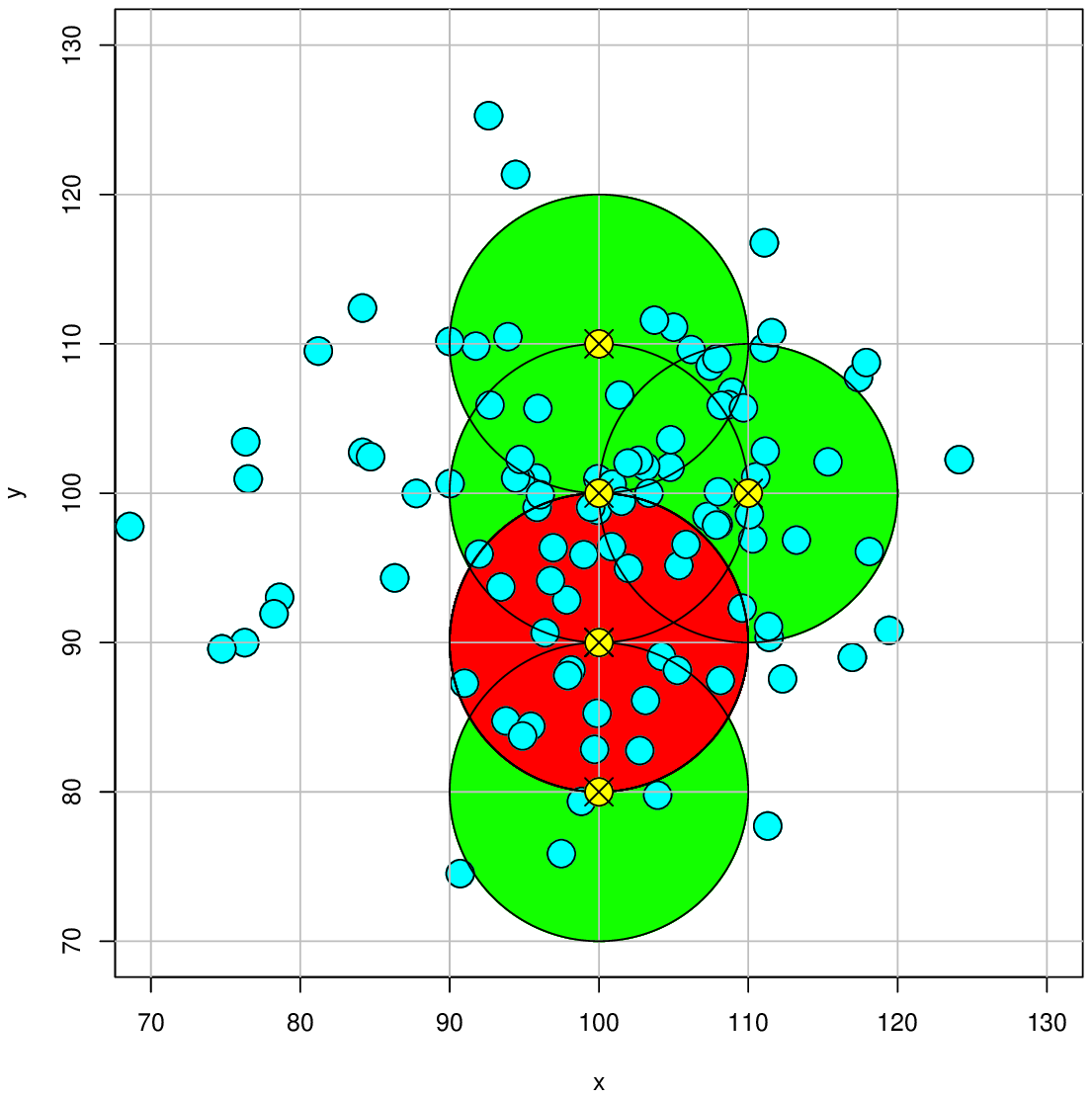}} &
\SetLabels
\L (0.11*0.825) \small{\textbf{(k):} optimal, $P_d=1.0$} \\
\endSetLabels
\psfrag{y}[]{\small{Northings (km)}}
\psfrag{x}[]{\small{Eastings (km)}}\hspace{0.5cm}
\strut\AffixLabels{\includegraphics[width=4cm]{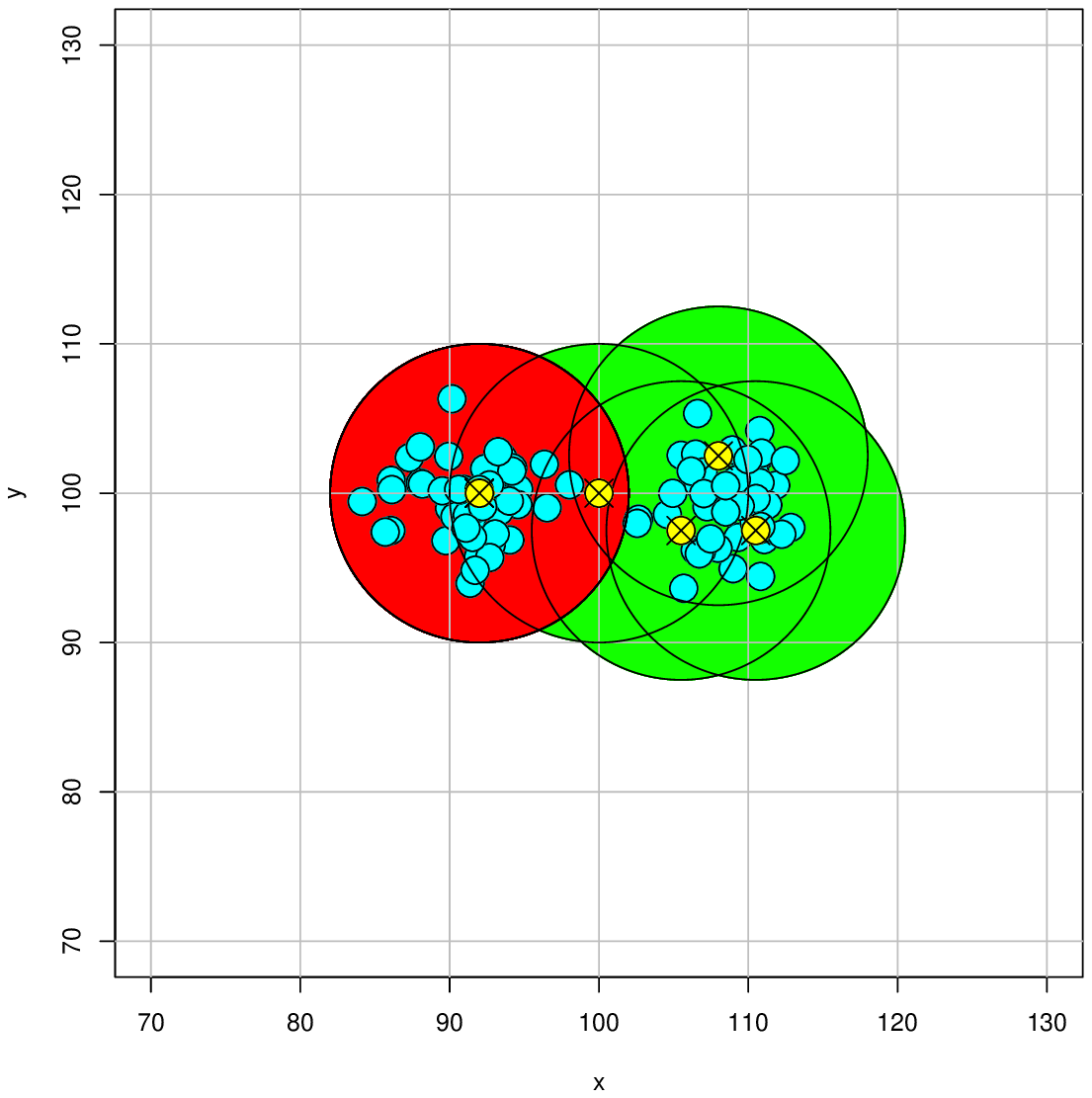}} &
\SetLabels
\L (0.11*0.825) \small{\textbf{(l):} optimal, $P_d=1.0$} \\
\endSetLabels
\psfrag{y}[]{\small{Northings (km)}}
\psfrag{x}[]{\small{Eastings (km)}}
\vspace{-1cm}\hspace{0.5cm}
\strut\AffixLabels{\includegraphics[width=4cm]{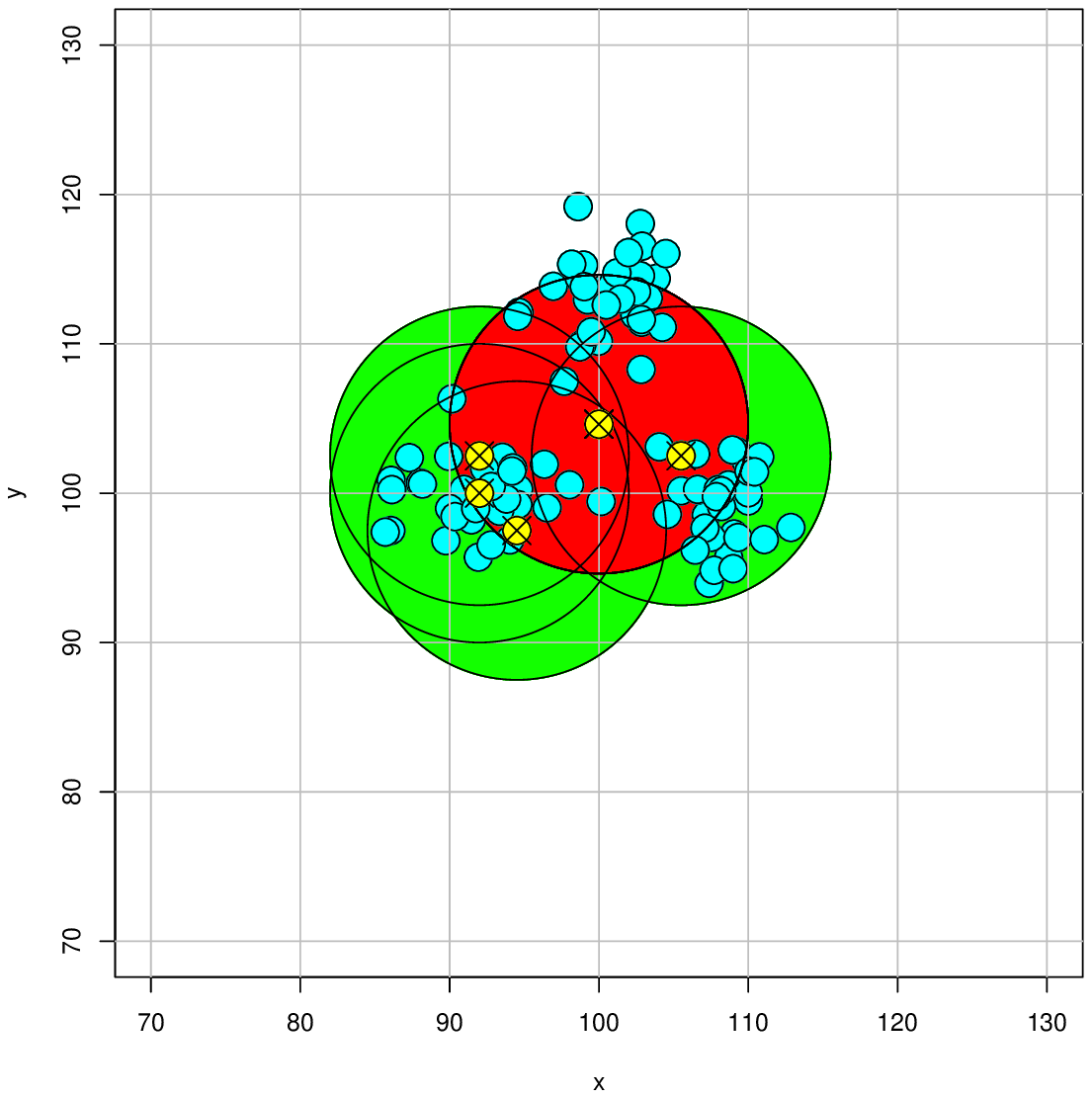}} \\
\vspace{0.3cm}
\end{array}
\]
\caption{Suboptimal and optimal actions for exemplar scenarios with clutter (i.e. $\lambda_{FA}=0.01$), $P_d \in\{0.6, 0.9, 1.0\}$, and  $T=2$. Left column: unimodal target prior distribution, middle column: bimodal distribution, right column: trimodal distribution. Red circles: FOV of 1st action,
green circles: FOV of 2nd action (for the optimal solution, the four most commonly occurring 2nd step actions are shown). Again, the optimal myopic action is always to observe the centre of the distribution (i.e. with a sensor spotlight centre as shown by the red circles in (a) -- (c)). \newline
$^a$The solution for $P_d=0.9$ has a slightly different 2nd action which also prioritises attempting to detect a target in the left most mode.
\newline
$^b$The solution for $P_d=0.9$ is virtually identical to that for $P_d=0.6$, but with a very slight difference to the most commonly occuring 2nd step actions.}
\label{figure_5}
\end{figure}

In each scenario, the optimal myopic strategy maximises the probability of
detecting the target via a single measurement, thereby always viewing the centre/midpoint of the target hypotheses (i.e. shown by the red circle in Figures \ref{figure_4}(a) -- \ref{figure_4}(c) and Figures \ref{figure_5}(a) -- \ref{figure_5}(c)). However, myopic planning does not have the foresight to appreciate that further measurements can be generated. Clearly, observing the centre point can be suboptimal in the multi-modal scenarios, e.g. for the bimodal distribution, it may then be necessary to make (at least) two further observations (one for each mode) to guarantee that the target is detected.

        For the non-myopic  approaches, it can be seen from Figures \ref{figure_4} -- \ref{figure_5} that:
\begin{enumerate}
\vspace{0.15cm}
\item When $P_d=\{0.6,0.9\}$, multi-step planning also almost always favours initially viewing the centre/midpoint of the target
    hypotheses (except in the scenarios shown in Figures \ref{figure_4}(d) and \ref{figure_4}(j)). Consequently, there is considerable overlap between the action spotlights on each time
step. This is because when the probability of detection is less than unity, multiple sensor observations in regions of high target probability mass provides an improved
opportunity to generate at least one detection of  a target.
    \vspace{0.15cm}
 \item   When $P_d=1.0$,  multi-step planning almost always avoids viewing the centre of the distribution of target
     hypotheses\footnote{Notable exceptions are shown in Figures \ref{figure_4}(f), \ref{figure_5}(f), and \ref{figure_5}(l), in which two-step planning cannot provide coverage of all three modes, and consequently prioritises attempting early detection by viewing the centre of the target distribution.}, and instead offsets each sensor spotlight to achieve
 optimal surveillance coverage over
     multiple time steps by observing each mode in turn, e.g. see Figures \ref{figure_4}(e),  \ref{figure_4}(k) and \ref{figure_4}(l)).
\vspace{0.15cm}
\item
When $P_d=1.0$ and $\lambda_{FA}=0$, one observation  is sufficient to
detect the target in any given region. Therefore, in these cases, a  ``cookie-cutter'' strategy with limited overlap between the sensor spotlights is optimal.
\end{enumerate}

\vspace{0.15cm}
\noindent
These  observations are  in agreement with an intuitive understanding of how taking into
account the ability to make  further observations should influence the first observation.

%\subsubsection{Computational Complexity}
%The run-time of the planning approaches is shown in Table \ref{table_runtime}. It can be seen that the run-time
%increases by circa two orders of magnitude for each additional step in the planning horizon. A one order of magnitude increase is as a result of each of (i): the increased complexity of the optimisation  (\ref{optimal_Bellman_act}), and: (ii): the increased computational expense of calculating the conditional AMMS-GOSPA over the longer time window.
%
%It is noted that the
%run-times are significantly lower when $P_d=1.0$ because the   algorithm  is then able to exploit the fact that if a
%target is not observed within a sensor spotlight, the target hypotheses within the spotlight can be removed from
%consideration on subsequent time steps. This effect is most evident when considering a bimodal target prior
%distribution.
%It is also noted that the run-time of the optimal control approach  increases as the number of action  hypotheses increases. Hence in the unimodal and trimodal scenarios (with 26 and 29 action hypotheses respectively), the run-times are greater than for the bimodal scenarios (which have 20 hypotheses).

\begin{table}[H]
\vspace{0.25cm}
\caption{Overall location RMSE and AMMM-GOSPA costs incurred (in km) by the baseline and optimal control approaches, for each of the scenarios considered with no clutter (i.e. $\lambda_{FA}=0$) and $P_d \in \{0.6, 0.9, 1.0\}$. Results are averaged over 20 runs, with the mean value $\pm$ one standard deviation shown. To remind the reader, the overall AMMS-GOSPA error is identical in the suboptimal and optimal approaches and given by $V^\star_{1:T}$ (equation (\ref{simon_10})) or  $\hat{V}^\star_{1:T}$ (equation (\ref{new_bellman_full})).}
\begin{center}
%\vspace{-0.5cm}
\begin{tabular}{|c|c|c||c|c|c|c|} \hline
  &     &  \textbf{ Control } & \multicolumn{2}{c|}{$\quad$\textbf{$T=2$}$\quad$}  &  \multicolumn{2}{c|}{$\quad$\textbf{$T=3$}$\quad$} \\
  \cline{4-7}
\raisebox{1.25ex}[0cm][0cm]{\textbf{ Target Prior }} &   \raisebox{1.25ex}[0cm][0cm]{$\quad P_d \quad$}  &\textbf{ Approach } &
\multicolumn{1}{c|}{$\qquad$\textbf{RMSE$\qquad$}} &  \multicolumn{1}{c|}{\textbf{AMMS-GOSPA}} &
\multicolumn{1}{c|}{$\qquad$\textbf{RMSE$\qquad$}} &  \multicolumn{1}{c|}{\textbf{AMMS-GOSPA}}
\\\hline\hline
&      & Baseline & 25.58 $\pm$ 1.01 & 65.53 $\pm$ 3.99 & 37.21 $\pm$ 1.57 & 92.57 $\pm$ 5.94
 \\ \cline{3-7}
 & \raisebox{1.25ex}[0cm][0cm]{0.6}    & Optimal &  26.72 $\pm$ 1.10
& 55.90 $\pm$ 1.85 & 38.99 $\pm$ 1.67 & 77.20 $\pm$ 2.65 \\ \cline{2-7}
&   & Baseline   & 23.77 $\pm$ 1.08  & 57.84 $\pm$ 5.85 & 33.61 $\pm$ 1.74  &  78.76 $\pm$ 7.70 \\ \cline{3-7}
\raisebox{1.25ex}[0cm][0cm]{Unimodal} & \raisebox{1.25ex}[0cm][0cm]{0.9}    & Optimal &
25.53 $\pm$ 1.14 & 45.53 $\pm$ 2.42 & 36.43 $\pm$ 1.96 & 60.00 $\pm$ 3.03 \\ \cline{2-7}
&     & Baseline & 23.05 $\pm$ 1.12 & 55.62 $\pm$ 6.24 & 32.09 $\pm$ 1.87 & 74.52 $\pm$ 7.66  \\ \cline{3-7}
 & \raisebox{1.25ex}[0cm][0cm]{1.0} & Optimal    &  25.17 $\pm$ 1.24  & 42.18 $\pm$ 2.57 & 35.41 $\pm$ 1.93
 & 54.56 $\pm$ 3.09
 \\ \hline\hline
 &      & Baseline &  11.64 $\pm$  0.30 &  49.43  $\pm$ 1.20 &  15.29 $\pm$ 0.40 &  69.55 $\pm$ 1.98
 \\ \cline{3-7}
 & \raisebox{1.25ex}[0cm][0cm]{0.6}    & Optimal &  12.87 $\pm$ 0.46
&  36.36 $\pm$ 2.20  &  17.36 $\pm$  0.66 &  45.14 $\pm$  2.43 \\ \cline{2-7}
&   & Baseline   & $\mbox{ \  }$6.70 $\pm$ 0.20  &  39.52 $\pm$ 1.31 &  $\mbox{ \  }$8.02 $\pm$ 0.25  &  41.59  $\pm$ 1.40 \\ \cline{3-7}
\raisebox{1.25ex}[0cm][0cm]{Bimodal} & \raisebox{1.25ex}[0cm][0cm]{0.9}    & Optimal &
$\mbox{ \  }$9.02 $\pm$ 0.66 &  19.29 $\pm$ 2.34  &  10.94 $\pm$ 0.81 &  20.56 $\pm$ 2.50 \\ \cline{2-7}
&     & Baseline &  $\mbox{ \  }$2.20 $\pm$  0.12 &  14.03 $\pm$ 0.52 &  $\mbox{ \  }$2.20 $\pm$  0.12 &  14.18 $\pm$  0.65 \\ \cline{3-7}
 & \raisebox{1.25ex}[0cm][0cm]{1.0} & Optimal    &   $\mbox{ \  }$4.01 $\pm$ 1.98   & 13.03 $\pm$ 1.39  & $\mbox{ \  }$4.01 $\pm$ 1.98
 & 13.03 $\pm$ 1.39
 \\ \hline\hline
 &      & Baseline &  15.68 $\pm$ 0.46  & 52.01 $\pm$ 3.21 &  21.85 $\pm$ 0.66 &  68.18 $\pm$ 4.88
 \\ \cline{3-7}
 & \raisebox{1.25ex}[0cm][0cm]{0.6}    & Optimal & 16.40  $\pm$ 0.70
& 46.57 $\pm$  2.57 &  22.79 $\pm$ 1.08 &  61.36 $\pm$  3.45 \\ \cline{2-7}
&   & Baseline   &  11.62 $\pm$ 0.36  &  40.30 $\pm$ 1.44  & 14.69  $\pm$  0.45 &  43.76 $\pm$  1.42 \\ \cline{3-7}
\raisebox{1.25ex}[0cm][0cm]{Trimodal} & \raisebox{1.25ex}[0cm][0cm]{0.9}    & Optimal &
14.13  $\pm$ 0.93  &  32.51 $\pm$ 3.38 & 16.77 $\pm$ 1.20 &  38.83 $\pm$  3.44 \\ \cline{2-7}
&     & Baseline &  \mbox{ \ }8.36 $\pm$ 0.38 &  34.91 $\pm$ 1.86 &  \mbox{ \ }8.36 $\pm$ 0.38  &  35.00 $\pm$  1.84 \\ \cline{3-7}
 & \raisebox{1.25ex}[0cm][0cm]{1.0} & Optimal    &  13.28  $\pm$  0.97  & 28.12  $\pm$  0.97 & 11.47 $\pm$ 2.42
 & 31.80 $\pm$ 3.06
 \\ \hline
\end{tabular}
\end{center}
%\vspace{0.2cm}
\label{table_costs}
\end{table}

\newpage

\begin{table}[H]
\caption{Overall errors (i.e. location RMSE and AMMS-GOSPA) incurred (in km) by the baseline, suboptimal and optimal control approaches. There is clutter (i.e. $\lambda_{FA}=0.01$) and $T=2$.
%To remind the reader, $V^\star_{1:T}$ is given by equation (\ref{simon_10}) and  $\hat{V}^\star_{1:T}$ is given by equation (\ref{new_bellman_full}).
Again, results are averaged over 20 runs, with the mean value $\pm$ one standard deviation shown.}
\begin{center}
%\vspace{-0.5cm}
%\hspace{-0.5cm}
\begin{tabular}{|c|c||c|c||c|c||c|c|} \hline
  &    &  \multicolumn{2}{c||}{$\quad$\textbf{Baseline Approach}$\quad$}  &  \multicolumn{2}{c||}{$\quad$\textbf{Suboptimal Control}$\quad$} &
  \multicolumn{2}{c|}{$\quad$\textbf{Optimal Control}$\quad$} \\
  \cline{3-8}
\raisebox{1.25ex}[0cm][0cm]{$\textbf{Target prior}$} &  \raisebox{1.25ex}[0cm][0cm]{\textbf{$P_d$}}  &
\multicolumn{1}{c|}{$\quad\mbox{\textbf{RMSE}}\quad$} & \mbox{\textbf{AMMS-GOSPA}} &
\multicolumn{1}{c|}{$\quad\mbox{\textbf{RMSE}}\quad$} & \mbox{\textbf{AMMS-GOSPA}} &
\multicolumn{1}{c|}{$\quad\mbox{\textbf{RMSE}}\quad$} & \mbox{\textbf{AMMS-GOSPA}}
\\\hline\hline
&   0.6   &  25.59  $\pm$ 1.01  &  74.65  $\pm$ 1.97  &  26.83 $\pm$ 1.11   & 69.08  $\pm$ 1.29  & 26.63 $\pm$  1.09 &
 67.74 $\pm$ 1.40
 \\ \cline{2-8}
Unimodal &   0.9   &  23.79 $\pm$ 1.08  &  68.16  $\pm$ 3.92 &  25.88 $\pm$  1.18  & 58.04 $\pm$  2.17 &  25.63 $\pm$ 1.21 &
 55.39 $\pm$ 2.34
 \\ \cline{2-8}
 & 1.0    &  23.07 $\pm$ 1.12 & 66.39 $\pm$ 4.30 & 25.71 $\pm$ 1.42 & 53.55 $\pm$ 2.60 & 25.23 $\pm$ 1.24 & 50.28 $\pm$ 2.54 \\ \hline\hline
 & 0.6    &   11.69 $\pm$ 0.30 & 64.49 $\pm$ 0.86 & 12.90 $\pm$ 0.51 & 50.22 $\pm$ 2.03  & 12.82 $\pm$ 0.49 &
49.03 $\pm$ 2.03 \\ \cline{2-8}
Bimodal &   0.9   & $\mbox{ \ }$6.79 $\pm$ 0.21 & 49.85 $\pm$ 3.39 & $\mbox{ \ }$9.08 $\pm$ 0.45 & 27.52 $\pm$ 2.73 & $\mbox{ \ }$8.92 $\pm $ 0.64 & 26.69 $\pm$ 2.67 \\ \cline{2-8}
 & 1.0  &   $\mbox{ \ }$2.33 $\pm$ 0.20 & 21.40 $\pm$ 0.76 & $\mbox{ \ }$5.48 $\pm$ 1.73 & 18.48 $\pm$ 2.45 & $\mbox{ \ }$5.52 $\pm$ 1.75 &
 18.22 $\pm$ 2.47 \\ \hline\hline
&   0.6   & 15.72 $\pm$ 0.46 & 66.48 $\pm$ 2.63 & 16.47 $\pm$ 0.67 & 61.82 $\pm$ 1.96 & 16.34 $\pm$ 0.78 & 60.45 $\pm$ 2.24 \\ \cline{2-8}
Trimodal &   0.9   & 11.67 $\pm$ 0.37 & 54.15 $\pm$ 1.74 & 14.17 $\pm$ 0.93 & 44.68 $\pm$ 3.36 & 14.05 $\pm$ 0.98 & 42.72 $\pm$ 3.60 \\ \cline{2-8}
 & 1.0    & $\mbox{ \ }$8.46 $\pm$ 0.36  & 48.32 $\pm$ 2.53 & 13.33 $\pm$ 1.01 & 37.88 $\pm$ 4.01 & 13.32 $\pm$ 0.98 & 35.37 $\pm$ 4.15 \\ \hline
\end{tabular}
\end{center}
%\vspace{0.2cm}
\label{table_costs_2}
\end{table}

\subsubsection{Concluding Remarks}

%The efficient sampling techniques developed in this paper allow GOSPA-based control approaches to be used to performance non-myopic sensor management even when a large number of potential actions are considered.

It is noted that multi-step planning often generates the same first action as myopic planning  when there is a low $P_d$. This is because there is no guarantee that the target will be detected across multiple time steps, and so regions of high probability mass offer the greatest opportunity achieve at least one target detection. By design, myopic planning always favours these high probability mass regions because its decision-making lacks the foresight to appreciate that further observations are possible.

Conversely, in scenarios with a high $P_d$, offsetting the first observation and using a cookie-cutter strategy often allows the multi-step approach to decrease the number of time steps required to provide complete surveillance of the region of interest. This approach should therefore be favoured if time constraints allow.

\section{Summary and Conclusions}
\label{sec:conc}

%This paper is concerned with sensor management for target search and  track applications. The generalised optimal subpattern assignment (GOSPA) metric is the performance metric of choice, and provides a unified mechanism for combining together the costs corresponding to localisation errors for properly detected targets, and cardinality errors for missed and false targets. These represent the primary metrics of interest in multiple target estimation.

%However, utilising the GOSPA metric to predict system performance in performing  sensor management is computationally challenging, because of the need to account for (and average across) uncertainties within the scenario, notably the number of targets, the locations of targets (that exist), and the measurements generated by the targets as a result of performing sensing actions. Consequently, to-date, GOSPA-based sensor management has been  demonstrated only in  simple scenarios with just a single action and a single hypothesis representing each  target location.

This paper has proposed a sample-based approach for myopic and non-myopic sensor management for a Bernoulli target using the GOSPA metric. We have provided the following contributions: analytical calculation of the MS-GOSPA for different measurements and actions, development of efficient sampling techniques to calculate the  AMMS-GOSPA error, and the development of an optimal non-myopic (Bellman type \cite{Bellman_1952}) planning recursion that exploits the conditional AMMS-GOSPA error.
Simulations demonstrate the approach in scenarios with: (i): missed detections, (ii): false alarms,
(iii): a high degree of uncertainty in the target location, with the prior distribution represented by large
    number of potential hypotheses, and: (iv): a planning horizon of up to three time
steps.

Various behavioural patterns are identified, notably
demonstrating the benefits of non-myopic planning, and in particular showing that optimal plans align with an
intuitive understanding of how taking into account the opportunity to make further observations should influence
the current action.
It is concluded that the GOSPA-based, non-myopic search and track algorithm offers a powerful mechanism for sensor
management in order to minimise estimation errors and errors due to missed and false targets in a unified way.

The current approach is directly applicable to multi-target scenarios with well-separated targets, due to the separability of the optimal actions when using the GOSPA metric \cite{garcia_fernandez_2021}. Future work will extend the approach to multi-sensor, multi-target scenarios in which targets may move in close proximity. We will also implement Monte Carlo roll-out as a mechanism for efficiently estimating the long-term impact of actions.
Furthermore, we will work to identify scenarios in which the suboptimal and optimal approaches generate  markedly different solutions, thereby highlighting the importance of accounting for the potential sequences of future measurements at each decision epoch.

\section*{Acknowledgements}

This research was funded  by the UK Ministry of Defence (MOD) through the Fusion and Information Theory (FIT) project, DSTLX-1000143908.
 The authors would like to thank Professor Gustav Hendeby (Link\"{o}ping University) and Professor Rahul Savani (University of Liverpool) for insightful discussions regarding sensor management.

\appendix

In this appendix, we prove the following proposition.

\vspace{0.25cm}
\noindent
\textit{Proposition}

\vspace{0.25cm}
\noindent
 \textit{Using the ideal measurement set assumptions: (i): $\Sigma=0$ (i.e. target generated measurements are error-free), and: (ii): $\lambda_{FA}=0$ (i.e. there are no false alarms), but allowing for missed detections (i.e. allowing $P_d< 1$), the suboptimal and optimal control approaches generate identical solutions (i.e. $\hat{a}^\star_1 = {a}^\star_1$) irrespective of the length of the planning horizon.}

\vspace{0.25cm}
\noindent
\textit{Proof}

\vspace{0.25cm}
\noindent
In this case,  MMS-GOSPA$(z_{1:t}, a_{1:t}) = 0$ unless $z_i=\phi$ for $i=1, \ldots, t$. This is because, with extremely accurate measurements, and no clutter, the presence of just a single measurement signifies that a target is present, and the accurate measurement allows the target to be geo-located without error. Hence, in determining the GOSPA-based cost function, it is necessary to only consider cases in which there are no previous measurements.

The optimal action (\ref{optimal_Bellman_act}) can then be
manipulated as follows:
\begin{eqnarray}
 \hat{a}^\star_1 & = & \argmin_{a_{1}}\bigg[p(z_1=\phi|a_1)\Big[r_1(z_{1}=\phi, a_{1}) +
\lambda\min_{a_{2}}p(z_2=\phi|z_1=\phi,a_{1:2})\big[r_{2}(z_{1:2}=\phi, a_{1:2}) \label{optimal_Bellman_act_app1}  \\
%&&  + \lambda\min_{a_{3}}p(z_3=\phi|z_{1:2}%=\phi,a_{1:3})\Big[ %r_{3}(z_{1:3}=\phi,a_{1:3})
&& \qquad + \ldots +
\lambda\min_{a_T}p(z_T=\phi|z_{1:T-1}=\phi,a_{1:T})\big[r_{T}(z_{1:T}=\phi, a_{1:T})\big]\Big]\bigg] \nonumber \\
&  = & \argmin_{a_{1}}\bigg[{\mathbb E}_{z_1}\left[r_1(z_{1}, a_{1})\right] +
\lambda\min_{a_{2}}\Big[ {\mathbb E}_{z_{1:2}}\left[r_{2}(z_{1:2}, a_{1:2})\right] \label{optimal_Bellman_act_app2}  + \ldots +
\lambda\min_{a_T}\big[{\mathbb E}_{z_{1:T}}\left[r_{T}(z_{1:T}, a_{1:T})\right]\big]\Big]\bigg] \qquad \\
& = &\argmin_{a_{1}}\Bigg[\min_{a_2:T}
\bigg[{\mathbb E}_{z_{1:T}}\Big[
\sum_{t=1}^T \lambda^{t-1}r_t(z_{1:t},a_{1:t})
\Big]\bigg]
\Bigg] \ \mbox{(as the conditional expectations have been removed)} \qquad \\
& = & {a}^\star_1 \qquad \mbox{(given by equation (\ref{optimal_action_multi}))}
\label{optimal_Bellman_act_app4}
\end{eqnarray}

\noindent
This completes the proof.
\vspace{-0.6cm}
\begin{flushright}
$\square$
\end{flushright}

\end{document}